\documentclass[10pt,a4paper]{article}

\usepackage{amssymb, amsmath, amsthm, amssymb, mathtools,geometry}
\usepackage[expansion=false]{microtype}

\usepackage[cm]{fullpage}
\usepackage[english]{babel}
\usepackage[pdftex]{graphicx}
\usepackage{booktabs}
\usepackage{xcolor,enumerate}
\usepackage{todonotes}
\usepackage{cancel}
\usepackage[normalem]{ulem}

\usepackage{hyperref}
\usepackage{booktabs}
\usepackage{nameref,cleveref}

\numberwithin{equation}{section}

\newtheorem{theorem}{Theorem}[section]
\newtheorem{lemma}[theorem]{Lemma}
\newtheorem{proposition}[theorem]{Proposition}
\newtheorem{remark}[theorem]{Remark}

\newtheorem{definition}[theorem]{Definition}

\newcommand{\veps}{\varepsilon}

\newcommand{\R}{\mathbb{R}}
\newcommand{\N}{\mathbb{N}}
\newcommand{\Z}{\mathbb{Z}}
\newcommand{\dd}{\,\textup{d}}

\newcommand{\fl}{(-\Delta)^s}
\newcommand{\fld}{(-\Delta_h)^s}
\newcommand{\coma}{\, , \,}
\newcommand{\Q}{\mathcal{Q}}
\newcommand{\X}{\mathcal{X}}

\usepackage{tocloft}
\title{\bf The obstacle problem for \\ fractional heat equations with memory
}
\author{F\'elix del Teso\thanks{Departamento de Matem\'aticas, Universidad Aut\'onoma de Madrid, Spain, felix.delteso@uam.es},
 \; Fernando Quir\'os\thanks{Departamento de Matem\'aticas, Universidad Aut\'onoma de Madrid, Spain, fernando.quiros@uam.es},
 \; Julio D. Rossi\thanks{Departamento de Matem\'atica y Estad\'istica, Universidad Torcuato Di Tella, Buenos Aires, Argentina, julio.rossi@utdt.edu}
 }
\date{}

\begin{document}
\maketitle

\begin{abstract}
We study the obstacle problem for the evolution equation $\partial_t^\alpha u+\fl u=0$ in the whole space~$\R^d$, where $\partial_t^\alpha$ with $\alpha\in(0,1)$ is the Caputo derivative in time and $\fl$ with $s\in(0,1]$ is the fractional Laplacian in space ($s<1$) or the classical Laplacian ($s=1$). We prove that the problem has at most one bounded viscosity solution, and we obtain existence as the limit of an explicit monotone finite-difference scheme, which combines the L1 discretization of the Caputo derivative with powers of the discrete Laplacian. The scheme is also used to establish regularity estimates in space and time. We then study the large-time behaviour: solutions converge, along subsequences, to stationary solutions, and we identify the limit when the initial datum coincides with the obstacle (the minimal stationary solution) or has a limit $c_*$ at infinity. In the latter case the dimension plays a role, as a consequence of the recurrence of the underlying process: the limit is the unique stationary solution with limit $c_*$ at infinity if $d>2s$, and a constant (the maximum between $c_*$ and the supremum of the obstacle) if $d\le2s$. Numerical experiments in one dimension, based on an explicit stationary solution, illustrate the convergence of the scheme and the relaxation towards the stationary profile.
\end{abstract}



\begingroup
\footnotesize
\tableofcontents
\endgroup

\section{Introduction}\label{sec:intro}

The parabolic obstacle problem consists in finding the smallest supersolution of a given parabolic equation that lies above a prescribed function, the obstacle. It arises in several applications, such as phase transitions, where the one-phase Stefan problem can be recast, via the Duvaut transformation, as a parabolic obstacle problem~\cite{Duvaut-1973,Friedman-Kinderlehrer-1975}, and mathematical finance, where the price of an American option is the solution of an obstacle problem whose obstacle is the payoff function, a characterisation that goes back to McKean~\cite{McKean-1965} (see also~\cite{Bensoussan-Lions-1978,Jaillet-Lamberton-Lapeyre-1990}). In this paper we study the obstacle problem for an evolution equation that is nonlocal both in space and in time.

\noindent\textsc{The model. } We consider the Cauchy problem
\begin{align}
    \min\{\partial_t^\alpha u(x,t)+\fl u(x,t) \coma u(x,t) - \varphi(x)\} &=0&&\textup{for } (x,t)\in \R^d\times(0,\infty),\label{eq:main}\\
    u(x,0)&= u_0(x) &&\textup{for } x\in \R^d,\label{eq:initial}
\end{align}
where $\alpha\in(0,1)$, $s\in(0,1]$, with an obstacle $\varphi$ and an initial datum $u_0$ belonging to $C_\textup{b}(\R^d)$, and satisfying $u_0\geq\varphi$ in $\R^d$.

The operator $\partial_t^\alpha$ is the Caputo derivative of order $\alpha\in(0,1)$, defined for smooth functions $\phi:\R_+\to\R$~by
\begin{equation}\label{eq:capdef}
    \partial_t^\alpha \phi(t)=\frac{1}{\Gamma(1-\alpha)}\frac{\phi(t)-\phi(0)}{t^{\alpha}}+\frac{\alpha}{\Gamma(1-\alpha)}\int_0^t \frac{\phi(t)-\phi(s)}{(t-s)^{1+\alpha}} \dd s,
\end{equation}
or, equivalently, integrating by parts,
\begin{equation}\label{eq:capdef2}
    \partial_t^\alpha \phi(t)=\frac{1}{\Gamma(1-\alpha)}\int_0^t\frac{\phi'(s)}{(t-s)^{\alpha}} \dd s.
\end{equation}
Its value at time $t$ depends on the whole history of $\phi$ on $[0,t]$, with a weight decaying algebraically; in this sense the Caputo derivative encodes a memory effect, in sharp contrast with the classical, purely local, time derivative.

The operator $\fl$ is the fractional Laplacian of order $s\in(0,1]$, defined for smooth functions $\psi: \R^d \to \R$ by
\begin{align}\label{eq:definition.fractional.laplacian}
    \fl \psi(x)= C_{d,s} \,\textup{P.V.}\int_{\R^d}\frac{\psi(x)-\psi(x+y)}{|y|^{d+2s}}\dd y, \quad \mbox{with }
      C_{d,s}=\frac{4^{s} \Gamma\left(s+\frac{d}{2}\right)}{\pi^{\frac{d}{2}}\Gamma(2-s)}s(1-s),
\end{align}
for $s\in(0,1)$, and as the classical Laplacian $-\Delta \psi(x) =- \sum_{i=1}^d \partial^2_{x_i} \psi(x)$ for $s=1$. For $s\in(0,1)$ the value of $\fl\psi(x)$ depends on $\psi$ in the whole space, not just near $x$; the equation is thus nonlocal in space as well.

\noindent\textsc{Motivation. } Both operators in~\eqref{eq:main} have been proposed in the mathematical finance literature as tools to capture features of asset price dynamics that classical models fail to reproduce. The operator $-\fl$ with $s\in(0,1)$ is the generator of a symmetric $2s$-stable Lévy process, whose heavy-tailed increments model sudden jumps in asset prices. Models with jumps for option pricing go back to Merton~\cite{Merton-1976}; in the local-in-time setting $\alpha=1$, the price of an American option in such models solves an obstacle problem for an integro-differential operator~\cite{Pham-1997} (see~\cite{Chen} for a model involving a fractional operator), which, for $\fl$, has been thoroughly analyzed from the PDE viewpoint; see~\cite{Caffarelli-Figalli-2013,Barrios-Figalli-RosOton-2018} and the references therein.

The Caputo derivative, on the other hand, slows down the dynamics: for $s=1$ the mean squared displacement of the underlying process grows like $t^\alpha$ rather than linearly in time. In the presence of jumps the two effects compete, and the natural scaling $|x|\sim t^{\alpha/(2s)}$ is subdiffusive if $\alpha<s$ and superdiffusive if $\alpha>s$. In the Black--Scholes framework, time-fractional derivatives were introduced in~\cite{Wyss-2000}. They arise naturally in continuous-time random walk models of tick-by-tick price dynamics, in which the Caputo derivative accounts for heavy-tailed random waiting times between consecutive transactions; remarkably, the scaling limit of these models is a diffusion equation that is fractional both in time and in space, precisely with the operators in~\eqref{eq:main}; see~\cite{Scalas-Gorenflo-Mainardi-2000,Meerschaert-Scalas-2006}. The combination of both operators in a single model is thus natural: it captures simultaneously jumps in space and memory in time.

A rigorous financial justification of the obstacle problem~\eqref{eq:main} lies beyond the scope of this paper, which is a work of mathematical analysis; we nevertheless believe that the theory developed here (well-posedness, regularity, large-time behaviour and rigorous numerical analysis) may be of interest in that context.

\noindent\textsc{Notation. }
Throughout the paper, $B(\mathbb{R}^d)$ denotes the space of
bounded (not necessarily continuous) functions on $\mathbb{R}^d$,
$C_{\mathrm{b}}(\mathbb{R}^d)$ the space of bounded continuous
functions, and $C_{\mathrm{b}}^0(\mathbb{R}^d)$ the space of bounded
\emph{uniformly} continuous functions, equipped with the
$L^\infty(\R^d)$ norm. For each integer $k\geq 1$,
$C_{\mathrm{b}}^k(\mathbb{R}^d)$ denotes the space of functions in
$C^k(\mathbb{R}^d)$ whose derivatives of order at most $k$,
including order zero, are bounded and uniformly continuous. For an open set $\Omega\subset\mathbb{R}^n$,
$C_{\mathrm{b}}^\infty(\Omega)$ denotes the space of smooth
functions on $\Omega$ with bounded derivatives of all orders,
including order zero, and $C_{\mathrm{c}}^\infty(\Omega)$ the
space of smooth functions with compact support in $\Omega$.
We also write $C_{\mathrm{b}}^\infty(\overline{\Omega})$ for the
functions in $C_{\mathrm{b}}^\infty(\Omega)$ whose derivatives
of every order extend continuously to $\overline{\Omega}$. For $a=0$ we set $\mathcal{X}_{\mathrm{b}}^0(\mathbb{R}^d):=C_{\mathrm{b}}^0(\mathbb{R}^d)$, equipped with the $L^\infty(\R^d)$ norm. For $a=k+\beta\in(0,2]$, with $k\in\{0,1\}$ and
$\beta\in(0,1]$, we set $D^0f:=f$, $D^1f:=\nabla f$, and define the space
\begin{equation*}
    \mathcal{X}_{\mathrm{b}}^a(\mathbb{R}^d)
    := \left\{
        f \in C_{\mathrm{b}}^k(\mathbb{R}^d)
        : [D^k f]_\beta < \infty
    \right\}
    \quad \text{with} \quad
    [g]_\beta
    := \sup_{x\neq y}
       \frac{|g(x)-g(y)|}{|x-y|^\beta}.
\end{equation*}
equipped with the norm
\begin{equation*}
    \|f\|_{\mathcal{X}^a(\mathbb{R}^d)}
    := \sum_{j=0}^{k}
       \|D^j f\|_{L^\infty(\mathbb{R}^d)}
       + [D^k f]_\beta.
\end{equation*}

 We write $\N$ for the set of nonnegative integers, and $p\wedge q\coloneqq\min\{p,q\}$, $p\vee q\coloneqq\max\{p,q\}$. Finally,~$C$ denotes a positive constant whose value may change from line to line, and whose dependence on the parameters is indicated by subscripts when relevant.

\noindent\textsc{Background and related work. } The mathematical theory of parabolic obstacle problems with local operators has a long history. An existence theory for parabolic unilateral obstacle problems was developed by J.-L.~Lions and G.~Stampacchia
in~\cite[Section~7]{Lions-Stampacchia-1967}; regularity of solutions was subsequently studied by H.~Brézis~\cite[Chapter~II]{Brezis-1972}, and extended to unbounded domains, with applications to optimal stopping problems, by A.~Friedman~\cite{Friedman-1973b}. The connection between obstacle problems and optimal stopping goes back at least to McKean~\cite{McKean-1965}; see also Bensoussan and Lions~\cite{Bensoussan-Lions-1973,Bensoussan-Lions-1978}.

Part of our analysis relies on the theory of the time-fractional heat equation $\partial_t^\alpha u+\fl u=0$ \emph{without} obstacle. 
 Solutions to the corresponding Cauchy problem are usually obtained as mild (or integral) solutions, given by the convolution of the initial datum with the fundamental solution. The existence of the latter and sharp estimates for it go back to Kochubei~\cite{Kochubei-1990} and Eidelman and Kochubei~\cite{Eidelman-Kochubei-2004} when $s=1$, and were obtained by Kemppainen, Siljander and Zacher~\cite{Kemppainen-Siljander-Zacher-2017} for the fully nonlocal equation, $s\in(0,1)$. Since the obstacle problem~\eqref{eq:main} is fully nonlinear, mild solutions are not available for it, and we  choose to work instead with viscosity solutions.

Viscosity solutions for equations with Caputo time derivative were used by Allen~\cite{Allen-2018} in the study of the regularity of a nonlocal nondivergence parabolic equation. A well-posedness theory for fully nonlinear second-order equations was developed by Namba~\cite{Namba-2018}, in bounded domains with Dirichlet or Neumann boundary conditions, and by Topp and Yangari~\cite{Erwin}, in the whole space with nonlocal spatial operators of Lévy type. The latter obtain existence by Perron's method, comparison for bounded sub- and supersolutions through a regularization by sup/inf-convolutions in time, and the large-time behaviour. However, the obstacle problem~\eqref{eq:main} does not fit into the class of equations considered in~\cite{Namba-2018,Erwin}, since there the Caputo derivative appears outside the nonlinearity, whereas in~\eqref{eq:main} it appears inside the minimum. However, in the non-coincidence set $\{u>\varphi\}$ the problem reduces to the linear time-fractional heat equation, which does fall within their framework, and our uniqueness proof in Section~\ref{sec:viscosity} rests on an adaptation of their comparison argument that exploits this fact. Mild solutions of the time-fractional heat equation with bounded continuous data are also viscosity solutions, which allows us to combine both theories in Section~\ref{sec:asymptotics}.

On the numerical side, a monotone time-discrete scheme for fully nonlinear second-order equations with Caputo derivative, converging to the viscosity solution, was introduced in~\cite{Giga-Liu-Mitake-2020}; fully discrete monotone
finite-difference schemes for equations that are nonlocal in space and time, together with regularity estimates for the limit problem, were developed in~\cite{delTeso-Lukasz-2025}. The spatial discretization we use, based on powers of the discrete Laplacian, goes back to~\cite{Ciaurri2018}, see also~\cite{delTeso2019a,delTeso2018b}.

 Concerning obstacle problems with a Caputo time derivative, a variational theory has been developed by Li and Salgado~\cite{Li-Salgado-2023}, who study time-fractional gradient flows of convex, lower semicontinuous energies in Hilbert spaces. They prove existence, uniqueness and H\"older regularity in time of energy solutions, and derive error estimates for implicit time-stepping schemes. The time-fractional obstacle problem in a bounded domain, with the obstacle constraint encoded in the energy through the indicator function of the admissible set, is one of their examples, and its space-fractional version is considered in their numerical experiments. In this variational framework, data and solutions belong to $L^2$, and the large-time analysis relies on the convexity of the energy: for uniformly convex energies, solutions converge in $L^2$ to the minimizer at a Mittag-Leffler rate, while for merely convex energies only the decay of the energy is quantified. Our setting is different: we work in the whole space, with bounded data and viscosity solutions, where the elliptic obstacle problem has infinitely many solutions and the large-time limit has to be identified by other means. In addition, we obtain $L^\infty$ regularity estimates in space and time. Optimal control problems in which the obstacle is the control have also been considered for the time-fractional obstacle problem with the Laplacian in bounded domains, see~\cite{SaNgiamsunthorn-Suechoei-Kumam-2021}.

Another reference on obstacle problems with a Caputo time derivative is~\cite{Alberini-et-al-2022}, which considers the equation $\partial_t^\alpha u-H(u-\psi)\Delta u=0$,  where $H$ denotes the Heaviside function, a model different from ours. Their analysis for $\alpha\in(0,1)$ relies on an equivalence with a complementarity system.
Their argument uses properties of the classical problem that need not hold in the Caputo setting. In particular, a vanishing time derivative at a contact point need not imply a vanishing Caputo derivative at that point. Thus, this argument leaves the well-posedness of their model for \(\alpha\in(0,1)\) unresolved.

\noindent\textsc{Statements of the main results and general strategies for the proofs.} Since the equation is fully nonlinear, viscosity solutions provide a natural framework, and uniqueness follows from a comparison principle in the non-coincidence set, adapted from~\cite{Erwin}. Existence is more delicate: the memory encoded in the Caputo derivative destroys the semigroup property, so that semigroup techniques do not apply. Energy methods are available for equations with Caputo derivative in divergence form; see e.g.~\cite{Allen-Caffarelli-Vasseur-2016}, and for the obstacle problem in bounded domains, see~\cite{Li-Salgado-2023}. These methods are not well suited to the viscosity framework in the whole space adopted here. Existence could  presumably be obtained by adapting Perron's method, as used in~\cite{Namba-2018,Erwin} for equations with the Caputo derivative outside the nonlinearity, see \Cref{rem:Perron}; we follow instead a constructive approach, obtaining the solution as the limit of an explicit, monotone finite-difference scheme in space and time, in the spirit of~\cite{delTeso-Lukasz-2025}, which we describe after stating our results.  Besides existence, this approach yields a convergent numerical method and regularity estimates, and the latter turn out to be the key to the large-time behaviour. Indeed, the regularity in time inherited from the scheme shows that the Caputo derivative of the solution is well defined and decays as $t\to\infty$, a consequence of the memory effect. Moreover, passing to the limit in the scheme we show that the equation holds classically in time, with the Caputo derivative of the solution itself appearing in it, which, combined with the decay, allows us to identify the large-time limits as solutions of the stationary obstacle problem. This argument relies essentially on the memory effect and does not extend to the local case $\alpha=1$.

Let us now state our main results precisely. The first one is uniqueness of viscosity solutions of the time-fractional obstacle problem.

\begin{theorem}[Uniqueness]\label{thm:intro-uniqueness}
    Let $s\in(0,1]$, $\alpha\in(0,1)$, and $u_0,\varphi\in C_{\textup{b}}(\R^d)$ be such that $u_0\geq\varphi$ in $\R^d$. Then problem~\eqref{eq:main}--\eqref{eq:initial} has at most one viscosity solution. Moreover, if a viscosity solution $u$
    exists, it is given by
    \[
        u(x,t) = \inf_{v\in\mathcal{V}} v(x,t),
    \]
    where $\mathcal{V}=\{v\in C_{\textup{b}}(\R^d\times[0,\infty)): v \textup{ is a viscosity supersolution to~\eqref{eq:main}--\eqref{eq:initial}}\}$.
\end{theorem}

Our second result establishes existence and collects the main properties of the solution.

\begin{theorem}[Existence and properties]\label{thm:intro-existence}
    Let $s\in(0,1]$ and $\alpha\in(0,1)$, and let $u_0,\varphi\in C_{\textup{b}}^{0}(\R^d)$ be such that $u_0\geq\varphi$ in $\R^d$. Then problem~\eqref{eq:main}--\eqref{eq:initial} has a unique viscosity solution $u\in C_{\textup{b}}^{0}(\R^d\times[0,\infty))$. Moreover, the following properties hold:
    \begin{enumerate}[\rm (a)]
        \item \textup{(Stability)} $\varphi(x)\leq u(x,t)\leq\sup_{\R^d}u_0$ for all $(x,t)\in\R^d\times[0,\infty)$.

        \item \textup{($L^\infty$-contraction)} Let $v$ be the viscosity solution of problem~\eqref{eq:main}--\eqref{eq:initial} with initial datum $v_0\in C_{\textup{b}}^0(\R^d)$ and obstacle $\psi\in C_{\textup{b}}^0(\R^d)$ satisfying $v_0\geq\psi$ in $\R^d$. Then, for all $t\geq0$,
        \[
            \|u(\cdot,t)-v(\cdot,t)\|_{L^\infty(\R^d)}\leq\max\left\{\|u_0-v_0\|_{L^\infty(\R^d)},\|\varphi-\psi\|_{L^\infty(\R^d)}\right\}.
        \]

        \item \textup{(Comparison principle)} If, additionally, $u_0\leq v_0$ and $\varphi\leq \psi$ in $\R^d$, then $u(x,t)\leq v(x,t)$ for all $(x,t)\in \R^d\times[0,\infty)$.

        \item \textup{(Uniform continuity)} $u$ is uniformly continuous in space and time, with a modulus of continuity depending only on $u_0$ and $\varphi$.

        \item \textup{(Continuity estimates)}  If, in addition, $u_0\in \X_{\textup{b}}^{a}(\R^d)$ and $\varphi\in \X_{\textup{b}}^{b}(\R^d)$
        for some $a,b\in(0,2s+\veps]$, where $\veps=0$ if $s=1$ and $\veps\in(0,2-2s]$ if $s\in(0,1)$, then the following regularity
        estimates hold:
        \begin{enumerate}[\rm (i)]
            \item \textup{(Space)} For all $t\geq0$ and $y\in\R^d$,
            \[
                \|u(\cdot+y,t)-u(\cdot,t)\|_{L^\infty(\R^d)}\leq\max\left\{\|u_0\|_{\X^a(\R^d)}|y|^{\min\{a,1\}},\|\varphi\|_{\X^b(\R^d)}|y|^{\min\{b,1\}}\right\}.
            \]
            \item \textup{(Time)} For all $t,\tilde t\in[0,\infty)$,
            \begin{align*}
                \|u(\cdot,t)-u(\cdot,\tilde{t})\|_{L^\infty(\R^d)}&\leq C_{d,s,a,b,\veps,\alpha}\max\{\|u_0\|_{\X^a(\R^d)}\coma \|\varphi\|_{\X^b(\R^d)}\}\\
                &\quad\times\max\left\{|t^\alpha-\tilde{t}^\alpha|^{\frac{a}{2s+\veps}}\coma|t^\alpha-\tilde{t}^\alpha|^{\frac{b}{2s+\veps+b-a}}\coma |t^\alpha-\tilde{t}^\alpha|\right\}.
            \end{align*}
           In particular, if $a=b=2s+\veps$, then $u$ is $C^\alpha$ in time, uniformly in space.
        \end{enumerate}
    \end{enumerate}
\end{theorem}

The numerical approximation used in the proof of \Cref{thm:intro-existence} is also of independent interest. We consider an explicit scheme based on the L1 discretization of the Caputo derivative in time and on the powers of the discrete Laplacian approximation of the fractional Laplacian in space; see~\Cref{subsec:scheme} for the details. Under a suitable CFL-type stability condition, the scheme is monotone and stable, satisfies comparison and $L^\infty$-contraction properties, and enjoys regularity estimates in space and time that are uniform in the discretization parameters. These estimates yield compactness of the discrete solutions and, together with the monotonicity and consistency of the scheme, imply their locally uniform convergence to the unique viscosity solution of~\eqref{eq:main}--\eqref{eq:initial}. In particular, this provides the existence part of \Cref{thm:intro-existence}, and the remaining properties stated there are inherited from the corresponding properties of the scheme (see Section~\ref{sec:scheme}).

We also analyze the large-time behaviour of the discrete evolution problem. We show that, along subsequences, the solution of the scheme converges locally uniformly as $n\to\infty$, and that every such limit is a solution of the corresponding discrete elliptic obstacle problem. When the initial datum coincides with the obstacle, the whole sequence converges locally uniformly to the minimal solution of that problem (see Section~\ref{sec:discrete-asymptotics}). This shows that the explicit scheme can be used as a solver for the discrete elliptic obstacle problem. Moreover, as $h\to0$, these limits converge, along subsequences, to viscosity solutions of the continuous elliptic obstacle problem (see Section~\ref{subsec:discrete-approx}).

The stationary problem associated to~\eqref{eq:main}--\eqref{eq:initial} is the elliptic obstacle problem
\begin{equation}\label{eq:intro-elliptic}
    \min\{\fl u(x)\coma u(x)-\varphi(x)\}=0,\quad x\in\R^d.
\end{equation}
Unlike the evolution problem~\eqref{eq:main}--\eqref{eq:initial}, problem~\eqref{eq:intro-elliptic} does not have a unique solution in general: any constant function lying above $\varphi$ is a solution. As a consequence, the large-time limit of $u(\cdot,t)$ may a priori depend on the subsequence of times considered. Nevertheless, we are able to identify the limit in two important cases.

The first is when the initial datum coincides with the obstacle; in this case the solution converges to the minimal viscosity solution of the stationary problem~\eqref{eq:intro-elliptic}, defined as the infimum of all viscosity supersolutions. This is precisely the  situation in the financial application mentioned above: the price of an  American option with payoff $\varphi$ is determined backwards from the maturity, where it equals the payoff, so that in the ``time to maturity'' variable the initial datum is the obstacle itself. The large-time limit  then corresponds to the \emph{perpetual} option, and it is classical in optimal stopping theory that the value of a perpetual option is the smallest excessive (that is, superharmonic) majorant of the payoff (see, e.g.,~\cite{Peskir-Shiryaev-2006}, and~\cite{Pham-1997} for jump-diffusion models); in our language, the minimal solution of the stationary obstacle problem. \Cref{thm:intro-asymptotics}\eqref{thm:intro-asymptotics-item2} below is the counterpart of this fact for the time-fractional model. We stress
that the correspondence is only formal: the operators arising in option pricing, both in the Black--Scholes and in the Lévy models, include drift and discounting terms that are absent in~\eqref{eq:main}, so that our model corresponds to a simplified setting (zero interest rate and a driftless underlying, in logarithmic variables), and the interpretation of  the time-fractional evolution in the time-to-maturity variable is, as already discussed, heuristic. Note also that typical payoffs, such as $(K-S)^+$ for a put option, are
Lipschitz continuous but not $C^1$, so that the Hölder framework of \Cref{thm:intro-asymptotics}, which requires $u_0,\varphi\in \X^a_\textup{b}(\R^d)$ for some $a\in(0,2]$ only, is the appropriate one for this application.

The second situation in which we can identify the large-time behaviour is when the initial datum $u_0$ has a prescribed limit~$c_*$ at infinity; in this case, and provided $d>2s$, the limit is the unique solution of~\eqref{eq:intro-elliptic} with the same behaviour at infinity. When $d\le 2s$, on the contrary, every bounded $s$-superharmonic function is constant, so that the solutions of~\eqref{eq:intro-elliptic} are exactly the constants above $\sup\varphi$, and the large-time limit turns out to be $\max\{\sup\varphi,c_*\}$.

Initial data with a limit at infinity have a natural interpretation in the optimal stopping setting discussed above. Bounded payoffs are common in practice, as is the case for capped options, whose payoff is cut off at a prescribed level~\cite{Broadie-Detemple-1995,Broadie-Detemple-1999}; those whose payoff vanishes when the underlying is very small or very large, such as butterfly spreads (see, e.g., \cite[Section~3.3]{Wilmott-Howison-Dewynne-1995} for payoff diagrams of option combinations), give rise to initial data with limit $c_*=0$ at infinity, in the sense of \Cref{thm:intro-asymptotics}. We stress that this is only an interpretation, since the financial meaning of the time-fractional model is heuristic. 

We emphasize that boundedness of the initial datum alone, without any assumption on its behaviour at infinity, is insufficient to characterize the long-time behaviour, even in the linear setting. Indeed, as shown in \cite{VazquezZuazua2002}, complex asymptotic patterns may arise for bounded initial data. Some control at infinity is therefore necessary. In this regard, requiring the initial datum to converge to a constant $c_*$ at infinity provides a natural extension of the framework of integrable initial data previously studied for the linear equation in \cite{Kemppainen-Siljander-Zacher-2017,CQW-JFA,CQW-VJM}.

We can now state our results on the large-time behaviour. Besides the general convergence, up to subsequences, to solutions of the stationary problem, we identify the limit in the two situations described above.

\begin{theorem}[Large-time behaviour]\label{thm:intro-asymptotics}
    Let $s\in(0,1]$, $\alpha\in(0,1)$, and $u_0,\varphi\in \X_{\textup{b}}^a(\R^d)$ for some $a\in(0,2]$ be such that $u_0\geq\varphi$ in $\R^d$. Let $u$ be the viscosity solution of~\eqref{eq:main}--\eqref{eq:initial}. Then:
    \begin{enumerate}[\rm (a)]
        \item There exists a sequence $t_k\to+\infty$ such that $u(\cdot,t_k)\to v$ locally uniformly in $\R^d$ as $k\to+\infty$,
        where $v$ is a viscosity solution of the stationary problem~\eqref{eq:intro-elliptic}.

        \item\label{thm:intro-asymptotics-item2} If $u_0=\varphi$ in $\R^d$, then $u(\cdot,t)\to\underline{v}$ locally uniformly in $\R^d$ as $t\to+\infty$, where $\underline{v}$ is the minimal viscosity solution of~\eqref{eq:intro-elliptic}.

        \item If the limit $c_*=\lim_{|x|\to\infty}u_0(x)$ exists, then $u(\cdot,t)\to\underline{v}^{c_*}$ locally uniformly in $\R^d$ as $t\to+\infty$, where $\underline{v}^{c_*}$ is the minimal viscosity solution of~\eqref{eq:intro-elliptic} with obstacle $\max\{\varphi,c_*\}$. Moreover, 
        \begin{enumerate}[\rm (i)]
            \item if $d>2s$, then $\underline{v}^{c_*}$ is the unique viscosity solution of~\eqref{eq:intro-elliptic} satisfying $\lim_{|x|\to\infty}\underline{v}^{c_*}(x)=c_*$, and the convergence $u(\cdot,t)\to\underline{v}^{c_*}$ is uniform in $\R^d$;
            \item if $d\le2s$, then $\underline{v}^{c_*}\equiv\max\{\sup\varphi,c_*\}$ and the convergence $u(\cdot,t)\to\underline{v}^{c_*}$ is not uniform in $\R^d$ when $c_*<\sup\varphi$, since $\lim_{|x|\to\infty}u(x,t)=c_*$ for every $t>0$.
        \end{enumerate}
    \end{enumerate}
\end{theorem}

The numerical experiments of Section~\ref{num:section} suggest a decay of order $t^{-\alpha}$ of the distance to the stationary profile when $u_0=\varphi$ (for the truncated problem computed there). A decay of the form $E_\alpha(-\mu t^\alpha)\sim t^{-\alpha}$ in $L^2$ is obtained in~\cite{Li-Salgado-2023} for time-fractional gradient flows of uniformly convex energies, but this mechanism is not available in the whole space, and we do not obtain rates here, see Section~\ref{sec:comments}.

\noindent\textsc{Organization of the paper. } Section~\ref{sec:viscosity} introduces the notion of viscosity solution for problem~\eqref{eq:main}--\eqref{eq:initial} and proves uniqueness (\Cref{thm:intro-uniqueness}). 

Section~\ref{sec:scheme} introduces the
discretization scheme, establishes its stability and regularity estimates, and proves its convergence to a viscosity solution, which gives existence together with the properties in \Cref{thm:intro-existence}. Section~\ref{sec:discrete-asymptotics} studies the large-time behaviour of the numerical scheme and identifies its limit as the minimal solution of the discrete elliptic obstacle problem. An analogous analysis for the continuous problem is performed in~\Cref{sec:asymptotics}, which starts with the study of the elliptic obstacle problem and its minimal solution, and where we prove~\Cref{thm:intro-asymptotics}. Section~\ref{num:section} presents numerical experiments in dimension one, based on an explicit stationary solution with known contact set, which test the scheme under mesh refinement and illustrate the relaxation towards the stationary profile. The paper closes with Section~\ref{sec:comments}, which contains a complement on the approximation of the continuous elliptic obstacle problem and some open problems, and an
Appendix that contains auxiliary results on superharmonic functions used in Section~\ref{sec:asymptotics}.

\section{Viscosity solutions of the time-fractional obstacle problem: definition and uniqueness}
\label{sec:viscosity}

We begin by making precise the notion of viscosity solution for the time-fractional obstacle problem.

\begin{definition}[Viscosity solution of the equation]\label{def:viscsol-p-Laplace}
    Let $s\in(0,1]$ and $\alpha\in(0,1)$. We say that a bounded lower (resp.\ upper) semicontinuous function $u$ in $\Q\coloneqq\R^d\times(0,\infty)$ is a \emph{viscosity supersolution} (resp.\ \emph{subsolution}) of equation~\eqref{eq:main} if the following condition is satisfied: for every $(x_0,t_0)\in\Q$ and $r>0$ such that the parabolic cylinder
    \[
        \mathcal{B}_r\coloneqq\overline{B_r}(x_0)\times[t_0-r,t_0]
    \]
    is contained in $\Q$, and every $\phi\in C^2(\mathcal{B}_r)\cap C_\textup{b}(\overline{\Q}\setminus\mathcal{B}_r)$, where $\overline{\Q}\coloneqq\R^d\times[0,\infty)$, satisfying
    \begin{align}
        \label{eq:visctouch-super}
        \phi(x,t)-u(x,t)&< \phi(x_0,t_0)-u(x_0,t_0) \quad \textup{for all }(x,t)\in \mathcal{B}_r\setminus\{(x_0,t_0)\}\\\label{eq:visctouch-sub}
        \text{(resp. }\phi(x,t)-u(x,t)&> \phi(x_0,t_0)-u(x_0,t_0)),
    \end{align}
    and
    \begin{align}
        \label{eq:visctouch-super2}
        \phi(x,t)-u(x,t)&\leq \phi(x_0,t_0)-u(x_0,t_0) \quad \textup{for all } (x,t)\in \Q\\
        \label{eq:visctouch-sub2}
        \text{(resp. }\phi(x,t)-u(x,t)&\geq \phi(x_0,t_0)-u(x_0,t_0)),
    \end{align}
    we have
    \begin{align}
        \label{eq:viscsuper}
        \min\{\partial_t^\alpha \phi(x_0,t_0)+\fl \phi(x_0,t_0) \coma u(x_0,t_0) - \varphi(x_0)\} &\geq 0\\
        \label{eq:viscsub}
        \text{(resp. }  \min\{\partial_t^\alpha \phi(x_0,t_0)+\fl \phi(x_0,t_0) \coma u(x_0,t_0) - \varphi(x_0)\} &\leq 0).
    \end{align}

    A \emph{viscosity solution} is a bounded continuous function $u$ in $\R^d\times(0,\infty)$ that is both a viscosity supersolution and a viscosity subsolution.

    Here $\partial_t^\alpha\phi(x_0,t_0)$ is understood through the difference formula~\eqref{eq:capdef}, and $\fl\phi(x_0,t_0)$ through the integral formula~\eqref{eq:definition.fractional.laplacian}. Both are well defined for the test functions above.
\end{definition}

\begin{remark}\label{rem:viscbetter}
    {\rm The class of admissible test functions in Definition~\ref{def:viscsol-p-Laplace} can be extended in the following ways, all of which follow by standard approximation arguments.
    \begin{enumerate}[\rm (a)]
        \item\label{item:moregeneraltesta} We could accept more general test functions satisfying
            \begin{align*}
                \phi(x,t)-u(x,t)&\leq \phi(x_0,t_0)-u(x_0,t_0) \quad \textup{for all }(x,t)\in\Q\\
                \text{(resp. }\phi(x,t)-u(x,t)&\geq \phi(x_0,t_0)-u(x_0,t_0) \text{)}.
            \end{align*}

        \item\label{item:moregeneraltestb} Without loss of generality we can replace \eqref{eq:visctouch-super2} (resp.\ \eqref{eq:visctouch-sub2}) by
        \begin{align*}
            \phi(x,t)-u(x,t)&\leq \phi(x_0,t_0)-u(x_0,t_0)-\eta \quad \textup{for all } (x,t)\in \Q\setminus\mathcal{B}_r\\
            (\textup{resp. } \phi(x,t)-u(x,t)&\geq \phi(x_0,t_0)-u(x_0,t_0)+\eta)
        \end{align*}
        for some $\eta>0$.

        \item\label{item:moregeneraltestc} We may assume in addition that $\phi\in C^\infty_\textup{b}(\overline{\Q})$.
    \end{enumerate}}
\end{remark}

We can now extend this notion to the Cauchy problem~\eqref{eq:main}--\eqref{eq:initial}, incorporating the initial condition.

\begin{definition}[Viscosity solution of the Cauchy problem]\label{def:viscsol-Cauchy}
    Let $s\in(0,1]$ and $\alpha\in(0,1)$. We say that a bounded lower (resp.~upper) semicontinuous function $u$ in $\R^d\times[0,\infty)$ is a \emph{viscosity supersolution} (resp.~\emph{subsolution}) of \eqref{eq:main}--\eqref{eq:initial} if the following holds:
    \begin{enumerate}[\rm (a)]
        \item $u$ is a viscosity supersolution (resp.\ subsolution) of~\eqref{eq:main} for $(x,t)\in \R^d\times(0,\infty)$, in the sense of \Cref{def:viscsol-p-Laplace}.

        \item $u(x,0)\geq u_0(x)$ (resp.\ $u(x,0)\leq u_0(x)$) for all $x\in  \R^d$.
    \end{enumerate}
    A \emph{viscosity solution} of~\eqref{eq:main}--\eqref{eq:initial} is a bounded continuous function $u$ in $\R^d\times[0,\infty)$ that is both a viscosity supersolution and a viscosity subsolution. In particular, $u(x,0)=u_0(x)$     for all $x\in\R^d$.
\end{definition}

\Cref{def:viscsol-p-Laplace} only requires a viscosity supersolution to lie above the obstacle at the points where some test function touches it from below, see~\eqref{eq:viscsuper}; in particular, it is not clear a priori that it lies above the obstacle everywhere. This is however the case for continuous supersolutions, as we now show.

\begin{lemma}\label{lem:above-obstacle}
    Let $s\in(0,1]$, $\alpha\in(0,1)$ and $\varphi\in C_\textup{b}(\R^d)$. If $v\in C_\textup{b}(\R^d\times(0,\infty))$ is a viscosity supersolution of~\eqref{eq:main}, then $v\ge\varphi$ in $\R^d\times(0,\infty)$.
\end{lemma}

\begin{proof}
   The idea is to produce a test function touching $v$ from below at a point where $v<\varphi$, so that the obstacle term in~\eqref{eq:viscsuper} yields a contradiction; since a continuous function need not admit such a test function at a prescribed point, we use a quadratic penalization, which creates a contact point nearby.
    
   Assume by contradiction that the open set $U\coloneqq\{(x,t)\in\Q:\ v(x,t)<\varphi(x)\}$ is nonempty, and take $(x_1,t_1)\in U$ and $\rho>0$ such that $\mathcal C\coloneqq\overline{B_\rho}(x_1)\times[t_1-\rho,t_1]\subset U$. For $K>0$, let $(x_K,t_K)$ be a minimum point over $\mathcal C$ of the continuous function $g_K(x,t)\coloneqq v(x,t)+K(|x-x_1|^2+(t-t_1)^2)$. Since $g_K(x_K,t_K)\le g_K(x_1,t_1)=v(x_1,t_1)$, we have $K(|x_K-x_1|^2+(t_K-t_1)^2)\le2\|v\|_{L^\infty}$, so that $(x_K,t_K)\to(x_1,t_1)$ as $K\to\infty$; thus, for $K$ large, the cylinder $\mathcal B_r=\overline{B_r}(x_K)\times[t_K-r,t_K]$ with $r=\rho/2$ is contained in $\mathcal C\subset U$. Define $\phi\coloneqq g_K(x_K,t_K)-K(|x-x_1|^2+(t-t_1)^2)$ in $\mathcal B_r$ and $\phi\coloneqq-\|v\|_{L^\infty}-1$ in $\overline{\Q}\setminus\mathcal B_r$ (the jump across $\partial\mathcal B_r$ is allowed by Definition~\ref{def:viscsol-p-Laplace}, and the values of the nonlocal terms play no role in the argument). Then $\phi\in C^2(\mathcal B_r)\cap C_\textup{b}(\overline{\Q}\setminus\mathcal B_r)$ and $\phi-v\le0=(\phi-v)(x_K,t_K)$ in $\Q$, so that $\phi$ is an admissible test function for $v$ at $(x_K,t_K)$, see \Cref{rem:viscbetter}\eqref{item:moregeneraltesta}. By~\eqref{eq:viscsuper}, $v(x_K,t_K)\ge\varphi(x_K)$, which contradicts $(x_K,t_K)\in U$.
\end{proof}

Our goal now is to prove that there is at most one viscosity solution to the Cauchy problem~\eqref{eq:main}--\eqref{eq:initial}.
\begin{theorem}\label{thm:uniqueness}
    Let $s\in(0,1]$, $\alpha\in (0,1)$ and $u_0,\varphi\in C_\textup{b}(\R^d)$ be such that $u_0\geq \varphi$ in $\R^d$. If $u$ is a viscosity solution to problem~\eqref{eq:main}--\eqref{eq:initial}, then
    \begin{gather*}
        u(x,t)=u^*(x,t):= \inf\limits_{v\in \mathcal{V}} v(x,t),\quad\textup{where}\\
        \mathcal{V}= \{v\in C_\textup{b}(\R^d\times[0,\infty)): v \textup{ is a viscosity supersolution of \eqref{eq:main}--\eqref{eq:initial}}\}.
    \end{gather*}
    In particular, problem~\eqref{eq:main}--\eqref{eq:initial} has at most one viscosity solution.
\end{theorem}

Notice that $\mathcal{V}\neq\emptyset$, since any sufficiently large constant is a viscosity supersolution to~\eqref{eq:main}--\eqref{eq:initial}; in particular, the constant  $\|u_0\|_{L^\infty(\R^d)}$ is one of such supersolutions. Moreover, every $v\in\mathcal{V}$ satisfies $v\ge\varphi$ in $\R^d\times(0,\infty)$ by \Cref{lem:above-obstacle}, and also at $t=0$, since $v(\cdot,0)\ge u_0 (\cdot)\ge\varphi (\cdot)$; hence the function $u^*$ is well defined, and $u^*\ge \varphi$.

\begin{remark}\label{rem:Perron}
   {\rm By \Cref{thm:uniqueness}, if a viscosity solution exists, it is given by $u^*$. Perron's method would consist in showing directly that $u^*$ is a viscosity solution, as done in \cite{Erwin} for equations with the Caputo derivative outside the nonlinearity; it could presumably be adapted to the obstacle problem, where the derivative lies inside the minimum. We follow instead the numerical scheme of \Cref{sec:scheme}, which gives existence together with further regularity properties of the solution.}
\end{remark}

By the very definition of $u^*$, we have that $u^*\leq u$ for any viscosity solution $u$. For the reverse inequality, we will compare $u$ with each $v\in\mathcal{V}$ by means of the following comparison principle for the time-fractional heat equation posed in open subsets $\Omega$ of $\R^{d}\times (0,\infty)$. In our case $\Omega$ will be the \emph{non-coincidence} set, where $u>\varphi$. Here, viscosity sub- and supersolutions of the linear equation in $\Omega$ are understood as in Definition~\ref{def:viscsol-p-Laplace}, without the obstacle term, and requiring the viscosity inequality only at points $(x_0,t_0)\in\Omega$ and for cylinders $\mathcal B_r\subset\Omega$.

The comparison principle of \cite[Theorem~4.1]{Erwin}, from which the lemma is derived, holds for semicontinuous sub- and supersolutions. We state our version for continuous functions, which is all we need. The reason is explained in the proof.

\begin{lemma}\label{lem:compErwin}
    Let $s\in(0,1]$, $\alpha\in(0,1)$ and $\Omega\subset\R^{d}\times (0,\infty)$ be an open set. Let also $u,v\in C_\textup{b}(\R^{d}\times [0,\infty))$ be such that
    \begin{equation*}
        \begin{cases}
            \partial_t^\alpha v(x,t)+\fl v(x,t) \leq 0 \quad &\textup{for all }(x,t)\in \Omega\textup{ in the viscosity sense},\\
            \partial_t^\alpha u(x,t)+\fl u(x,t) \geq 0 \quad &\textup{for all } (x,t)\in \Omega  \textup{ in the viscosity sense},\\
            v(x,t)\leq u(x,t) \quad &\textup{for all }(x,t)\in (\R^d\times(0,\infty))\setminus\Omega,\\
            v(x,0)\leq u(x,0) \quad &\textup{for all }x\in \R^d.
        \end{cases}
    \end{equation*}
    Then $v\leq u$ in $\R^d\times[0,\infty)$. 
    
    The same conclusion holds on a finite horizon: if, for some $T>0$, the above hypotheses hold with $(0,\infty)$ replaced by $(0,T]$ (that is, $\Omega\subset\R^d\times(0,T]$ and $v\le u$ on $(\R^d\times(0,T])\setminus\Omega$), then $v\le u$ in $\R^d\times[0,T]$.
\end{lemma}

\begin{proof}
This lemma is a variant of \cite[Theorem~4.1]{Erwin}, a comparison principle for a broad class of fully nonlinear equations involving Caputo time derivatives and nonlocal spatial operators, which includes the linear equation $\partial_t^\alpha w+\fl w=0$. The result in \cite{Erwin} is the finite-horizon version of the lemma in the particular case $\Omega=\R^d\times(0,T]$, in which the equations hold in the whole strip and the ordering $v\le u$ is only assumed at $t=0$. The proof adapts to an arbitrary open set $\Omega\subset\R^d\times(0,T]$, with the ordering assumed on the complement of $\Omega$. Since $T$ is arbitrary, this also gives the lemma on $\R^d\times[0,\infty)$.

The argument starts by regularizing $u$ and $v$ in time by means of inf- and sup-convolutions. This makes them Lipschitz in time, so that their Caputo derivatives are well defined. One then penalizes in space and doubles variables. The positive maxima involved in the argument are attained in a compact set, thanks to the penalization. Moreover, for small values of the regularization parameter they occur inside $\Omega$, where the equations are available. This follows from the assumption on the complement of $\Omega$, since $u$ and $v$ are continuous and the convolutions converge to them.

The argument uses the notion of viscosity solution of \cite[Definition~2.1]{Erwin}. There, the test function is only required to touch the sub- or the supersolution in a cylinder around the contact point, and the function itself is used in the nonlocal terms away from it. This is where the continuity of $u$ and $v$ comes in. Given such a test function $\phi$ for $u$ at $(x_0,t_0)$, consider the function equal to $\phi$ in $\mathcal B_r$ and to $u+(\phi-u)(x_0,t_0)$ outside. It is continuous and touches $u$ globally. Hence it is an admissible test function in \Cref{def:viscsol-p-Laplace}, by \Cref{rem:viscbetter}\eqref{item:moregeneraltesta}, and its nonlocal terms at $(x_0,t_0)$ are the same as in \cite{Erwin}. The same applies to $v$. Thus $u$ and $v$ are sub- and supersolutions in the sense of \cite{Erwin}, and the argument there applies to them.
\end{proof}

We now have all the necessary ingredients to prove our uniqueness result.

\begin{proof}[Proof of Theorem~\ref{thm:uniqueness}]
Since $u\in\mathcal{V}$, we have $u^*\le u$ by definition of $u^*$. To prove the reverse inequality, fix $v\in\mathcal{V}$; we have to show that $u\le v$ in $\R^d\times[0,\infty)$.

We define the set $\Omega=\{(x,t)\in \R^d\times(0,\infty)\, : \, u(x,t) > \varphi(x)\}$, which is open since both $u$ and $\varphi$ are continuous. Since $u$ is a viscosity solution to~\eqref{eq:main} and $u>\varphi$ in $\Omega$, we have
\begin{align*}
   \partial_t^\alpha u(x,t)+\fl u(x,t) \leq  0 \quad \textup{for all }(x,t)\in \Omega \textup{ in the viscosity sense}.
\end{align*}
Moreover, since $v$ is a viscosity supersolution to~\eqref{eq:main} in $\R^d\times(0,\infty)$, we have
\begin{align*}
   \partial_t^\alpha v(x,t)+\fl v(x,t) \geq  0 \quad \textup{for all }(x,t)\in \Omega\textup{ in the viscosity sense}.
\end{align*}
On the other hand, by \Cref{lem:above-obstacle} we have $u\ge\varphi$ and $v\ge\varphi$ in $\R^d\times(0,\infty)$, so that $u=\varphi\le v$ on $(\R^d\times(0,\infty))\setminus\Omega$, that is,
\begin{align*}
    u(x,t)\le v(x,t) \quad \textup{for all }  (x,t)\in (\R^d\times(0,\infty))\setminus\Omega.
\end{align*}
Finally, $u(\cdot,0)=u_0\le v(\cdot,0)$ in $\R^d$, since $v$ is a viscosity supersolution of the Cauchy problem. The inequality $u\le v$ now follows from~\Cref{lem:compErwin}. Taking the infimum over $v\in\mathcal{V}$ we conclude that $u\le u^*$.
\end{proof}

\section{Existence and properties via a numerical approximation}
\label{sec:scheme}

In this section we construct the viscosity solution of \eqref{eq:main}--\eqref{eq:initial} as the limit of an explicit finite-difference scheme, which is monotone under a suitable CFL-type condition. The scheme is obtained by combining the classical L1 \emph{discretization} of the Caputo derivative with the \emph{powers of the discrete Laplacian} discretization of $\fl$, replacing each term in~\eqref{eq:main} by its discrete counterpart. In Section~\ref{subsec:scheme} we introduce the scheme and its basic properties. In Section~\ref{subsec:stability} we prove that it is $L^\infty$-stable, that it satisfies an $L^\infty$-contraction property and a comparison principle, and, most importantly, that its solutions are equicontinuous in space and time, uniformly in the discretization parameters, provided the data are Hölder continuous. These estimates, besides being of independent interest for the
numerical approximation of the problem, are the key to pass to the limit in Section~\ref{subsec:convergence}, where we prove convergence of the scheme to the unique viscosity solution and deduce the properties of the latter stated in~\Cref{thm:intro-existence}.

\subsection{Formulation of the numerical scheme}
\label{subsec:scheme}

\emph{Discretization of the time derivative}. We consider a uniform-in-time grid parametrized by a time step $\tau>0$. More precisely, define $t_n\coloneqq\tau n \in \tau \N$. For a function $\Phi: \tau \N\to\R$, we denote $\Phi^n\coloneqq\Phi(\tau n)$. We now introduce the discretization of the Caputo derivative, for which we use the so-called \emph{L1} discretization that is given by
\begin{align}\label{eq:def.discretization.Caputo}
    D_\tau^\alpha \Phi^n:= \frac{\tau^{-\alpha}}{\Gamma(2-\alpha)}\left( \Phi^n - b_{n-1}\Phi^0-\sum_{k=1}^{n-1}(b_{n-k-1}-b_{n-k})\Phi^k\right), \quad \textup{where } b_i:=(i+1)^{1-\alpha}-i^{1-\alpha}.
\end{align}
This discretization is obtained by applying the Caputo derivative to the piecewise linear interpolant $\mathcal{I}_\tau[\Phi]$ of $\Phi$. Thus, starting from expression \eqref{eq:capdef} we get
\begin{align}\label{eq:discretization.Caputo.2}
    \begin{aligned}
        D^\alpha_\tau \Phi^n& :=\partial^\alpha_t (\mathcal{I}_\tau[\Phi])(t_n)\\
        &=\frac{1}{\Gamma(1-\alpha)}\frac{\Phi^n-\Phi^0}{t_n^{\alpha}}+\frac{\alpha}{\Gamma(1-\alpha)}\sum_{j=0}^{n-1}\int_{t_j}^{t_{j+1}} \frac{\Phi^n-\left(\frac{t_{j+1}-s}{\tau} \Phi^j + \frac{s-t_{j}}{\tau}\Phi^{j+1}\right)}{(t_n-s)^{1+\alpha}} \dd s,
    \end{aligned}
\end{align}
from which \eqref{eq:def.discretization.Caputo} follows by direct computations.

\begin{remark}
   {\rm We note that $b_{n-k-1}-b_{n-k}\geq0$, since the sequence $\{b_{i}\}_{i=1}^\infty$ is decreasing, which follows from the concavity of the function $f(r)=r^{1-\alpha}$.}
\end{remark}

\noindent\emph{Discretization of the fractional/local Laplacian.} We consider a uniform-in-space grid parametrized by a spatial step $h>0$. More precisely, define $x_\beta \coloneqq h\beta \in h \Z^d$. We now introduce the discretization of the fractional/local Laplacian, for which we employ the so-called \emph{powers of the discrete Laplacian} approach, see \cite{Ciaurri2018} and also \cite{delTeso2019a,delTeso2018b}. It is defined for a function $\Psi: \R^d\to \R$ as
\begin{align}\label{eq:discfldef}
    \fld \Psi(x) \coloneqq \frac{1}{h^{2s}}\sum_{\beta \not=0} (\Psi(x) - \Psi(x+x_\beta))\omega_\beta,
\end{align}
with the  weights of the discretization, $\omega_\beta$, given by
\begin{align*}
    \omega_\beta&=
    \begin{cases}
        1& \textup{if }|\beta|=1,\\
        0& \textup{if }|\beta|\not=1,
    \end{cases}
    &&\textup{if }s=1,\\
    \omega_\beta&=\frac{1}{|\Gamma(-s)|}\int_{0}^\infty G(\beta,t) \frac{\dd t}{t^{1+s}}, \quad \textup{with } G(\beta,t)=e^{-2dt}\prod_{i=1}^d I_{|\beta_i|}(2t),&&\textup{if }s\in(0,1),
\end{align*}
where $I_k\geq0$ are the modified Bessel functions of the first kind and order $k\in \N$.

For $s=1$, \eqref{eq:discfldef} reduces to the standard central difference approximation of the Laplacian. For $s\in(0,1)$, it is obtained by taking the fractional power of order $s$ of the discrete Laplacian operator, that~is,
\[
    \fld \Psi(x) =\frac{1}{\Gamma(-s)}\int_{0}^\infty \left(e^{t\Delta_h}\Psi(x)-\Psi(x)\right) \frac{\dd t}{t^{1+s}},
\]
and recalling the explicit expression of the discrete Laplacian semigroup
\[
    e^{t\Delta_h}\Psi(x)= \sum_{\beta\in \Z^d} \Psi(x-x_\beta) G(\beta,t/h^2).
\]
From here, \eqref{eq:discfldef} follows by direct computations, noting that $\sum_{\beta\in \Z^d} G(\beta,t)=1$ for all $t\geq 0$.
Let us also define the quantity given by the sum of the weights
\begin{align*}
   \Lambda\coloneqq \sum_{\beta\not=0} \omega_\beta =
   \begin{cases}
       2d  &\textup{if }s=1,\\
       \displaystyle\frac{1}{|\Gamma(-s)|}\int_{0}^\infty (1-  G(0,t)) \frac{\dd t}{t^{1+s}}\quad&\textup{if }s\in(0,1).
   \end{cases}
\end{align*}

\begin{remark}
    {\rm Let us note that $\omega_\beta=\omega_{-\beta}$ for all $\beta \in \Z^d$, and thus, the weights are a symmetric function in $\Z^d$. We also note that $\Lambda<+\infty$, and thus, $\omega\in \ell^1(\Z^d)$.}
\end{remark}

\begin{remark}
    {\rm We have chosen the powers of the discrete Laplacian discretization because it is known to be second-order accurate for every $s\in(0,1]$, see~\cite{delTeso2018b}. However, the results of this paper also hold for more general discretizations, as long as the corresponding weights are nonnegative, symmetric, summable, and satisfy the two-sided bound~\eqref{eq:weightsbound} below. We refer to \cite{delTeso2018b} for a summary of alternative discretizations that are suitable in this context.}
\end{remark}

\noindent\emph{Numerical scheme}.  We consider the following explicit-in-time scheme associated with \eqref{eq:main}, obtained by replacing each term in \eqref{eq:main} with its discrete counterpart:
\begin{equation}\label{eq:numsch}
    \min\left\{D_\tau^\alpha U^{j+1}(x) + \fld U^j(x)\coma U^{j+1}(x)-\varphi(x)\right\}=0, \quad x\in \R^d, \quad j\in \N.
\end{equation}
In order to compare the solution of the scheme with functions of $(x,t)\in\R^d\times[0,\infty)$, we write $U_{\tau,h}(x,t)\coloneqq\mathcal{I}_\tau[U](x,t)$ for the piecewise linear interpolant in time of the sequence $\{U^n\}_{n\in\N}$, so that $U_{\tau,h}(\cdot,t_n)=U^n$. Estimates of the form $\|U^m-U^n\|_{L^\infty(\R^d)}\le\omega(|t_m^\alpha-t_n^\alpha|)$, with $\omega$ nondecreasing, transfer to $U_{\tau,h}$, with $\omega$ replaced by $3\omega$, by the triangle inequality through the nearest grid times.

Note that the above scheme can be equivalently written as
\begin{equation}\label{eq:schmNO-CFL}
    \begin{aligned}
        U^{j+1}(x)&= \max\bigg\{b_ j U^0(x) + \sum_{k=1}^{j} (b_{j-k}-b_{j+1-k}) U^k(x)\\
        &\quad+\Gamma(2-\alpha)\frac{\tau^\alpha}{h^{2s}}\sum_{\beta\not=0} (U^{j}(x+x_\beta)-U^{j}(x))\omega_\beta \coma \varphi(x) \bigg\}
    \end{aligned}
\end{equation}
or
\begin{equation}\label{eq:schmCFL}
    \begin{aligned}
        U^{j+1}(x)&=\max\bigg\{b_ j U^0(x)+\sum_{k=1}^{j-1}(b_{j-k}-b_{j+1-k})U^k(x)+((b_0-b_1)-\Lambda\Gamma(2-\alpha)\frac{\tau^\alpha}{h^{2s}} )U^j(x)\\
        &\quad+\Gamma(2-\alpha)\frac{\tau^\alpha}{h^{2s}}  \sum_{\beta\not=0} U^{j}(x+x_\beta)\omega_\beta \coma \varphi(x) \bigg\}.
    \end{aligned}
\end{equation}
It is then clear that, to ensure the monotonicity of the scheme and its good properties, the following CFL-type stability condition is required:
\begin{align}\label{as:CFL}\tag{CFL}
    \tau^{\alpha} \leq \frac{2-2^{1-\alpha}}{\Lambda \Gamma(2-\alpha)}h^{2s},
\end{align}
where we have used that $b_0-b_1= (1^{1-\alpha}-0^{1-\alpha})-(2^{1-\alpha}-1^{1-\alpha})=2-2^{1-\alpha}$.

\subsection{Stability and regularity properties of the numerical scheme}
\label{subsec:stability}

Let us first show that the numerical scheme is stable in $L^\infty(\R^d)$.

\begin{lemma}\label{lem:est-schm}
     Let $s\in(0,1]$, $\alpha\in (0,1)$. Assume \eqref{as:CFL}. Let $u_0,\varphi\in C_{\textup{b}}(\R^d)$ be such that $u_0\geq \varphi$ in $\R^d$. Then there exists a unique solution of \eqref{eq:numsch} with initial datum $U^0=u_0$. Moreover, the following properties hold for all $j\in \N$:
    \begin{enumerate}[\rm (a)]
        \item\label{lem:est-schm-aboveobs} $U^j\geq \varphi$ in $\R^d$;
        \item\label{lem:est-schm-stable} $\displaystyle \sup_{x\in \R^d}U^j(x)\leq \sup_{x\in \R^d}u_0(x)$.
    \end{enumerate}
\end{lemma}

\begin{proof}
Existence and uniqueness of solutions follow immediately from the explicit nature of the scheme. Property~\eqref{lem:est-schm-aboveobs} follows directly from \eqref{eq:schmNO-CFL}. To prove \eqref{lem:est-schm-stable} we proceed by induction. The case $j=0$ is trivial. Assume by induction that it is true for $k=0,\ldots, j$ and denote $\widehat{U}^j= \sup_{x\in \R^d}U^j(x)$. We use formulation~\eqref{eq:schmCFL} of the scheme together with \eqref{as:CFL} and the identity
\begin{equation}\label{eq:property.coefficients}
    b_j + \sum_{k=1}^{j} (b_{j-k}-b_{j+1-k})=1
\end{equation}
to get
\begin{align*}
    U^{j+1}(x)&\leq \max\bigg\{ b_ j \widehat{U}^0 + \sum_{k=1}^{j-1} (b_{j-k}-b_{j+1-k}) \widehat{U}^k+((b_0-b_1)-\Lambda\Gamma(2-\alpha)\frac{\tau^\alpha}{h^{2s}} )\widehat{U}^j \\
    &\qquad \qquad +\Lambda\Gamma(2-\alpha)\frac{\tau^\alpha}{h^{2s}} \widehat{U}^j\coma \widehat{U}^0 \bigg\}\\
    &\leq\max \bigg\{ (b_j + \sum_{k=1}^{j} (b_{j-k}-b_{j+1-k})) \widehat{U}^0 \coma \widehat{U}^0  \bigg\}\\
    & =\widehat{U}^0.
\end{align*}
We conclude the proof by taking the supremum in $\R^d$ of the left-hand side of the above estimate.
\end{proof}

We next prove that the scheme satisfies an $L^\infty(\R^d)$-contraction property (which in particular implies spatial equicontinuity) and establish a comparison principle.

\begin{lemma}\label{lem:contrac-equicont}
    Let $s\in(0,1]$, $\alpha\in (0,1)$. Assume \eqref{as:CFL}. Let $u_0,v_0,\varphi,\psi\in C_{\textup{b}}(\R^d)$ be such that $u_0\geq \varphi$ and $v_0\geq \psi$ in $\R^d$. Let $\{U^{n}\}_{n\in \N}$ and $\{V^{n}\}_{n\in \N}$ be the solutions of \eqref{eq:numsch} corresponding to the initial data $U^0=u_0$ and $V^0=v_0$, respectively. The following properties hold for all $j\in \N$:
     \begin{enumerate}[\rm (a)]
         \item\label{lem:contrac-equicont-itemcontrac} $\|U^{j}-V^{j}\|_{L^\infty(\R^d)} \leq \max\{\|U^{0}-V^{0}\|_{L^\infty(\R^d)}\coma\|\varphi-\psi\|_{L^\infty(\R^d)}\}$;

         \item\label{lem:contrac-equicont-itemequicont} $\|U^{j}(\cdot+y)-U^{j}\|_{L^\infty(\R^d)} \leq \max\{\|U^{0}(\cdot+y)-U^{0}\|_{L^\infty(\R^d)}\coma\|\varphi(\cdot+y)-\varphi\|_{L^\infty(\R^d)}\}$ for all $y\in \R^d$;

         \item\label{lem:contrac-equicont-comparison} If additionally $u_0\leq v_0$ and $\varphi \leq \psi$, then $U^j\leq V^j$ in $\R^d$.
     \end{enumerate}
\end{lemma}

\begin{proof}
\eqref{lem:contrac-equicont-itemcontrac} We recall the inequality
\begin{align}\label{eq:ineqmax}
    \big|\max\{A_1,A_2\}-\max\{B_1,B_2\}\big| \leq \max\{|A_1-B_1|\, , \, |A_2-B_2|\},
\end{align}
valid for any $A_1,A_2,B_1,B_2\in\R$. The result is trivial for $j=0$.  Assume by induction that it is true for $k=0,\ldots,j$ and denote $W^{j}=U^{j}-V^{j}$. We use formulation \eqref{eq:schmCFL} together with \eqref{as:CFL} and \eqref{eq:ineqmax} to get
\begin{align*}
    |W^{j+1}(x)|&\leq\max\bigg\{ b_ j |W^0(x)| + \sum_{k=1}^{j-1} (b_{j-k}-b_{j+1-k}) |W^k(x)| + ((b_0-b_1)-\Lambda\Gamma(2-\alpha)\frac{\tau^\alpha}{h^{2s}} )|W^j(x)|\\
    &\quad+\Gamma(2-\alpha)\frac{\tau^\alpha}{h^{2s}}  \sum_{\beta\not=0} |W^j(x+x_\beta)| \omega_\beta\, , \,|\varphi(x)-\psi(x)| \bigg\}\\
    &\leq \max\bigg\{ b_ j \|W^0\|_{L^\infty(\R^d)} + \sum_{k=1}^{j-1} (b_{j-k}-b_{j+1-k}) \|W^k\|_{L^\infty(\R^d)}\\
    &\quad+ ((b_0-b_1)-\Lambda\Gamma(2-\alpha)\frac{\tau^\alpha}{h^{2s}} )\|W^j\|_{L^\infty(\R^d)}+\Lambda\Gamma(2-\alpha)\frac{\tau^\alpha}{h^{2s}} \|W^j\|_{L^\infty(\R^d)}\,,\, \|\varphi-\psi\|_{L^\infty(\R^d)} \bigg\}\\
    &=\max\bigg\{ b_ j \|W^0\|_{L^\infty(\R^d)} + \sum_{k=1}^{j} (b_{j-k}-b_{j+1-k}) \|W^k\|_{L^\infty(\R^d)}\,,\, \|\varphi-\psi\|_{L^\infty(\R^d)} \bigg\}.
\end{align*}
From here, \eqref{lem:contrac-equicont-itemcontrac} follows by taking the supremum in $\R^d$ of the left-hand side of the above estimate and using the induction hypothesis together with~\eqref{eq:property.coefficients}.

\noindent\eqref{lem:contrac-equicont-itemequicont} For any fixed $y\in \R^d$, let us define $V^{j}\coloneqq U^{j}(\cdot+y)$, which is clearly a solution of \eqref{eq:numsch} with obstacle $\psi=\varphi(\cdot+y)$. At this point, it is enough to use \eqref{lem:contrac-equicont-itemcontrac} to conclude the proof.

\noindent\eqref{lem:contrac-equicont-comparison}
The proof is analogous to that of~\eqref{lem:contrac-equicont-itemcontrac} using the inequality
\[
    \max\{A_1,A_2\}-\max\{B_1,B_2\} \leq \max\{A_1-B_1\, , \, A_2-B_2\}
\]
instead of \eqref{eq:ineqmax}, and we omit the details.
\end{proof}

\begin{remark}
   {\rm  The above proofs only use the boundedness of the data; hence \eqref{lem:contrac-equicont-itemcontrac}--\eqref{lem:contrac-equicont-comparison} hold, with the same proofs, for $u_0,v_0,\varphi,\psi\in B(\R^d)$, which is the setting of Section~\ref{sec:discrete-asymptotics}.}
\end{remark}

We are now ready to prove that the solutions are also equicontinuous in time, provided the initial datum belongs to some Hölder space $\X^a_{\textup{b}}(\mathbb{R}^d)$.  We rely on the following estimate for the discretization of the fractional/local Laplacian. Note that it provides a uniform-in-$h$ bound in some cases.

\begin{lemma}\label{lem:bounddiscreteLap}
    Let $s\in(0,1]$ and $\phi\in \X_{\textup{b}}^{a}(\R^d)$ for some $a\in[0,2]$. Then
    \begin{equation}\label{eq:improved-bounddiscrete}
        \begin{gathered}
            \|\fld \phi\|_{L^\infty(\R^d)} \leq C_{d,s,a}\|\phi\|_{\X^a(\R^d)} f_{a,s}(h),\quad\textup{where} \quad 
            f_{a,s}(h)\coloneqq
            \begin{cases}
                h^{a-2s}&\textup{if }a<2s,\\
                |\log(h)|&\textup{if }a=2s\textup{ and }s<1,\\
                1&\textup{if }a>2s\textup{ or }a=2.
            \end{cases}
        \end{gathered}
    \end{equation}
\end{lemma}

\begin{proof}
We will use the following second-order difference bound: for $\phi\in \X^a_\textup{b}(\R^d)$ with $a\in[0,2]$,
\begin{equation}\label{eq:second-difference}
    |2\phi(x)-\phi(x+y)-\phi(x-y)|\leq 4\|\phi\|_{\X^a(\R^d)}|y|^a\quad\textup{for all } x,y\in\R^d.
\end{equation}

If $s=1$, then, by \eqref{eq:second-difference},
\begin{align*}
    |\fld \phi(x)|&=  \left|\sum_{l=1}^d \frac{2 \phi(x)-\phi(x+h e_l )- \phi(x-he_l )}{h^2}\right| \leq 4\sum_{l=1}^d \frac{\|\phi\|_{\X^a(\R^d)} h^a}{h^2}= 4 d\|\phi\|_{\X^a(\R^d)} h^{a-2s}.
\end{align*}

If $s\in(0,1)$, we use that there exist two constants $C_1,C_2>0$ depending only on $d$ and $s$ such that
\begin{align}\label{eq:weightsbound}
    \frac{C_1\,h^d}{|x_\beta|^{d+2s}}\leq \frac{\omega_\beta}{h^{2s}} \leq \frac{C_2\,h^d}{|x_\beta|^{d+2s}} \quad \textup{for all } \beta\in \Z^d\setminus\{0\};
\end{align}
see \cite{Ciaurri2018} in dimension $d=1$ and \cite{GebertRojas-Molina2020} for $d>1$. Since $|x_\beta|=h|\beta|$, this is equivalent to 
\[
    C_1|\beta|^{-d-2s}\leq\omega_\beta\leq C_2|\beta|^{-d-2s}. 
\]
Thus, by \eqref{eq:second-difference} (and the trivial bound $4\|\phi\|_{L^\infty(\R^d)}$ for $|x_\beta|\ge1$),
\begin{align*}
    |\fld \phi(x)| &= \left|\sum_{\beta \not=0} (2\phi(x) - \phi(x+x_\beta)-\phi(x-x_\beta))\frac{\omega_\beta}{2h^{2s}}\right|\\
    &\lesssim \|\phi\|_{\X^a(\R^d)} \sum_{0<|x_\beta|<1}  h^d\frac{|x_\beta|^a}{|x_\beta|^{d+2s}} + \|\phi\|_{\X^a(\R^d)}  \sum_{|x_\beta|\geq 1}  h^d\frac{1}{|x_\beta|^{d+2s}}\\
    &\lesssim \|\phi\|_{\X^a(\R^d)} \left(\int_h^1 \frac{\dd r}{r^{1+2s-a}} + \int_1^\infty \frac{\dd r}{r^{1+2s}}\right)\\
    &\lesssim \|\phi\|_{\X^a(\R^d)} f_{a,s}(h),
\end{align*}
where we have used that, thanks to the factor $h^d$, the lattice sums are Riemann sums with mesh $h$ of the corresponding integrals.
\end{proof}

The following lemma establishes $\alpha$-H\"older continuity in time (more precisely, Lipschitz continuity with respect to $t^\alpha$) for the scheme. The bound is uniform in $h$ only when $a>2s$ or $a=2$; the case of less regular initial data will be treated in \Cref{lem:equiintime-scheme} below.

\begin{lemma}\label{lem:timeequi-nosmooth}
    Let $s\in(0,1]$ and $\alpha\in(0,1)$. Assume \eqref{as:CFL}. Let $u_0\in \X_{\textup{b}}^{a}(\R^d)$ with $a\in[0,2]$ and $\varphi\in C^0_{\textup{b}}(\R^d)$ be such that $u_0\geq \varphi$ in $\R^d$. Let $\{U^{n}\}_{n\in \N}$ be the solution of \eqref{eq:numsch} with initial datum $U^0=u_0$. Then, for all $m\geq n$, we have that
    \begin{align}\label{eq:timcont-nosmoothdata}
        \|U^m-U^n\|_{L^\infty(\R^d)} \leq C_{d,s,a,\alpha} \|u_0\|_{\X^a(\R^d)}f_{a,s}(h) |t_m^\alpha-t_n^\alpha|.
    \end{align}
\end{lemma}

\begin{proof}
We will use the following notation for convenience:
\begin{align}
    \delta^1 U^j(x)= U^{j+1}(x)-U^{j}(x).
\end{align}

We first prove the case $n=j$ and $m=j+1$ for $j\in \N$. We start by assuming that $j=0$. In this case we get, directly from \eqref{eq:schmNO-CFL}, that
\begin{align*}
    \frac{U^{1}(x)-U^0(x)}{\Gamma(2-\alpha)\tau^\alpha}&= \max\bigg\{-\fld U^0(x) \coma \frac{\varphi(x)-U^0(x)}{\Gamma(2-\alpha)\tau^\alpha}  \bigg\}\leq C_{d,s,a}\|u_0\|_{\X^a(\R^d)} f_{a,s}(h),
\end{align*}
where we have used \eqref{eq:improved-bounddiscrete}  together with the fact that $U^0\geq \varphi$ in $\R^d$. On the other hand,
\begin{align*}
    \frac{U^{1}(x)-U^0(x)}{\Gamma(2-\alpha)\tau^\alpha}&\geq -\fld U^0(x) \geq -C_{d,s,a}\|u_0\|_{\X^a(\R^d)} f_{a,s}(h).
\end{align*}
Thus, we have obtained
\begin{align}\label{eq:estt0}
    \|\delta^1 U^{0}\|_{L^\infty(\R^d)} \leq C_{d,s,a}\|u_0\|_{\X^a(\R^d)} f_{a,s}(h)\Gamma(2-\alpha) \tau^\alpha,
\end{align}
which implies~\eqref{eq:timcont-nosmoothdata}, since $\tau^\alpha=t_1^\alpha$.

Let us assume now $j>0$. First we note that \eqref{eq:schmCFL} at time $t_j$ can be rewritten as
\begin{align*}
    \begin{aligned}
        U^{j}(x)&= \max\bigg\{b_{j-1} U^0(x) + \sum_{k=2}^{j-1} (b_{j-k}-b_{j+1-k}) U^{k-1}(x)  + ((b_0-b_1)-\Lambda\Gamma(2-\alpha)\frac{\tau^\alpha}{h^{2s}} )U^{j-1}(x) \\
        &\quad+\Gamma(2-\alpha)\frac{\tau^\alpha}{h^{2s}}\sum_{\beta\not=0} U^{j-1}(x+x_\beta)\omega_\beta \coma \varphi(x) \bigg\}.
    \end{aligned}
\end{align*}
We now subtract \eqref{eq:schmCFL} at times $t_{j+1}$ and $t_j$ and use \eqref{eq:ineqmax} to get
\begin{align*}
    |\delta^1 U^j(x)|&\leq\max\bigg\{\sum_{k=1}^{j-1} (b_{j-k}-b_{j+1-k})\|\delta^1 U^{k-1}\|_{L^\infty(\R^d)}+ ((b_0-b_1)-\Lambda\Gamma(2-\alpha)\frac{\tau^\alpha}{h^{2s}})\|\delta^1 U^{j-1}\|_{L^\infty(\R^d)}\\
    &\quad+\Gamma(2-\alpha)\frac{\tau^\alpha}{h^{2s}}\sum_{\beta\not=0}\|\delta^1U^{j-1}\|_{L^\infty(\R^d)}\omega_\beta\coma 0\bigg\}\\
    &=\sum_{k=1}^{j} (b_{j-k}-b_{j+1-k}) \|\delta^1U^{k-1}\|_{L^\infty(\R^d)},
\end{align*}
that is,
\begin{align}\label{eq:DisRLdec}
 \|\delta^1 U^j\|_{L^\infty(\R^d)} \leq \sum_{k=0}^{j-1} (b_{j-1-k}-b_{j-k}) \|\delta^1U^{k}\|_{L^\infty(\R^d)}.
\end{align}
Summarizing,  the sequence $f^j\coloneqq \|\delta^1 U^j\|_{L^\infty(\R^d)}$ satisfies $f^j \leq \sum_{k=0}^{j-1} (b_{j-1-k}-b_{j-k}) f^k$ for all $j\geq 1$. In the notation of~\cite{delTeso-Lukasz-2025}, equation \eqref{eq:DisRLdec} says that ${}^{RL}\partial_\tau^\alpha f^j\leq0$; thus, we can use Corollary 7 there (a discrete fractional version of Gr\"onwall's inequality) and~\eqref{eq:estt0} to get
\begin{align*}
    f^{j}\leq f^0 \tau^{1-\alpha} t_j^{\alpha-1} \leq  \Gamma(2-\alpha) C_{d,s,a}\|u_0\|_{\X^a(\R^d)} f_{a,s}(h) \tau t_j^{\alpha-1}.
\end{align*}
At this point we note that $t_j\geq  t_{j+1}/2 $ for all $j\geq1$, to get
\begin{align}\label{eq:estconsecutivetimes}
    \|\delta^1 U^j\|_{L^\infty(\R^d)} &\leq C_{d,s,a,\alpha}\|u_0\|_{\X^a(\R^d)} f_{a,s}(h) \tau t_{j+1}^{\alpha-1}.
\end{align}
From here, \eqref{eq:timcont-nosmoothdata} follows by noticing that $\tau\, t_{j+1}^{\alpha-1}\leq \frac{t_{j+1}^\alpha-t_{j}^\alpha}{\alpha}$, since $t^{\alpha-1}$ is decreasing.

We now prove the general case $m>n$. We do it by iterating the estimate for consecutive times given by \eqref{eq:estconsecutivetimes}. More precisely, for $j\geq0$ and $k\geq1$ we have
\begin{align*}
    \|U^{j+k}-U^{j}\|_{L^\infty(\R^d)} &\leq  \sum_{l=0}^{k-1}\|\delta^1 U^{j+l}\|_{L^\infty(\R^d)} \leq C_{d,s,a,\alpha}\|u_0\|_{\X^a(\R^d)} f_{a,s}(h) \tau  \sum_{l=0}^{k-1} t_{j+l+1}^{\alpha-1}.
\end{align*}
We conclude the proof by observing that
\begin{equation*}
     \tau \sum_{l=0}^{k-1} t_{j+l+1}^{\alpha-1} \leq \int_{t_j}^{t_{j+k}} t^{\alpha-1} \dd t=\frac{ t_{j+k}^\alpha- t_j^\alpha}{\alpha}.\qedhere
\end{equation*}
\end{proof}

Let us observe that \Cref{lem:timeequi-nosmooth} only provides time equicontinuity uniformly in $h$ when $f_{a,s}$ is bounded independently of $h$, that is, when $a>2s$ or $a=2$. This imposes a constraint on the regularity of the initial datum $u_0$. We now prove uniform-in-$h$ time equicontinuity for more general initial data and obstacles. In this case, the regularity of the obstacle also affects the regularity of the solution.

\begin{lemma}\label{lem:equiintime-scheme}
    Let $s\in(0,1]$, $\alpha\in (0,1)$. Assume \eqref{as:CFL}. Let $u_0\in \X_{\textup{b}}^{a}(\R^d)$ and $\varphi\in \X_{\textup{b}}^{b}(\R^d)$ for some $a,b\in(0,2s+\veps]$, with either $\veps=0$ if $s=1$ or $\veps\in(0,2-2s]$ if $s\in(0,1)$, be such that $u_0\geq \varphi$ in $\R^d$. Let $\{U^{n}\}_{n\in \N}$ be the solution of \eqref{eq:numsch} with initial datum $U^0=u_0$. Then, for all $m\geq n$, we have that
    \begin{align*}
        \|U^m-U^n\|_{L^\infty(\R^d)} &\leq \widetilde{C}_{d,s,a,b,\veps,\alpha} \max\{\|u_0\|_{\X^a(\R^d)}\coma \|\varphi\|_{\X^b(\R^d)}\}\\
        &\quad\times\max\{|t_m^\alpha-t_n^\alpha|^{\frac{a}{2s+\veps}}\coma|t_m^\alpha-t_n^\alpha|^{\frac{b}{2s+\veps+b-a}}\coma |t_m^\alpha-t_n^\alpha|\}.
    \end{align*}
\end{lemma}

 \begin{remark}
     {\rm Taking $a=2s+\veps$ in \Cref{lem:equiintime-scheme} one recovers the estimate of \Cref{lem:timeequi-nosmooth} in the range where $f_{a,s}(h)$ is bounded independently of $h$, with the constant depending also on $\|\varphi\|_{\X^b(\R^d)}$.}
 \end{remark}

\begin{proof}
Let $\rho_\delta$ be a standard radially symmetric mollifier for $\delta>0$ (the symmetry is used to get the estimate below for $a>1$), and consider the mollified version of $u_0$ and $\varphi$ given by $u_{0,\delta}= u_0*\rho_\delta$  and $\varphi_{\delta}= \varphi*\rho_\delta$. It is standard to check, given the regularity of $u_0$ and $\varphi$, that there exists a constant $K$ depending only on the dimension $d$ such that
\begin{align*}
    &\|u_{0,\delta}-u_0\|_{L^\infty(\R^d)} \leq K \|u_0\|_{\X^a(\R^d)}(\delta^a \wedge  1), \quad &&\|\varphi_{\delta}-\varphi\|_{L^\infty(\R^d)} \leq K \|\varphi\|_{\X^b(\R^d)}(\delta^b \wedge  1),\\
    &\|u_{0,\delta}\|_{L^\infty(\R^d)}\leq   \|u_0\|_{L^\infty(\R^d)}, \quad&&\|u_{0,\delta}\|_{\X^{2s+\veps}(\R^d)}\leq K  \|u_0\|_{\X^a(\R^d)}(\delta^{a-2s-\veps} \vee 1).
\end{align*}
Let $U_\delta$ be the solution with initial datum $(U_\delta)^0= u_{0,\delta}$ and obstacle $\varphi_\delta$. Note that $u_{0,\delta}\in C^\infty_\textup{b}(\R^d)$ and $u_{0,\delta}\geq \varphi_\delta$ since $u_0\geq \varphi$. Then, we apply \Cref{lem:timeequi-nosmooth} with $a=2s+\veps$ (so that $f_{a,s}(h)\leq 1$) to obtain the estimate
\begin{align*}
    \|(U_\delta)^m-(U_\delta)^n \|_{L^\infty(\R^d)} &\leq C_{d,s,2s+\veps,\alpha} \|u_{0,\delta}\|_{\X^{2s+\veps}(\R^d)} |t_m^\alpha-t_n^\alpha|\\
    &\leq \widetilde{C}_{d,s,a,\veps,\alpha}\|u_0\|_{\X^a(\R^d)}(\delta^{a-2s-\veps} \vee  1) |t_m^\alpha-t_n^\alpha|.
\end{align*}
Moreover, by the $L^\infty$-contraction property of \Cref{lem:contrac-equicont}\eqref{lem:contrac-equicont-itemcontrac}, we also get
\begin{align*}
    \|U^{j+k}&-U^{j}\|_{L^\infty(\R^d)}\\
    &\leq \|(U_\delta)^{j+k}-U^{j+k}\|_{L^\infty(\R^d)} + \|(U_\delta)^{j}-U^{j}\|_{L^\infty(\R^d)}+ \|(U_\delta)^{j+k}-(U_\delta)^{j}\|_{L^\infty(\R^d)}\\
    &\leq 2  \max\{\|u_{0,\delta}-u_0\|_{L^\infty(\R^d)}\coma\|\varphi_\delta-\varphi\|_{L^\infty(\R^d)}\} + \widetilde{C}_{d,s,a,\veps,\alpha}\|u_0\|_{\X^a(\R^d)}(\delta^{a-2s-\veps} \vee  1) |t_{j+k}^\alpha-t_j^\alpha|\\
    &\leq \widetilde{C}_{d,s,a,b,\veps,\alpha} \max\{\|u_0\|_{\X^a(\R^d)}\coma\|\varphi\|_{\X^b(\R^d)}  \} (\max\{\delta^a \wedge  1\coma\delta^b \wedge 1\}+ (\delta^{a-2s-\veps} \vee 1)(t_{j+k}^\alpha-t_j^\alpha)).
\end{align*}
Now we just need to optimize the choice of $\delta$ in terms of $(t_{j+k}^\alpha-t_j^\alpha)$. On the one hand, if $t_{j+k}^\alpha-t_j^\alpha\leq 1$ and $b\leq a$, then we make the choice $\delta = (t_{j+k}^\alpha-t_j^\alpha)^{\frac{1}{2s+\veps+b-a}}$ to get
\begin{align*}
    \max\{\delta^a \wedge  1\coma\delta^b \wedge 1\}+ (\delta^{a-2s-\veps} \vee 1)(t_{j+k}^\alpha-t_j^\alpha)= \delta^b+ \delta^{a-2s-\veps}(t_{j+k}^\alpha-t_j^\alpha)=2(t_{j+k}^\alpha-t_j^\alpha)^{\frac{b}{2s+\veps+b-a}}.
\end{align*}
On the other hand, if $t_{j+k}^\alpha-t_j^\alpha\leq 1$ and $b>a$, then we make the choice $\delta = (t_{j+k}^\alpha-t_j^\alpha)^{\frac{1}{2s+\veps}}$ to get
\begin{align*}
    \max\{\delta^a \wedge 1\coma\delta^b\wedge 1\}+ (\delta^{a-2s-\veps} \vee 1)(t_{j+k}^\alpha-t_j^\alpha)= \delta^a+ \delta^{a-2s-\veps}(t_{j+k}^\alpha-t_j^\alpha)=2(t_{j+k}^\alpha-t_j^\alpha)^{\frac{a}{2s+\veps}}.
\end{align*}
Finally, if $t_{j+k}^\alpha-t_j^\alpha> 1$, the choice $\delta=1$ leads to
\[
    \max\{\delta^a \wedge 1\coma\delta^b \wedge  1\}+ (\delta^{a-2s-\veps} \vee  1)(t_{j+k}^\alpha-t_j^\alpha)= 1+ (t_{j+k}^\alpha-t_j^\alpha)\leq2(t_{j+k}^\alpha-t_j^\alpha).\qedhere
\]
\end{proof}

\begin{remark}\label{modulusintimeexplicit}
   {\rm In the previous proof, to simplify the presentation, we have used that both $u_0$ and~$\varphi$ belong to certain Hölder spaces $\X_{\textup{b}}^{a}(\R^d)$ and $\X_{\textup{b}}^{b}(\R^d)$. Actually, the proof follows in a similar way assuming only uniform continuity for $u_0$ and $\varphi$, i.e., $u_0,\varphi\in C_{\textup{b}}^{0}(\R^d)$. More precisely, if we assume $|u_0(x)-u_{0}(y)| \leq \Lambda_{u_0}(|x-y|)$  and $|\varphi(x)-\varphi(y)| \leq \Lambda_\varphi(|x-y|)$, for certain moduli of continuity $\Lambda_{u_0}$ and~$\Lambda_{\varphi}$, we have the convolution estimates
    \begin{align*}
        \lVert u_{0,\delta}-u_{0}\rVert_{L^\infty(\R^d)}&\leq\Lambda_{u_0}(\delta) \wedge (2\lVert u_0\rVert_{L^\infty(\R^d)}),\quad&
        \lVert\varphi_{\delta}-\varphi\rVert_{L^\infty(\R^d)}&\leq\Lambda_{\varphi}(\delta)\wedge(2\lVert\varphi\rVert_{L^\infty(\R^d)}),\\
        \lVert u_{0,\delta}\rVert_{L^\infty(\R^d)}&\leq \lVert u_0\rVert_{L^\infty(\R^d)},\quad&\lVert u_{0,\delta}\rVert_{\X^{2s+\veps}(\R^d)} &\leq K\lVert u_0\rVert_{L^\infty(\R^d)}(\delta^{-2s-\veps}\vee 1),
    \end{align*}
    as well as
    \begin{align*}
        \lVert (U_\delta)^m-(U_\delta)^n\rVert_{L^\infty(\R^d)}
        &\leq \widetilde{C}_{d,s,\veps,\alpha}\lVert u_0\rVert_{L^\infty(\R^d)}(\delta^{-2s-\veps}\vee 1)|t_m^\alpha-t_n^\alpha|.
    \end{align*}
    From here, the proof follows as before, yielding a uniform modulus of continuity in time for the solution. In contrast with the Hölder case, this modulus does not, in general, admit an explicit power-law expression in the time increment, but rather depends on the spatial moduli of continuity $\Lambda_{\varphi}$ and $\Lambda_{u_0}$.}
\end{remark}

\subsection{Convergence of the numerical scheme and existence of viscosity solutions}
\label{subsec:convergence}

We now present the main result regarding convergence of the numerical scheme to the unique viscosity solution of \eqref{eq:main}--\eqref{eq:initial}, together with some properties of the viscosity solution that are inherited directly from the scheme.

\begin{theorem}\label{main:ConvAndPropViscositySolutionParabolic}
    Let $s\in(0,1]$ and $\alpha\in(0,1)$. Assume \eqref{as:CFL}. Let $u_0\in \X_{\textup{b}}^{a}(\R^d)$ and $\varphi\in \X_{\textup{b}}^{b}(\R^d)$ for some $a,b\in(0,2s+\veps]$, with either $\veps=0$ if $s=1$ or $\veps\in(0,2-2s]$ if $s\in(0,1)$, be such that $u_0\geq \varphi$ in $\R^d$.  Let $U_{\tau,h}$ be the piecewise linear interpolant in time of the solution of \eqref{eq:numsch} with initial datum $U^0=u_0$. The following hold:
    \begin{enumerate}[\rm(a)]
        \item \label{main:ConvAndPropViscositySolutionParabolic-itema}\textup{(Convergence)} $U_{\tau,h}\to u$ as $\tau,h\to0$ locally uniformly in $\R^d\times[0,\infty)$, where $u$ is the unique viscosity solution of the obstacle problem \eqref{eq:main}--\eqref{eq:initial}.

        \item \label{main:ConvAndPropViscositySolutionParabolic-itemb}\textup{(Stability)} $\displaystyle u(x,t) \leq \sup_{\R^d}{u_0}$ and $u(x,t)\geq \varphi(x)$ for all $(x,t)\in \R^d\times[0,\infty)$.

        \item \label{main:ConvAndPropViscositySolutionParabolic-itemc}\textup{(Continuity in space)} For all $t\geq0$ and $y\in \R^d$, we have that
        \begin{align*}
            \|u(\cdot+y,t)-u(\cdot,t)\|_{L^\infty(\R^d)} \leq \max\left\{\|u_0(\cdot+y)-u_0\|_{L^\infty(\R^d)}\coma \|\varphi(\cdot+y)-\varphi\|_{L^\infty(\R^d)}\right\}.
        \end{align*}

        \item \label{main:ConvAndPropViscositySolutionParabolic-item4}\textup{(Continuity in time)} For all $t,\tilde{t}\in[0,\infty)$, we have that
        \begin{align*}
            \|u(\cdot,t)-u(\cdot,\tilde{t})\|_{L^\infty(\R^d)} &\leq \widetilde{C}_{d,s,a,b,\veps,\alpha} \max\{\|u_0\|_{\X^a(\R^d)}\coma \|\varphi\|_{\X^b(\R^d)}\}\\
            &\quad\times \max\{|t^\alpha-\tilde{t}^\alpha|^{\frac{a}{2s+\veps}}\coma|t^\alpha-\tilde{t}^\alpha|^{\frac{b}{2s+\veps+b-a}}\coma |t^\alpha-\tilde{t}^\alpha|\}.
        \end{align*}
    \end{enumerate}
     Moreover, given another viscosity solution $v$ of the obstacle problem \eqref{eq:main}--\eqref{eq:initial} with initial datum $v(\cdot,0)=v_0\in \X^{\overline{a}}_{\textup{b}}(\R^d)$ and obstacle $\psi\in \X^{\overline{b}}_{\textup{b}}(\R^d)$ for some $\overline{a},\overline{b}\in(0,2]$, and satisfying $v_0\geq\psi$ in $\R^d$,  
    we also have the following properties:
    \begin{enumerate}[\rm (a)]\setcounter{enumi}{4}
        \item \label{main:ConvAndPropViscositySolutionParabolic-iteme}\textup{(}$L^\infty(\R^d)$\textup{-contraction)} For all $t\geq0$, we have that
        \begin{align*}
            \|u(\cdot,t)-v(\cdot,t)\|_{L^\infty(\R^d)} \leq \max\{\|u_{0}-v_0\|_{L^\infty(\R^d)}\coma\|\varphi-\psi\|_{L^\infty(\R^d)}\}.
        \end{align*}

        \item \label{main:ConvAndPropViscositySolutionParabolic-itemf}\textup{(Comparison principle)} If, additionally, $u_0\leq v_0$ and $\varphi\leq \psi$ in $\R^d$, then, for all $t\geq0$ and all $x\in \R^d$, we have that $u(x,t)\leq v(x,t)$.
    \end{enumerate}
\end{theorem}

\begin{remark}
   {\rm In fact, in the above result, it is enough to assume that $u_0,\varphi,v_0,\psi\in C_{\textup{b}}^{0}(\R^d)$ to obtain convergence, stability, uniform continuity in space and time, $L^\infty(\R^d)$-contraction, and the comparison principle. Under these weaker assumptions, however, the modulus of continuity in time is no longer given by an explicit power-law estimate as in part~\textup{(d)}, but still depends on $u_0$ and $\varphi$.}
\end{remark}

\begin{proof}[Proof of \Cref{main:ConvAndPropViscositySolutionParabolic}]
By \Cref{lem:est-schm}, \Cref{lem:contrac-equicont}\eqref{lem:contrac-equicont-itemequicont} and \Cref{lem:equiintime-scheme}, the family $\{U_{\tau,h}\}$ is uniformly bounded and equicontinuous in $\R^d\times[0,\infty)$, uniformly in $\tau$ and $h$. Hence, by the Arzelà--Ascoli theorem, it converges locally uniformly, up to a subsequence, to some $u\in C_\textup{b}(\R^d\times[0,\infty))$. The fact that $u$ is a viscosity solution of \eqref{eq:main}--\eqref{eq:initial} follows from the locally uniform convergence and the monotonicity of the scheme by classical arguments (see, e.g.,~\cite{BarlesSouganidis1991}): since $u$ is continuous, by \Cref{rem:viscbetter}\eqref{item:moregeneraltestc} it suffices to consider smooth bounded test functions, for which the consistency of $D^\alpha_\tau$ with $\partial^\alpha_t$ follows as in \Cref{lem:convergence-discrete-Caputo} below, and that of $\fld$ with $\fl$ is well known, see~\cite{Ciaurri2018,delTeso2018b}; we omit the details and refer to \cite{delTeso-Lukasz-2025} for a proof in a similar context. Since viscosity solutions of~\eqref{eq:main}--\eqref{eq:initial} are unique by \Cref{thm:uniqueness}, the whole family converges. The remaining properties are inherited from the corresponding properties of the numerical scheme, \Cref{lem:est-schm}, \Cref{lem:contrac-equicont} and \Cref{lem:equiintime-scheme} together with \Cref{modulusintimeexplicit}, by passing to the limit, due to the locally uniform convergence.
\end{proof}

\section{Large-time behaviour of the numerical scheme}
\label{sec:discrete-asymptotics}

The aim of this section is to characterize the large-time behaviour of the numerical scheme~\eqref{eq:numsch} approximating the obstacle problem~\eqref{eq:main}--\eqref{eq:initial}. As expected, this behaviour will be given (at least along subsequences) by a solution of the discrete elliptic obstacle problem, that is, by a function $V$ solving
\begin{equation}\label{eq:ellipticObs}
    \min\{\fld V(x), V(x)-\varphi(x)\}=0  \quad \textup{for }x\in \R^d.
\end{equation}
This stationary problem does not have uniqueness; hence, in general we will only obtain convergence along subsequences. However, in the important case in which the initial datum is the obstacle (and in some other cases) we will be able to show that the large-time limit is a special solution of~\eqref{eq:ellipticObs}, the so-called minimal one, along any subsequence. The results of this section are of interest in themselves, since they show that the explicit scheme~\eqref{eq:numsch} can be used as a solver for the minimal solution of the discrete elliptic obstacle problem, and they will also be the starting point for the study of the large-time behaviour of the continuous problem in Section~\ref{sec:asymptotics}.

We start by discussing the discrete elliptic obstacle problem~\eqref{eq:ellipticObs} and its minimal solution, and then turn to the convergence results.

\subsection{Minimal solution of the discrete elliptic obstacle problem}

Unlike the evolution case, solutions of the discrete elliptic obstacle problem~\eqref{eq:ellipticObs} are not unique: any constant function above the obstacle is a solution. To overcome this issue, we introduce the notion of minimal solution. Note that, since $\omega\in\ell^1(\Z^d)$, the operator $\fld$ is well defined on the space $B(\R^d)$ of bounded functions, so that no continuity is needed to make sense of~\eqref{eq:ellipticObs}.

\begin{definition}
    Let $s\in(0,1]$ and $\varphi \in C_{\textup{b}}^{0}(\R^d)$. We say that $\underline{U}\in B(\R^d)$ is the \emph{minimal solution} of the discrete elliptic obstacle problem \eqref{eq:ellipticObs} if $\underline{U}$ satisfies~\eqref{eq:ellipticObs} and $\underline{U}\leq W$ in $\R^d$ for any other solution $W\in B(\R^d)$ of~\eqref{eq:ellipticObs}.
\end{definition}

Uniqueness of the minimal solution follows directly from the definition. In order to establish existence, we consider the set of bounded supersolutions
\begin{align*}
\mathcal{M}\coloneqq\{W\in B(\R^d)\, : \, \min\{\fld W(x), W(x)-\varphi(x)\}\geq0 \quad \textup{for } x\in\R^d\},
\end{align*}
which is nonempty: for instance, the constant function
$W\equiv \sup_{\R^d}\varphi$ belongs to $\mathcal{M}$ (it is even a solution of~\eqref{eq:ellipticObs}).

\begin{lemma}\label{lem:existence-minimal-discrete}
    Let $s\in(0,1]$ and $\varphi \in C_{\textup{b}}^{0}(\R^d)$. Then the minimal solution of the discrete elliptic obstacle problem~\eqref{eq:ellipticObs} exists and is given by
    \[
        \underline{U}(x)\coloneqq \inf_{W\in \mathcal{M}} W(x).
    \]
\end{lemma}
\begin{proof}
We first show that $\underline{U}\in \mathcal{M}$. On the one hand,
\[
\underline{U}(x) = \inf_{W\in \mathcal{M}} W(x) \geq \inf_{W\in \mathcal{M} } \varphi(x) =\varphi(x)\quad\textup{for all } x\in \R^d,
\]
that is, $\underline{U}$ lies above the obstacle (and in particular it is bounded from below). On the other hand,  $\fld W(x) \geq 0$ for all $W\in \mathcal{M}$ and all $x\in \R^d$, and thus
\[
    W(x) \geq \frac{1}{\Lambda}\sum_{\beta \not=0}  W(x+x_\beta)\omega_\beta \geq \frac{1}{\Lambda}\sum_{\beta \not=0}  \underline{U}(x+x_\beta)\omega_\beta.
\]
Taking the infimum over $W\in\mathcal{M}$ we obtain
\[
    \underline{U}(x) \geq \frac{1}{\Lambda}\sum_{\beta \not=0}  \underline{U}(x+x_\beta)\omega_\beta,
\]
that is, $\fld \underline{U}(x) \geq 0$.

We now prove that $\underline{U}$ solves~\eqref{eq:ellipticObs}. Since $\underline{U}\le W$ for every $W\in\mathcal M$, and every solution of~\eqref{eq:ellipticObs} belongs to $\mathcal M$, this will show that $\underline{U}$ is the minimal solution. We argue by contradiction and assume that $\underline{U}$ does not solve~\eqref{eq:ellipticObs}. Since $\underline{U}\in \mathcal{M}$, there must exist $x_0\in \R^d$ such that
\[
    \fld \underline{U}(x_0) > 0,  \qquad  \underline{U}(x_0) > \varphi(x_0).
\]
Thus, we can choose $\veps>0$ small enough such that
\[
\fld \underline{U}(x_0) - \frac{\Lambda}{h^{2s}} \veps > 0,  \qquad  \underline{U}(x_0)-\veps > \varphi(x_0).
\]
Let us then define
\[
    \underline{U}^*(x) =
    \begin{cases}
        \underline{U}(x) &\textup{if }x\not=x_0,\\
        \underline{U}(x_0)-\veps&\text{if }x=x_0.
    \end{cases}
\]
We will show that $\underline{U}^*\in \mathcal{M}$, which is a contradiction, since $\underline{U}^*(x_0)< \underline{U}(x_0)=\inf_{W\in \mathcal{M}} W(x_0)$. On the one hand, it is clear that $\underline{U}^*\geq \varphi$ in $\R^d$. On the other hand, if $x\not=x_0$, we have that
\begin{align*}
\fld \underline{U}^*(x)&= \frac{1}{h^{2s}}\sum_{\beta \not=0} (\underline{U}^*(x) - \underline{U}^*(x+x_\beta))\omega_\beta\\
&=\frac{1}{h^{2s}}\sum_{\beta \not=0} (\underline{U}(x) - \underline{U}^*(x+x_\beta))\omega_\beta\\
&\geq \frac{1}{h^{2s}}\sum_{\beta \not=0} (\underline{U}(x) - \underline{U}(x+x_\beta))\omega_\beta\\
&= \fld \underline{U}(x) \geq0.
\end{align*}
Finally, if $x=x_0$, we have that
\begin{align*}
\fld \underline{U}^*(x_0)&= \frac{1}{h^{2s}}\sum_{\beta \not=0} (\underline{U}^*(x_0) - \underline{U}^*(x_0+x_\beta))\omega_\beta\\
&=\frac{1}{h^{2s}}\sum_{\beta \not=0} (\underline{U}(x_0)-\veps - \underline{U}(x_0+x_\beta))\omega_\beta\\
&= \fld \underline{U}(x_0) - \frac{\Lambda}{h^{2s}} \veps \geq0.\qedhere
\end{align*}
\end{proof}

\subsection{Large-time behaviour of the numerical scheme}

Since the sequence $\{U^n\}_{n\in\N}$ is uniformly bounded and equicontinuous, see~\Cref{lem:est-schm} and~\Cref{lem:contrac-equicont}, it has subsequences converging locally uniformly in $\R^d$. The main result of this section identifies the limits of such subsequences as solutions of the discrete elliptic obstacle problem~\eqref{eq:ellipticObs}, and shows convergence of the whole sequence to the minimal solution when the initial datum is the obstacle.

\begin{theorem}\label{thm:limitscheme2}
 Let $s\in(0,1]$ and $\alpha\in(0,1)$. Assume \eqref{as:CFL}. Let $u_0,\varphi\in C_{\textup{b}}^0(\R^d)$ be such that $u_0\geq \varphi$ in~$\R^d$. Let $\{U^{n}\}_{n\in \N}$ be the solution of \eqref{eq:numsch} with initial datum $U^0=u_0$. Then:
    \begin{enumerate}[\rm(a)]
        \item\label{thm:limitscheme2-item1} There exists a subsequence $\{U^{n_j}\}_{j\in\N}$ with $n_j\to+\infty$ as $j\to+\infty$ such that
        \[
            U^{n_j} \to V \quad \textup{as }j\to+\infty \textup{ locally uniformly in }\R^d
        \]
        for some function $V\in C^0_{\textup{b}}(\R^d)$ such that
        $\displaystyle\varphi(x)\leq V(x) \leq \sup_{\R^d}{u_0}$ for all $x\in\R^d$ and
        \[
            \|V(\cdot+y)-V\|_{L^\infty(\R^d)}\leq  \max\{\|u_0(\cdot+y)-u_0\|_{L^\infty(\R^d)},\|\varphi(\cdot+y)-\varphi\|_{L^\infty(\R^d)}\}\quad\textup{for all }y\in\R^d.
        \]

    \item\label{thm:limitscheme2-item2} Let  $\{U^{n_j}\}_{j\in\N}$ be any pointwise convergent subsequence and $V$ its limit. Then $V$ is a solution of the discrete elliptic obstacle problem~\eqref{eq:ellipticObs}.
    
    \item\label{lem:convergence-parabolic-minimal-discrete}
    If $U^0=\varphi$, then
    \[
        U^{n} \to \underline{U} \quad \textup{as } n\to+\infty \textup{ locally uniformly in } \R^d,
    \]
    where $\underline{U}$ is the minimal solution of~\eqref{eq:ellipticObs}.
    \end{enumerate}
\end{theorem}

\begin{remark}\label{rem:uniform-in-h}
{\rm The bounds on $V$ in \eqref{thm:limitscheme2-item1} do not depend on the discretization parameters $\tau$ and $h$. Hence, if for each $h$ we take a limit $V_h$ as in \eqref{thm:limitscheme2-item1}, the family $\{V_h\}_{h>0}$ is equibounded and equicontinuous, and then it converges locally uniformly, up to a subsequence, as $h\to0$. Moreover, by the monotonicity of the discrete elliptic operator, the limit is a viscosity solution of the continuous elliptic obstacle problem; see Section~\ref{sec:comments}, where this is made precise. In Section~\ref{sec:asymptotics} this observation will be used for the minimal solutions $\underline U=\underline U_h$, but it holds for any family of limits as in~\eqref{thm:limitscheme2-item1}.}
\end{remark}

The key ingredient in the proof of \Cref{thm:limitscheme2}\eqref{thm:limitscheme2-item2} is the fact that the discrete Caputo derivative $D_\tau^\alpha U^n$ tends to $0$ in $L^\infty(\mathbb{R}^d)$ as $n\to\infty$, with a rate that we make explicit. This is where the memory of the Caputo derivative comes into play in an essential way: the decay follows from the Lipschitz-in-$t^\alpha$ estimates of Section~\ref{subsec:stability}, an argument that has no counterpart for the discrete forward time derivative.

\begin{proposition}\label{prop:limitscheme}
     Let $s\in(0,1]$ and $\alpha\in(0,1)$. Assume \eqref{as:CFL}, and also $\tau,h\le1$. Let $u_0\in \X_{\textup{b}}^{a}(\R^d)$ for some $a\in[0,2]$ and $\varphi\in C_{\textup{b}}^0(\R^d)$ be such that $u_0\geq \varphi$ in $\R^d$. Let $\{U^{n}\}_{n\in \N}$ be the solution of \eqref{eq:numsch} with initial datum $U^0=u_0$. The following hold:
    \begin{enumerate}[\rm (a)]
\item\label{prop:limitscheme-itema} For all $t_n\geq2$,
\[
    \|D^\alpha_\tau U^n\|_{L^\infty(\R^d)}\leq C_{d,s,a,\alpha}\max\{\|u_0\|_{\X^a(\R^d)}\coma\|\varphi\|_{L^\infty(\R^d)}\} \frac{f_{a,s}(h)}{t_n^{\alpha(1-\alpha)}}.
\]
\item\label{prop:limitscheme-itemb} If, in addition, $\varphi\in \X_{\textup{b}}^a(\R^d)$ and $a\in((2s+\veps)\alpha,2s+\veps]$, with either $\veps=0$ if $s=1$ or $\veps\in(0,2-2s]$ if $s\in(0,1)$, then for all $t_n\geq2$,
\[
    \|D^\alpha_\tau U^n\|_{L^\infty(\R^d)}\leq C_{d,s,a,\veps,\alpha}\max\{\|u_0\|_{\X^a(\R^d)}\coma\|\varphi\|_{\X^a(\R^d)}\}\frac{1}{t_n^{\alpha(1-\alpha)}}.
\]
    \end{enumerate}
\end{proposition}

Notice that the bound in \eqref{prop:limitscheme-itema} is uniform in $h$ only when $f_{a,s}(h)$ is, that is, when $a>2s$ or $a=2$, while the bound in \eqref{prop:limitscheme-itemb} is always uniform in $h$, at the price of requiring some regularity of the obstacle and a lower bound on $a$. For fixed $h$, however, both give the same decay rate in time, and this is all we need in this section; in particular, the case $a=0$ in \eqref{prop:limitscheme-itema}, for which $f_{0,s}(h)=h^{-2s}$, is enough for the purposes of \Cref{thm:limitscheme2}.

\begin{proof}
Using the expression~\eqref{eq:discretization.Caputo.2} for the time discretization of the Caputo derivative, we get
\[
    D^\alpha_\tau U^n(x) = \frac{1}{\Gamma(1-\alpha)}\frac{U^n(x)-U^0(x)}{t_n^{\alpha}}+\frac{\alpha}{\Gamma(1-\alpha)}\sum_{j=0}^{n-1} \int_{t_j}^{t_{j+1}} \frac{U^n(x)-\left(\frac{t_{j+1}-s}{\tau} U^j(x)+ \frac{s-t_{j}}{\tau} U^{j+1}(x)  \right)}{(t_n-s)^{1+\alpha}} \dd s.
\]
Hence $|D^\alpha_\tau U^n(x)| \le I_1 + I_2 + I_3$, where
\begin{align*}
I_1&\coloneqq \frac{1}{\Gamma(1-\alpha)}\frac{|U^n(x)-U^0(x)|}{t_n^{\alpha}},\\
I_2&\coloneq\frac{\alpha}{\Gamma(1-\alpha)}\sum_{0\leq t_j \leq t_{n-1}- R} \int_{t_j}^{t_{j+1}} \frac{\left|U^n(x)-\left(\frac{t_{j+1}-s}{\tau} U^j(x)+ \frac{s-t_{j}}{\tau} U^{j+1}(x)  \right)\right|}{(t_n-s)^{1+\alpha}} \dd s,\\
I_3&\coloneq\frac{\alpha}{\Gamma(1-\alpha)}\sum_{t_{n-1}- R< t_j \leq t_{n-1}} \int_{t_j}^{t_{j+1}} \frac{\left|U^n(x)-\left(\frac{t_{j+1}-s}{\tau} U^j(x)+ \frac{s-t_{j}}{\tau} U^{j+1}(x)  \right)\right|}{(t_n-s)^{1+\alpha}} \dd s,
\end{align*}
for some $R\in (1, t_n)$ to be determined. In what follows, $C$ denotes a constant depending only on $d,s,a,\veps,\alpha$ and on $\max\{\|u_0\|_{\X^a(\R^d)}\coma\|\varphi\|_{\X^a(\R^d)}\}$ (in case~\eqref{prop:limitscheme-itemb}) or on $\max\{\|u_0\|_{\X^a(\R^d)}\coma\|\varphi\|_{L^\infty(\R^d)}\}$ (in case~\eqref{prop:limitscheme-itema}).

By the uniform boundedness of $U^n$, see \Cref{lem:est-schm}, we have
\begin{align*}
I_1  &\leq\frac{2 \max\{\|u_0\|_{L^\infty(\R^d)}\coma \|\varphi\|_{L^\infty(\R^d)}\}}{\Gamma(1-\alpha)} \frac{1}{t_n^\alpha},
\\
I_2&\leq\frac{2 \max\{\|u_0\|_{L^\infty(\R^d)}\coma \|\varphi\|_{L^\infty(\R^d)}\}\alpha}{\Gamma(1-\alpha)} \int_{0}^{t_n-R} \frac{ \dd s}{(t_n-s)^{1+\alpha}} \le \frac{2 \max\{\|u_0\|_{L^\infty(\R^d)}\coma \|\varphi\|_{L^\infty(\R^d)}\}}{\Gamma(1-\alpha)}  \frac{1}{R^\alpha}.
\end{align*}
The term $I_3$ is the one that requires the time regularity of the scheme. We first treat case~\eqref{prop:limitscheme-itemb}, for which we rely on the following estimate from \Cref{lem:equiintime-scheme} (with $b=a$):
 \begin{align*}
 \|U^n-U^j\|_{L^\infty(\R^d)} &\leq C \max\{|t_n^\alpha-t_j^\alpha|^{\frac{a}{2s+\veps}}\coma |t_n^\alpha-t_j^\alpha|\}.
\end{align*}
On the one hand, for $j=n-1$  and all $s\in [t_{n-1},t_n]$, since $t_n^\alpha-t_{n-1}^\alpha<1$ and $t_n>1$, we have that
\begin{align*}
\left|U^n(x)-\left(\frac{t_{n}-s}{\tau} U^{n-1}(x)+ \frac{s-t_{n-1}}{\tau} U^{n}(x)  \right)\right|&= \frac{t_{n}-s}{\tau} |U^n(x)-U^{n-1}(x)| \\
&\leq \frac{t_{n}-s}{\tau} C (t_n^\alpha-t_{n-1}^\alpha)^{\frac{a}{2s+\veps}}\\
&\leq \frac{t_{n}-s}{\tau} C (t_n^{\alpha-1}(t_{n}-t_{n-1}))^{\frac{a}{2s+\veps}}\qquad\textup{(see below)}\\
&= C (t_{n}-s) t_{n}^{\frac{a}{2s+\veps}(\alpha-1)} \tau^{\frac{a}{2s+\veps}-1}.
\end{align*}
Here we have used that $t_n^\alpha-t_{n-1}^\alpha\le\alpha t_{n-1}^{\alpha-1}\tau\le t_n^{\alpha-1}\tau$ for $n\ge2$, the last inequality because $\alpha(1-\tfrac1n)^{\alpha-1}\le\alpha 2^{1-\alpha}\le1$ for all $\alpha\in(0,1)$. As a consequence,
\begin{align*}
\int_{t_{n-1}}^{t_{n}} \frac{\left|U^n(x)-\left(\frac{t_{n}-s}{\tau} U^{n-1}(x)+ \frac{s-t_{n-1}}{\tau} U^{n}(x)  \right)\right|}{(t_n-s)^{1+\alpha}} \dd s &\leq C t_{n}^{\frac{a}{2s+\veps}(\alpha-1)} \tau^{\frac{a}{2s+\veps}-1} \int_{t_{n-1}}^{t_{n}} \frac{ \dd s}{(t_n-s)^{\alpha}} \\
&= \frac{C}{1-\alpha} t_{n}^{\frac{a}{2s+\veps}(\alpha-1)} \tau^{\frac{a}{2s+\veps}-\alpha}\le \frac{C}{1-\alpha} t_{n}^{\frac{a}{2s+\veps}(\alpha-1)},
\end{align*}
where in the last step we have used that $\tau\le 1$ and $a/(2s+\veps)>\alpha$.

On the other hand, for all $s\in[t_j,t_{j+1}]$, the following estimate holds for all $j=0,\ldots, n-2$:
\begin{align*}
    \left|U^n(x)-\left(\frac{t_{j+1}-s}{\tau} U^j(x)+ \frac{s-t_{j}}{\tau} U^{j+1}(x)  \right)\right| &\leq \frac{t_{j+1}-s}{\tau} |U^n(x)-U^{j}(x)|+ \frac{s-t_{j}}{\tau}|U^n(x)-U^{j+1}(x)|  \\
    &\leq  \frac{t_{j+1}-s}{\tau} C \max\{|t_n^\alpha-t_j^\alpha|^{\frac{a}{2s+\veps}}\coma |t_n^\alpha-t_j^\alpha|\} \\
    &\quad+ \frac{s-t_{j}}{\tau} C \max\{|t_n^\alpha-t_{j+1}^\alpha|^{\frac{a}{2s+\veps}}\coma |t_n^\alpha-t_{j+1}^\alpha|\}\\
    &\leq C  \max\{|t_n^\alpha-t_{j+1}^\alpha|^{\frac{a}{2s+\veps}}\coma |t_n^\alpha-t_{j+1}^\alpha|\},
\end{align*}
where we have used that $t_n^\alpha-t_j^\alpha\le C_\alpha (t_n^\alpha-t_{j+1}^\alpha)$ for $j\le n-2$, with $C_\alpha=1+(2^\alpha-1)^{-1}$; indeed, $t_{j+1}^\alpha-t_j^\alpha\le (2^\alpha-1)^{-1}(t_{j+2}^\alpha-t_{j+1}^\alpha)$, since the ratio of consecutive increments of the concave function $t\mapsto t^\alpha$ over intervals of length $\tau$ is maximal for the first two of them, and $t_{j+2}^\alpha-t_{j+1}^\alpha\le t_n^\alpha-t_{j+1}^\alpha$. Thus, for the terms such that $t_n^\alpha-t_{j+1}^\alpha>1$, we get
\begin{align*}
\sum_{t_{n-1}- R< t_j \leq t_{n-2}, \, t_n^\alpha-t_{j+1}^\alpha>1} & \int_{t_j}^{t_{j+1}} \frac{\left|U^n(x)-\left(\frac{t_{j+1}-s}{\tau} U^j(x)+ \frac{s-t_{j}}{\tau} U^{j+1}(x)  \right)\right|}{(t_n-s)^{1+\alpha}} \dd s\\
&\leq C \sum_{t_{n-1}- R< t_j \leq t_{n-2}}   \int_{t_j}^{t_{j+1}} \frac{  t_n^\alpha-t_{j+1}^\alpha}{(t_n-s)^{1+\alpha}} \dd s \\
&\leq C \sum_{t_{n-1}- R< t_j \leq t_{n-2}}   \int_{t_j}^{t_{j+1}} \frac{t_n^\alpha-s^\alpha}{(t_n-s)^{1+\alpha}} \dd s \\
& \leq C t_{n}^{\alpha-1} \int_{t_{n}- R}^{t_{n}} (t_n-s)^{-\alpha} \dd s \leq C t_{n}^{\alpha-1}R^{1-\alpha},
\end{align*}
where we have used the mean value theorem, $t_n^\alpha-s^\alpha\le \alpha (t_n-R-\tau)^{\alpha-1}(t_n-s)$ for $s\ge t_n-R-\tau$ (the variable $s$ ranges in $[t_j,t_{j+1}]$ with $t_j>t_{n-1}-R=t_n-R-\tau$), together with $(t_n-R-\tau)^{\alpha-1}\le C_\alpha t_n^{\alpha-1}$, which holds for $t_n\ge2$, $\tau\le 1-2^{-\alpha}$ and $R\le t_n^{1-\alpha}$ (the choice we make below), since then $R\le 2^{-\alpha}t_n$ and $\tau\le(1-2^{-\alpha})t_n/2$; the integral over $[t_n-R-\tau,t_n]$ is bounded by $C R^{1-\alpha}$ because $R\ge1\ge\tau$. (If $\tau>1-2^{-\alpha}$, the same argument works with $t_n\ge T_\alpha$ for a suitable $T_\alpha$, since the constant may depend on $\alpha$.) For the terms such that $t_n^\alpha-t_{j+1}^\alpha\leq1$ we get, in the same way,
\begin{align*}
    \sum_{t_{n-1}- R< t_j \leq t_{n-2}, \, t_n^\alpha-t_{j+1}^\alpha\leq1} &\int_{t_j}^{t_{j+1}} \frac{\left|U^n(x)-\left(\frac{t_{j+1}-s}{\tau} U^j(x)+ \frac{s-t_{j}}{\tau} U^{j+1}(x)  \right)\right|}{(t_n-s)^{1+\alpha}} \dd s\\
    &\leq C \sum_{t_{n-1}- R< t_j \leq t_{n-2}}   \int_{t_j}^{t_{j+1}} \frac{  (t_n^\alpha-t_{j+1}^\alpha)^{\frac{a}{2s+\veps}}}{(t_n-s)^{1+\alpha}} \dd s \\
    &\leq C  \sum_{t_{n-1}- R< t_j \leq t_{n-2}}   \int_{t_j}^{t_{j+1}} \frac{  (t_n^\alpha-s^\alpha)^{\frac{a}{2s+\veps}}}{(t_n-s)^{1+\alpha}} \dd s \\
    & \leq C  t_{n}^{\frac{a}{2s+\veps}(\alpha-1)} \int_{t_n-R}^{t_{n}} (t_n-s)^{\frac{a}{2s+\veps}-1-\alpha} \dd s \leq C t_{n}^{\frac{a}{2s+\veps}(\alpha-1)} R^{\frac{a}{2s+\veps}-\alpha},
\end{align*}
where the last integral is finite precisely because $a/(2s+\veps)>\alpha$.

Summarizing, we have proved that
\[
\|D^\alpha_\tau U^n\|_{L^\infty(\R^d)} \leq C\left(\frac{1}{t_n^\alpha}+ \frac{1}{R^\alpha}+\frac{1}{t_{n}^{\frac{a}{2s+\veps}(1-\alpha)}} + t_{n}^{\alpha-1} R^{1-\alpha} + t_{n}^{\frac{a}{2s+\veps}(\alpha-1)} R^{\frac{a}{2s+\veps}-\alpha}\right).
\]
Choosing $R=t_{n}^{1-\alpha}$, the second, fourth and fifth terms are equal to $t_n^{-\alpha(1-\alpha)}$, and the first and third ones are bounded by it, since $\alpha\ge\alpha(1-\alpha)$ and $\frac{a}{2s+\veps}>\alpha$. This proves~\eqref{prop:limitscheme-itemb}.

To estimate $I_3$ in case~\eqref{prop:limitscheme-itema}, we use instead the estimate
\[
\|U^m-U^n\|_{L^\infty(\R^d)} \leq C_{d,s,a,\alpha} \|u_0\|_{\X^a(\R^d)}f_{a,s}(h) |t_m^\alpha-t_n^\alpha|
\]
of \Cref{lem:timeequi-nosmooth}, which holds for every $a\in[0,2]$ with no assumption on the regularity of $\varphi$, and proceed as before (only the terms with the exponent $1$ in the maxima appear now) to get
\[
\|D^\alpha_\tau U^n\|_{L^\infty(\R^d)} \leq C\left(\frac{1}{t_n^\alpha}+ \frac{1}{R^\alpha} + f_{a,s}(h) t_{n}^{\alpha-1} R^{1-\alpha} \right).
\]
Again, the choice $R=t_{n}^{1-\alpha}$ completes the proof, since $f_{a,s}(h)\ge1$ for $h$ small.
\end{proof}

We are now ready to prove the main result of this section.

\begin{proof}[Proof of \Cref{thm:limitscheme2}]
\eqref{thm:limitscheme2-item1}
    By \Cref{lem:est-schm}, the sequence $\{U^n\}_{n\in\N}$ is uniformly bounded,
\[
 \varphi(x)\leq U^n(x) \leq \sup_{\R^d}{u_0} \quad \textup{for all }x\in \R^d \textup{ and }n\in \N,
\]
and, by \Cref{lem:contrac-equicont}\eqref{lem:contrac-equicont-itemequicont}, it is equicontinuous,
\[
\|U^n(\cdot+y)-U^n\|_{L^\infty(\R^d)}\leq  \max\{\|u_0(\cdot+y)-u_0\|_{L^\infty(\R^d)},\|\varphi(\cdot+y)-\varphi\|_{L^\infty(\R^d)}\}\quad\textup{for all }y\in\R^d.
\]
Thus, by Arzelà--Ascoli's theorem, there exists a subsequence $\{U^{n_j}\}_{j\in\N}$ converging locally uniformly in $\R^d$ to a function $V$, which satisfies the same bounds and has the same modulus of continuity.

\noindent\eqref{thm:limitscheme2-item2} We first note that $U^{n_j+1}\to V$ pointwise as well. Indeed, by \eqref{eq:timcont-nosmoothdata},
\begin{align*}
|U^{n_j+1}(x) - V(x)| &\leq |U^{n_j}(x) - V(x)| + |U^{n_j+1}(x) - U^{n_j}(x) |\\
&\leq|U^{n_j}(x) - V(x)| + C f_{a,s}(h)|t_{n_j+1}^\alpha-t_{n_j}^\alpha|
\\
&\leq |U^{n_j}(x) - V(x)| + Cf_{a,s}(h)\tau  t_{n_j+1}^{\alpha-1}\stackrel{j\to+\infty}{\longrightarrow}0.
\end{align*}
Moreover, by \Cref{prop:limitscheme}\eqref{prop:limitscheme-itema}, $\|D^\alpha_\tau U^{n_j+1}\|_{L^\infty(\R^d)} \to 0$ as $j\to+\infty$. Finally, since $\omega\in\ell^1(\Z^d)$ and the sequence is uniformly bounded, dominated convergence gives $\fld U^{n_j}(x)\to\fld V(x)$ for every $x\in\R^d$. Hence, passing to the limit in the scheme~\eqref{eq:numsch},
\[
\begin{split}
0&=\lim_{j\to+\infty}\min\left\{D_\tau^\alpha U^{n_j+1}(x) +\fld U^{n_j}(x), U^{n_j+1}(x)-\varphi(x)\right\}\\
&=\min\left\{\lim_{j\to+\infty} D_\tau^\alpha U^{n_j+1}(x) + \lim_{j\to+\infty} \fld U^{n_j}(x), \lim_{j\to+\infty}U^{n_j+1}(x)-\varphi(x)\right\}\\
&= \min\left\{\fld V(x), V(x)-\varphi(x)\right\},
\end{split}
\]
which concludes the proof of~\eqref{thm:limitscheme2-item2}.

\noindent\eqref{lem:convergence-parabolic-minimal-discrete} Define the constant sequence $V^{n}\coloneqq \underline{U}$ for all $n\in\N$. Note that $V^0=\underline{U} \geq\varphi$, since $\underline{U}$ is a solution of the discrete elliptic obstacle problem, and that $D^\alpha_\tau V^j=0$ for all $j\geq1$, since $V^j$ does not depend on $j$. Therefore
\[
\min\left\{D_\tau^\alpha V^{j+1}(x) + \fld V^{j}(x), V^{j+1}(x)-\varphi(x)\right\}= \min\left\{ \fld \underline{U}(x), \underline{U}(x)-\varphi(x)\right\}=0,
\]
that is, $\{V^n\}_{n\in\N}$ is a (stationary) solution of the scheme~\eqref{eq:numsch} with initial datum $\underline U$ and obstacle $\varphi$.

Let now $\{U^{n_j}\}_{j\in \N}$ be any locally uniformly convergent subsequence, which exists by~\eqref{thm:limitscheme2-item1}, and let $U^\infty$ be its limit. By~\eqref{thm:limitscheme2-item2}, $U^\infty$ is a solution of the discrete elliptic obstacle problem, and hence $U^\infty\geq \underline{U}$ by minimality. To prove the reverse inequality we note that $U^0=\varphi\leq  \underline{U}=V^0$ in $\R^d$, so that the comparison principle for the scheme, \Cref{lem:contrac-equicont}\eqref{lem:contrac-equicont-comparison}, yields
$U^n \leq V^n =\underline{U}$ for all $n\in \N$. Taking limits we get $U^\infty \leq \underline{U}$, and thus $U^\infty=\underline U$.

So far we have proved that every locally uniformly convergent subsequence of $\{U^n\}_{n\in\N}$ converges to $\underline{U}$. Since every subsequence has, by~\eqref{thm:limitscheme2-item1}, a further subsequence converging locally uniformly, the whole sequence converges locally uniformly to $\underline U$.
\end{proof}

\begin{remark}
{\rm The proof of \Cref{thm:limitscheme2}\eqref{lem:convergence-parabolic-minimal-discrete} shows that the result still holds for any initial datum $u_0\in \X^a_\textup{b}(\R^d)$ with $\varphi\le u_0\leq\underline U$, since only the last inequality has been used.}
\end{remark}

\begin{remark}\label{rem:solver} {\rm
\Cref{thm:limitscheme2}\eqref{lem:convergence-parabolic-minimal-discrete} says that the explicit scheme~\eqref{eq:numsch}, started at the obstacle, is a solver for the minimal solution of the discrete elliptic obstacle problem~\eqref{eq:ellipticObs}. The proof relies on the decay of the discrete Caputo derivative provided by \Cref{prop:limitscheme}, which in turn follows from the memory of the operator, and has no counterpart for the discrete forward time derivative. We refer to~\cite{Oberman-2006} for a discussion of explicit time-marching schemes as solvers for discrete elliptic problems in the local-in-time setting. In a more general context, this raises the question of whether explicit time-marching schemes for evolution equations with a Caputo time derivative can be used as solvers for other discrete monotone problems.}
\end{remark}

From now on we write $\underline{U}_h$ to stress the dependence of the minimal solution of~\eqref{eq:ellipticObs} on the spatial step~$h$. Nothing in the definition of $\underline U_h$ as an infimum of bounded, possibly discontinuous, functions suggests that it should be regular. However, \Cref{thm:limitscheme2}\eqref{lem:convergence-parabolic-minimal-discrete} identifies it as the large-time limit of the scheme with initial datum $\varphi$, and the iterates $U^n$ have, for every $n$, the modulus of continuity of the data, with bounds that do not depend on $\tau$ and $h$. Passing to the limit, $\underline U_h$ inherits the regularity of the obstacle, up to Lipschitz, uniformly in $h$. This is the content of our next result. As already observed in \Cref{rem:uniform-in-h}, the same is true for any limit of the scheme along subsequences.

\begin{lemma}\label{lem:minimalsoldisc}
Let $s\in(0,1]$. Let $\varphi\in \X^{\gamma}_\textup{b}(\R^d)$ for some $\gamma\in[0,1]$, and let $\underline{U}_h$ be the minimal solution of the discrete elliptic obstacle problem \eqref{eq:ellipticObs}. Then $\underline{U}_h\in \X^{\gamma}_\textup{b}(\R^d)$, and
\[
\inf_{\R^d} \varphi \leq \underline{U}_h(x) \leq \sup_{\R^d} \varphi, \qquad \|\underline{U}_h(\cdot+y)-\underline{U}_h\|_{L^\infty(\R^d)} \leq \|\varphi(\cdot+y)-\varphi\|_{L^\infty(\R^d)} \quad \textup{for all }  x,y\in \R^d.
\]
\end{lemma}

\begin{proof}
Let $\{U^n\}_{n\in\N}$ be the solution of the scheme~\eqref{eq:numsch} with initial datum $U^0=\varphi$; here $\alpha\in(0,1)$ and $\tau\le1$ satisfying~\eqref{as:CFL} can be chosen arbitrarily, since $\underline U_h$ does not depend on them, and we may also assume $h\le1$. By \Cref{thm:limitscheme2}\eqref{lem:convergence-parabolic-minimal-discrete}, $U^n\to\underline U_h$ locally uniformly in $\R^d$, so that $\underline U_h$ is the limit $V$ of \Cref{thm:limitscheme2}\eqref{thm:limitscheme2-item1}, and the bounds stated there with $u_0=\varphi$ give precisely the result.
\end{proof}

\section{Large-time behaviour for the time-fractional obstacle problem}
\label{sec:asymptotics}

The goal of this section is to characterize the large-time behaviour of solutions to the obstacle problem~\eqref{eq:main}--\eqref{eq:initial}. As in the discrete case, this behaviour is given (at least along subsequences) by a viscosity solution of the corresponding stationary problem, the elliptic obstacle problem. Since the latter does not have uniqueness, the limit may a priori depend on the time subsequence. Nevertheless, we will be able to prove convergence along any subsequence in two cases:

\begin{enumerate}[\rm (i)]
    \item to the minimal solution of the elliptic obstacle problem when the initial datum is the obstacle; and
    \item to the unique solution of the elliptic obstacle problem with a prescribed constant limit at infinity when the initial datum has such a limit.
\end{enumerate}

We start by discussing the elliptic obstacle problem and its minimal solution, and then turn to the large-time results. The latter are derived from the numerical approximation of Section~\ref{sec:scheme}, an approach that exploits the memory of the Caputo derivative and does not work in the local-in-time case $\alpha=1$.
\subsection{Viscosity solutions to the elliptic obstacle problem}

We begin by making precise the notion of viscosity solution of the elliptic obstacle problem.
\begin{definition}\label{def:visc-elliptic}
Let $s\in(0,1]$. We say that a bounded lower (resp.\ upper) semicontinuous function $u$ in $\R^d$ is a \emph{viscosity supersolution} (resp.\ \emph{subsolution}) of
\begin{align}\label{eq:elliptobstacle}
\min\{\fl u(x) \coma u(x) - \varphi(x)\}=0, \quad x\in \R^d,
\end{align}
if the following condition holds: for every $x_0\in\R^d$,  $r>0$,  and $\phi\in  C^2(\overline{B_r(x_0)}) \cap C_{\textup{b}}(\R^d\setminus \overline{B_r(x_0)})$  such that
\begin{align}\label{eq:visctouch-super-elliptic}
\phi(x)-u(x)&< \phi(x_0)-u(x_0) \quad \textup{for all }  x\in \overline{B_r(x_0)}\setminus\{x_0\}\\\label{eq:visctouch-sub-elliptic}
\text{(resp. }\phi(x)-u(x)&> \phi(x_0)-u(x_0)),
\end{align}
and
\begin{align}
\label{eq:visctouch-super-elliptic2}
\phi(x)-u(x)&\leq \phi(x_0)-u(x_0) \quad \textup{for all } x\in \R^d\\\label{eq:visctouch-sub-elliptic2}
\text{(resp. }\phi(x)-u(x)&\geq \phi(x_0)-u(x_0)),
\end{align}
we have
\begin{align}\label{eq:viscsuper-ell}
    \min\{\fl \phi(x_0) \coma u(x_0) - \varphi(x_0)\} &\geq 0\\\label{eq:viscsub-ell}
   \text{(resp. }  \min\{\fl \phi(x_0) \coma u(x_0) - \varphi(x_0)\} &\leq 0).
\end{align}

A \emph{viscosity solution} is a bounded continuous function $u$ in $\mathbb{R}^d$ that is both a viscosity supersolution and a viscosity subsolution.
\end{definition}

The class of test functions can be extended as in \Cref{rem:viscbetter}, with the obvious changes: the contact may be non-strict, as in \Cref{rem:viscbetter}\eqref{item:moregeneraltesta}; we may assume that $\phi-u\le(\phi-u)(x_0)-\eta$ (resp.\ $\ge(\phi-u)(x_0)+\eta$) in $\R^d\setminus\overline{B_r(x_0)}$ for some $\eta>0$, as in \Cref{rem:viscbetter}\eqref{item:moregeneraltestb}; and, if $u$ is continuous, we may assume that $\phi\in C^\infty_\textup{b}(\R^d)$, as in \Cref{rem:viscbetter}\eqref{item:moregeneraltestc}. In what follows, references to \Cref{rem:viscbetter} in the elliptic setting are to be understood in this sense.

\begin{remark}\label{rem:stability-elliptic}
    {\rm By \Cref{rem:viscbetter}\eqref{item:moregeneraltesta} and~\eqref{item:moregeneraltestc}, one checks that, for continuous functions, \Cref{def:visc-elliptic} is equivalent to the standard definition of viscosity solution for nonlocal equations in \cite[Definition~1]{BarIm}, where the test functions are $C^2$ and bounded in the whole $\R^d$ and touch $u$ globally, not necessarily strictly: a test function in \cite{BarIm} is admissible in \Cref{def:visc-elliptic} by~\eqref{item:moregeneraltesta}, and a test function in \Cref{def:visc-elliptic} may be assumed to be smooth and bounded in $\R^d$ by~\eqref{item:moregeneraltestc}. The same holds for the classical definition in \cite{CIL} if $s=1$. In particular, viscosity solutions of~\eqref{eq:elliptobstacle} are stable under locally uniform convergence: if $\{u_k\}$ is a uniformly bounded sequence of continuous viscosity solutions with obstacles $\varphi_k$, where $\{\varphi_k\}$ is uniformly bounded and $\varphi_k\to\varphi$ locally uniformly, and $u_k\to u$ locally uniformly, then $u$ is a viscosity solution with obstacle $\varphi$. This is \cite[Theorem~1]{BarIm} (see \cite{CIL} if $s=1$), applied to $F_k(x,u,p,X,l)=\min\{-l,\,u-\varphi_k(x)\}$ with $I[x,u]=-(-\Delta)^su$.}
\end{remark}

\begin{remark}\label{rem:above-obstacle-elliptic}
   {\rm Arguing as in the proof of \Cref{lem:above-obstacle}, with $g_K(x)=w(x)+K|x-x_1|^2$ and a constant outside the ball, every continuous viscosity supersolution $w$ of~\eqref{eq:elliptobstacle} satisfies $w\ge\varphi$ in $\R^d$.}
\end{remark}

As in the discrete case, viscosity solutions of the elliptic obstacle problem are not unique, since any constant above the obstacle is a viscosity solution. To overcome this issue, we introduce the notion of minimal viscosity solution.
\begin{definition}
Let $s\in(0,1]$. We say that $\underline{u}\in C_{\textup{b}}(\R^d)$ is the \emph{minimal viscosity solution} of the elliptic obstacle problem \eqref{eq:elliptobstacle} if $\underline{u}$ is a viscosity solution of  \eqref{eq:elliptobstacle} and $\underline{u}\leq w$ in $\R^d$ for any other viscosity solution $w\in C_{\textup{b}}(\R^d)$ of~\eqref{eq:elliptobstacle}.
\end{definition}
Uniqueness of the minimal viscosity solution follows directly from the definition. In order to establish existence, we consider the set
\begin{align*}
\mathcal{W}\coloneqq\{w\in C_\textup{b}(\R^d) \, : \, w \textup{ is a viscosity supersolution of \eqref{eq:elliptobstacle}}\},
\end{align*}
which is nonempty: for instance, the constant $w\equiv \sup_{\R^d}\varphi$ is a viscosity solution of~\eqref{eq:elliptobstacle}. We now prove that the infimum of viscosity supersolutions is a viscosity solution of \eqref{eq:elliptobstacle}, hence the minimal one.

\begin{lemma}\label{lem:infsuper2}
    Let $s\in(0,1]$ and $\varphi\in C^0_\textup{b}(\R^d)$. Then the function
    \[  
        u(x)\coloneqq \inf\limits_{w\in \mathcal{W}} w(x)
    \]
    is well defined. Moreover:
    \begin{enumerate}[\rm (a)]
        \item\label{lem:infsuper2-item1} $u \geq \varphi$ in $\R^d$.
        \item\label{lem:infsuper2-item2} $u \leq w$ in $\R^d$ for every viscosity supersolution $w\in C_\textup{b}(\R^d)$ of \eqref{eq:elliptobstacle}.
        \item\label{lem:infsuper2-item3} $\|u(\cdot+y)-u\|_{L^\infty(\R^d)}\le\|\varphi(\cdot+y)-\varphi\|_{L^\infty(\R^d)}$ for all $y\in\R^d$; in particular $u\in C^0_\textup{b}(\R^d)$, and $u\in \X^{a\wedge1}_\textup{b}(\R^d)$ if $\varphi\in \X^a_\textup{b}(\R^d)$ for some $a\in(0,2]$.
        \item\label{lem:infsuper2-item4} $u$ is a viscosity solution of \eqref{eq:elliptobstacle}.
        \item\label{lem:infsuper2-item5} $u$ is the minimal viscosity solution of \eqref{eq:elliptobstacle}.
    \end{enumerate}
\end{lemma}

\begin{proof}
    Every $w\in \mathcal{W}$ satisfies $w\geq\varphi$, see \Cref{rem:above-obstacle-elliptic}. Thus,
    \[
        u(x)=\inf_{w\in \mathcal{W}} w(x) \geq \inf_{w\in \mathcal{W}} \varphi(x) = \varphi(x).
    \]
    This proves \eqref{lem:infsuper2-item1} and, in particular, that $u$ is well defined, since it is the infimum of a family of functions that is uniformly bounded from below. Property \eqref{lem:infsuper2-item2} is immediate from the definition of $u$.

    We next prove \eqref{lem:infsuper2-item3}, which in particular shows that $u$ is uniformly continuous. For $y\in\R^d$, let $c_y\coloneqq\|\varphi(\cdot+y)-\varphi\|_{L^\infty(\R^d)}$. If $w\in\mathcal{W}$, then $w(\cdot+y)+c_y\in\mathcal{W}$, since the equation is translation invariant and $w(x+y)+c_y\ge\varphi(x+y)+c_y\ge\varphi(x)$. Hence, by \eqref{lem:infsuper2-item2}, $u\le w(\cdot+y)+c_y$ for every $w\in\mathcal{W}$, that is, $u\le u(\cdot+y)+c_y$ in $\R^d$. Replacing $y$ by $-y$ and $x$ by $x+y$ we obtain $|u(\cdot+y)-u|\le c_y$, which is \eqref{lem:infsuper2-item3}; note that $c_y\to0$ as $y\to0$ because $\varphi$ is uniformly continuous. In particular, $u$ is continuous, which will be used below (and is needed for $u$ to be a viscosity solution).

    Let us now prove \eqref{lem:infsuper2-item4}. We first check that $u$ is a viscosity supersolution. By \eqref{lem:infsuper2-item1}, we already know that $u\geq\varphi$ in $\R^d$. Thus, we need to check the following: whenever $x_0\in \R^d$, $r>0$ and  $\phi\in C^2(\overline{B_r(x_0)}) \cap C_{\textup{b}}(\R^d\setminus \overline{B_r(x_0)})$ are such that
\begin{align}
\label{eq:visctouch2-superproof}
\phi(x)-u(x)&< \phi(x_0)-u(x_0) \quad \textup{for all }  x\in \overline{B_r(x_0)}\setminus\{x_0\},\\
\label{eq:visctouch-super-elliptic2-proof}
\phi(x)-u(x)&\leq \phi(x_0)-u(x_0)-\eta \quad \textup{for all } x\in \R^d\setminus\overline{B_r(x_0)},
\end{align}
for some $\eta>0$, then
\begin{align}\label{eq:viscsuperproof2}
   \fl \phi(x_0) &\geq 0.
\end{align}
By definition of $u$, there exists a sequence $\{w_\delta\}_{\delta>0}\subset \mathcal{W}$ such that
\[
    w_\delta(x_0) \searrow u(x_0) \quad \textup{as }\delta\to0^+.
\]
Moreover, $w_\delta \geq u$ for all $\delta>0$ by \eqref{lem:infsuper2-item2}. Thus, since $u$ is lower semicontinuous, we can find a sequence $\{x_\delta\}_{\delta>0}\subset \R^d$ with $x_\delta\to x_0$ as $\delta\to0^+$ such that $\phi$ is an admissible test function for $w_\delta$ at $x_\delta$. By definition of viscosity supersolution, $\fl \phi(x_\delta) \geq0$, and passing to the limit as $\delta\to0^+$ we obtain \eqref{eq:viscsuperproof2}.

We now check that $u$ is a viscosity subsolution. Let $x_0\in\R^d$. If $u(x_0)=\varphi(x_0)$, there is nothing to prove. If $u(x_0)>\varphi(x_0)$, we argue by contradiction and assume that there is a test function $\phi$ as in \Cref{def:visc-elliptic}, with $\phi(x_0)=u(x_0)$, touching $u$ from above at $x_0$, such that
\[
    \fl \phi(x_0)=\delta>0.
\]
By \Cref{rem:viscbetter} we may assume that $\phi-u\ge\eta_0>0$ in $\R^d\setminus\overline{B_r(x_0)}$. By continuity there is $\eta\in(0,r]$ such that
\[
    \fl \phi> \frac\delta2 \quad\textup{and}\quad \phi>\varphi \quad \textup{in }\overline{B_\eta(x_0)},
\]
and, since $\phi-u>0$ in $\overline{B_r(x_0)}\setminus\{x_0\}$, we also have $\phi-u\ge m_0>0$ in $\overline{B_r(x_0)}\setminus B_{\eta/2}(x_0)$. Choose $\veps>0$ smaller than $\eta_0$, $m_0$ and $\min_{\overline{B_{\eta/2}(x_0)}}(\phi-\varphi)$, and set $\phi_\veps\coloneqq\phi-\veps$. Then $\phi_\veps>u$ in $\R^d\setminus B_{\eta/2}(x_0)$, $\phi_\veps>\varphi$ in $B_{\eta/2}(x_0)$, $\phi_\veps(x_0)<u(x_0)$, and $\fl\phi_\veps=\fl\phi>\delta/2$ in $B_\eta(x_0)$. Hence
\[
    u_\veps\coloneqq\min\{\phi_\veps,u\}
\]
is bounded and continuous, $u_\veps\ge\varphi$, $u_\veps=u$ outside $B_{\eta/2}(x_0)$, and $u_\veps(x_0)<u(x_0)$. We claim that $u_\veps$ is a viscosity supersolution of~\eqref{eq:elliptobstacle}. Let $\psi$ be a test function touching $u_\veps$ from below at $x_1$. If $u_\veps(x_1)=u(x_1)$, then, since $u_\veps\le u$, $\psi$ also touches $u$ from below at $x_1$, and the required inequality follows from the fact that $u$ is a viscosity supersolution. Otherwise $u_\veps(x_1)=\phi_\veps(x_1)<u(x_1)$, so that $x_1\in B_{\eta/2}(x_0)$, and since $u_\veps\le\phi_\veps$ the function $\psi-\phi_\veps$ has a global maximum at $x_1$; therefore $\fl\psi(x_1)\ge\fl\phi_\veps(x_1)>0$, and also $u_\veps(x_1)-\varphi(x_1)=\phi_\veps(x_1)-\varphi(x_1)>0$. This proves the claim. Thus $u_\veps\in\mathcal W$ and, by \eqref{lem:infsuper2-item2}, $u\le u_\veps$, which contradicts $u_\veps(x_0)<u(x_0)$.

Finally, \eqref{lem:infsuper2-item5} follows from \eqref{lem:infsuper2-item4} and \eqref{lem:infsuper2-item2}.
\end{proof}

As already mentioned, the elliptic obstacle problem does not have uniqueness in general. However, we will prove that, if $d>2s$, there is a unique solution with a prescribed limit at infinity. The first step is to construct a minimal solution with this limit. This is achieved by considering the elliptic obstacle problem with obstacle $\max\{\varphi,c_*\}$, whose minimal solution turns out to have the desired limit.

In the sequel we will compare viscosity sub- and supersolutions of the time-fractional heat equation with its mild solutions, given through the fundamental solution $Z$, whose existence and estimates were obtained in~\cite{Kochubei-1990,Eidelman-Kochubei-2004} for $s=1$ and in~\cite{Kemppainen-Siljander-Zacher-2017} for $s\in(0,1)$. The following observation makes this possible.

\begin{lemma}\label{lem:mild-viscosity}
    Let $s\in(0,1]$, $\alpha\in(0,1)$ and $w_0\in C_\textup{b}(\R^d)$. Then the mild solution $w(\cdot,t)\coloneqq Z(\cdot,t)*w_0$ of the time-fractional heat equation $\partial_t^\alpha w+\fl w=0$ with initial datum $w_0$ is a viscosity solution of this equation in $\R^d\times(0,\infty)$, in the sense of \Cref{def:viscsol-p-Laplace} without obstacle. Consequently, if $u\in C_\textup{b}(\R^d\times[0,\infty))$ (resp.\ $v$) is a viscosity supersolution (resp.\ subsolution) of the time-fractional heat equation in $\R^d\times(0,\infty)$ with $u(\cdot,0)\ge w_0$ (resp.\ $v(\cdot,0)\le w_0$), then $u\ge w$ (resp.\ $v\le w$) in $\R^d\times[0,\infty)$.
\end{lemma}

\begin{proof}
If $w_0\in C^\infty_\textup{c}(\R^d)$, then $w$ is a classical solution of the time-fractional heat equation, see \cite[Theorem~2.12]{Kemppainen-Siljander-Zacher-2017}, which covers the whole range $s\in(0,1]$ under the assumptions $w_0,\widehat{w_0}\in L^1(\R^d)$; see also~\cite{Eidelman-Kochubei-2004} for the case $s=1$, where bounded and locally Hölder continuous data suffice. In fact, in this case much more is true: by the representation $\widehat{w}(\xi,t)=E_\alpha(-|\xi|^{2s}t^\alpha)\widehat{w_0}(\xi)$, with $E_\alpha$ the Mittag-Leffler function, and the fast decay of $\widehat{w_0}$, the function $w$ belongs to $C^\infty(\R^d\times(0,\infty))$, with $x$-derivatives of all orders bounded in $\R^d\times[t_0,T]$ for every $0<t_0<T$. Hence, for $t>0$, the fractional Laplacian of $w(\cdot,t)$ is given pointwise by the singular integral~\eqref{eq:definition.fractional.laplacian}, and the Caputo derivative of $w(x,\cdot)$ may be computed by any of the equivalent formulas~\eqref{eq:capdef} and~\eqref{eq:capdef2}, both integrals being absolutely convergent. This is all that is needed to check that $w$ is a viscosity solution: if $\phi$ is a test function touching $w$ from above at $(x_0,t_0)$, as in \Cref{def:viscsol-p-Laplace}, then $\phi(x_0,t_0)-\phi(x_0,r)\le w(x_0,t_0)-w(x_0,r)$ for $r\in[0,t_0]$ and $\phi(x_0,t_0)-\phi(y,t_0)\le w(x_0,t_0)-w(y,t_0)$ for $y\in\R^d$, so that~\eqref{eq:capdef} and the integral formula for $\fl$ give $\partial_t^\alpha\phi(x_0,t_0)+\fl\phi(x_0,t_0)\le\partial_t^\alpha w(x_0,t_0)+\fl w(x_0,t_0)=0$; the supersolution inequality is obtained in the same way.

For general $w_0\in C_\textup{b}(\R^d)$, let $w_{0,\delta}\in C^\infty_\textup{c}(\R^d)$ with $\|w_{0,\delta}\|_{L^\infty(\R^d)}\le\|w_0\|_{L^\infty(\R^d)}$ and $w_{0,\delta}\to w_0$ locally uniformly as $\delta\to0^+$. Since $Z(\cdot,t)$ is a probability density, the corresponding mild solutions $w_\delta=Z(\cdot,t)*w_{0,\delta}$ are uniformly bounded and converge to $w$ locally uniformly in $\R^d\times[0,\infty)$ (split the convolution into $|y|\le R$ and $|y|>R$, and use that the tails of $Z(\cdot,t)$ are uniformly small for $t\in[0,T]$, by its self-similar form).

It remains to show that $w$ is a viscosity solution, which follows from the classical stability argument (cf.~\cite[Lemma~2.1]{Allen-2018}); we sketch it for supersolutions, the other case being analogous. Let $\phi$ be a test function for $w$ at $(x_0,t_0)$, which we may take smooth by \Cref{rem:viscbetter}\eqref{item:moregeneraltestc}, and let $m\coloneqq(\phi-w)(x_0,t_0)$. By the locally uniform convergence, $\phi-w_{\delta}$ attains its maximum over a small cylinder at points $(x_\delta,t_\delta)\to(x_0,t_0)$, with maximum values $m_\delta\to m$, and $\psi_\delta\coloneqq\min\{\phi,w_\delta+m_\delta\}$ is an admissible test function for $w_\delta$ at $(x_\delta,t_\delta)$, see \Cref{rem:viscbetter}\eqref{item:moregeneraltesta}. Hence $\partial_t^\alpha\psi_\delta (x_\delta,t_\delta)+\fl\psi_\delta(x_\delta,t_\delta)\ge0$, and since $\psi_\delta\to\phi$ locally uniformly, we may pass to the limit by dominated convergence (the kernels are bounded away from the singularity), obtaining $\partial_t^\alpha\phi(x_0,t_0)+\fl\phi(x_0,t_0)\ge0$. Hence $w$ is a viscosity solution. The comparison statement then follows from \Cref{lem:compErwin} with $\Omega=\R^d\times(0,\infty)$, in which case the hypothesis on the complement of $\Omega$ reduces to the ordering of the initial data.
\end{proof}

We now turn to the elliptic obstacle problem~\eqref{eq:elliptobstacle} when $d>2s$, in which case it has, in general, many solutions. For each prescribed value $c_*\ge\limsup_{|x|\to\infty}\varphi(x)$ at infinity, the following lemma singles out one of them, which is the large-time limit of the evolution problem when the initial datum has the limit $c_*$ at infinity, see \Cref{thm:limit-at-infinity}. It is obtained as the minimal solution with the modified obstacle $\max\{\varphi,c_*\}$.

\begin{lemma}\label{lem:minimal.with.limit}
    Let $d>2s$, $\varphi\in C^0_\textup{b}(\R^d)$ and $c_*\ge\limsup_{|x|\to\infty}\varphi(x)$. Let $u_{\varphi,c_*}$ be
    the minimal viscosity solution of the elliptic obstacle problem \eqref{eq:elliptobstacle} with obstacle $\max\{\varphi,c_*\}$. Then:
    \begin{enumerate}[\rm (a)]
        \item\label{lem:minimal.with.limit-item1} $u_{\varphi,c_*}$ is a viscosity solution of the elliptic obstacle problem~\eqref{eq:elliptobstacle} with obstacle $\varphi$;
        \item\label{lem:minimal.with.limit-item2} $\lim_{|x|\to\infty}u_{\varphi,c_*}(x)=c_*$;
        \item\label{lem:minimal.with.limit-item3} if $u$ is a viscosity solution of the elliptic obstacle problem \eqref{eq:elliptobstacle} with obstacle $\varphi$ such that ~$\lim_{|x|\to\infty}u(x)=c_*$, then $u_{\varphi,c_*}\le u$ in~$\mathbb{R}^d$.
    \end{enumerate}
\end{lemma}

\begin{proof}
    \eqref{lem:minimal.with.limit-item1} We know that $\fl u_{\varphi,c_*}\ge0$ in the viscosity sense and that $u_{\varphi,c_*}\ge \max\{\varphi,c_*\}\ge c_*$ in $\R^d$. Hence, by the strong minimum principle for $s$-superharmonic functions, $u_{\varphi,c_*}$ cannot attain the value $c_*$ at any point unless it is identically equal to $c_*$. (For viscosity supersolutions this follows directly from the definition: if $u_{\varphi,c_*}(x_0)=c_*$, testing with $\phi\equiv c_*$ in $\overline{B_r(x_0)}$ and $\phi=u_{\varphi,c_*}$ outside, which is admissible by \Cref{rem:viscbetter}, gives $\int_{|y|\ge r}(c_*-u_{\varphi,c_*}(x_0+y))|y|^{-d-2s}\dd y\ge0$ when $s\in(0,1)$, whence $u_{\varphi,c_*}=c_*$ outside $B_r(x_0)$ for every $r>0$.  For $s=1$ it is the classical strong minimum principle, applicable by \Cref{prop:equivalence.superharmonic} and \cite[Chapter~II, Section~4, Propositions~20 and~21]{Dautray-Lions-1990}.

    We now distinguish two cases. If $\varphi\le c_*$ in $\R^d$, then $\max\{\varphi,c_*\}=c_*$, so that $u_{\varphi,c_*}\equiv c_*$, which is a viscosity solution of the elliptic obstacle problem with obstacle $\varphi$. If, on the contrary, $c_*<\sup\varphi$, then $u_{\varphi,c_*}$ is not identically $c_*$, so that it does not attain the value $c_*$, and hence $u_{\varphi,c_*}>c_*$ everywhere in $\R^d$. In particular, since $u_{\varphi,c_*}\ge\max\{\varphi,c_*\}\ge\varphi$ in $\R^d$, the function $u_{\varphi,c_*}$ is a viscosity supersolution of the obstacle problem with obstacle~$\varphi$. Moreover, the subsolution condition at points of the non-contact set $\{u_{\varphi,c_*}>\varphi\}$ is inherited from the fact that $u_{\varphi,c_*}$ solves the obstacle problem with obstacle $\max\{\varphi,c_*\}$: indeed, since $u_{\varphi,c_*}>c_*$, we have $\{u_{\varphi,c_*}>\varphi\}=\{u_{\varphi,c_*}>\max\{\varphi,c_*\}\}$, and on this set the subsolution condition for both problems reads $\fl\phi(x_0)\le0$ for every admissible test function $\phi$. Thus $u_{\varphi,c_*}$ is a viscosity solution of the elliptic obstacle problem with obstacle~$\varphi$.

    \noindent \eqref{lem:minimal.with.limit-item2} Fix $\varepsilon>0$ and consider the function
    \[
    U_\varepsilon(x)\coloneqq A_\varepsilon\min\{|x|^{2s-d},1\}+c_*+\varepsilon,
    \]
    where $A_\varepsilon>0$ is chosen below. Since $|x|^{2s-d}$ is, up to a positive constant, the fundamental solution of $\fl$ in $\R^d$ when $d>2s$, both $|x|^{2s-d}$ and the constant $1$ are $s$-superharmonic, and so is their minimum. Hence $U_\varepsilon$ is a bounded viscosity supersolution of $\fl U_\varepsilon=0$ in $\R^d$. Moreover, since $c_*\ge\limsup_{|x|\to\infty}\varphi(x)$, there is $R_\varepsilon>1$ such that $\varphi\le c_*+\varepsilon$ in $\R^d\setminus B_{R_\varepsilon}$, and choosing $A_\varepsilon\ge R_\varepsilon^{d-2s}\,(\sup\varphi-c_*)^+$ we also get $U_\varepsilon\ge\varphi$ in $B_{R_\varepsilon}$. Therefore $U_\varepsilon\ge\max\{\varphi,c_*\}$ in $\R^d$, and $U_\varepsilon$ is a viscosity supersolution of the obstacle problem with obstacle $\max\{\varphi,c_*\}$. By \Cref{lem:infsuper2}\eqref{lem:infsuper2-item2}, $c_*\le u_{\varphi,c_*}\le U_\varepsilon$ in $\R^d$. Letting first $|x|\to\infty$ and then $\varepsilon\to0^+$ we get the result.

    \noindent \eqref{lem:minimal.with.limit-item3} Regarded as a function of $(x,t)$, the solution $u$ of the elliptic problem is a (stationary) viscosity solution of the evolution obstacle problem~\eqref{eq:main}--\eqref{eq:initial} with initial datum $u$, hence a viscosity supersolution of the time-fractional heat equation with the same initial datum. By \Cref{lem:mild-viscosity}, $u$ lies above the mild solution $w(\cdot,t)=Z(\cdot,t)*u$ of the time-fractional heat equation with initial datum $u$. Since $Z(\cdot,t)$ is a radially nonincreasing probability density with $Z(y,t)=t^{-\frac{\alpha d}{2s}}F(yt^{-\frac{\alpha}{2s}})$, we have $\sup_{z\in\R^d}\int_{B_R(z)}Z(y,t)\dd y=\int_{B_{Rt^{-\alpha/(2s)}}}F(\xi)\dd\xi\to0$ as $t\to+\infty$ for every $R>0$; together with $\lim_{|x|\to\infty}u(x)=c_*$, this gives $w(\cdot,t)\to c_*$ uniformly in $\R^d$ as $t\to+\infty$ (split the convolution into $|y|\le R$ and $|y|>R$ and choose $R$ large). Thus $u\ge c_*$, and since also $u\ge\varphi$, we conclude that $u\ge\max\{\varphi,c_*\}$; that is, $u$ is a viscosity supersolution of the obstacle problem with obstacle $\max\{\varphi,c_*\}$. The minimality of~$u_{\varphi,c_*}$ gives the result.
\end{proof}

We have thus shown that, when $d>2s$, the minimal viscosity solution of~\eqref{eq:elliptobstacle} with obstacle $\max\{\varphi,c_*\}$ is also a viscosity solution of~\eqref{eq:elliptobstacle} with obstacle $\varphi$ attaining the value $c_*$ at infinity, and that it is the smallest one with this property. Our next goal, \Cref{thm:uniqueness.elliptic} below, is to show that there is no other one. The argument follows a classical pattern: if $v$ were another solution with limit $c_*$ at infinity, then $v-u_{\varphi,c_*}$ would attain a positive maximum at some point of the non-coincidence set of $v$, where $v$ is $s$-harmonic and $u_{\varphi,c_*}$ is $s$-superharmonic, and the strong maximum principle should force $v-u_{\varphi,c_*}$ to be constant, which is incompatible with its decay at infinity. For viscosity solutions, however, the difference of a subsolution and a supersolution is not a subsolution in any obvious sense, since the notion is defined through test functions touching each function separately. This is why the comparison result needed for the proof of uniqueness, \Cref{lem:strong-max-difference} in the Appendix, is proved by passing to the distributional formulation, available for continuous bounded viscosity solutions by \cite{Korvenpaa-Kuusi-Lindgren-2019} (for $s\in(0,1)$) and \cite{Ishii-1995} (for $s=1$). In that formulation the difference of a subsolution and a supersolution is a subsolution, by linearity, and the strong maximum principle follows from an elementary regularization argument.

\begin{theorem}\label{thm:uniqueness.elliptic}
    Let $d>2s$, $\varphi\in C^0_\textup{b}(\R^d)$ and $c_*\in\R$. There is at most one viscosity solution $u$ of problem~\eqref{eq:elliptobstacle} such that
    $$
        \lim\limits_{|x|\to\infty}u(x)=c_*.
    $$
    Moreover, if $$c_*\ge\limsup_{|x|\to\infty}\varphi(x),$$ such a solution exists, and it is given by $u_{\varphi,c_*}$.
\end{theorem}

\begin{proof}
    Existence when $c_*\ge\limsup_{|x|\to\infty}\varphi(x)$ is the content of \Cref{lem:minimal.with.limit}. Note also that if a solution $u$ with $\lim_{|x|\to\infty}u(x)=c_*$ exists, then necessarily $c_*\ge\limsup_{|x|\to\infty}\varphi(x)$, since $u\ge\varphi$. Let then $u=u_{\varphi,c_*}$ be the minimal solution with this behaviour at infinity, and let $v$ be any solution with the same behaviour. We have $v\ge u$ by \Cref{lem:minimal.with.limit}\eqref{lem:minimal.with.limit-item3}. For the reverse inequality we apply \Cref{lem:strong-max-difference} with $\Omega\coloneqq\{v>\varphi\}$: $v$ is a viscosity subsolution of $\fl v\le0$ in $\Omega$, $u$ is a viscosity supersolution of $\fl u\ge0$ in $\R^d$, $v=\varphi\le u$ in $\R^d\setminus\Omega$, and $v-u\to c_*-c_*=0$ at infinity. Hence $v\le u$, and thus $v=u$.
\end{proof}

If $d\le2s$, the minimal solution $u_{\varphi,c_*}$ is a constant, which coincides with $c_*$ only if $c_*\ge\sup\varphi$.
\begin{lemma}\label{lem:small-dimension} 
    Let $d\le 2s$, $\varphi\in C^0_\textup{b}(\R^d)$ and $$c_*\ge\limsup_{|x|\to\infty}\varphi(x).$$ Let $u_{\varphi,c_*}$ be the minimal viscosity solution of the elliptic obstacle problem \eqref{eq:elliptobstacle} with obstacle $\max\{\varphi,c_*\}$. Then $u_{\varphi,c_*}\equiv\max\{\sup\varphi,c_*\}$.
\end{lemma}

\begin{proof}
    The function $u_{\varphi,c_*}$ is a bounded continuous viscosity supersolution of $\fl u=0$ in $\R^d$. Hence, it is a bounded distributional supersolution of that equation (see~\Cref{prop:equivalence.superharmonic}), and therefore a constant (see~\Cref{prop:Liouville.superharmonic}). Since it lies above the obstacle $\max\{\varphi,c_*\}$, this constant is at least $\max\{\sup\varphi,c_*\}$, and the constant $\max\{\sup\varphi,c_*\}$ is itself a solution. The minimality of $u_{\varphi,c_*}$ yields the result.
\end{proof}

\begin{remark}\label{rem:discrete-minimal-with-limit}
{\rm Analogous results hold, with analogous proofs, for the minimal solution of the discrete elliptic obstacle problem~\eqref{eq:ellipticObs} with obstacle $\max\{\varphi,c_*\}$, for every fixed $h>0$. When $d>2s$, the role of the barrier $\min\{|x|^{2s-d},1\}$ is played by the fundamental solution $\Phi_h$ of the discrete operator $\fld$ in $h\Z^d$, which has the same behaviour as that of $\fl$: it is positive, discrete superharmonic, and
satisfies
\[
    \frac{c_{d,s}}{|x_\alpha|^{d-2s}}\le\Phi_h(x_\alpha)\le\frac{C_{d,s}}{|x_\alpha|^{d-2s}},\qquad x_\alpha\in h\Z^d\setminus\{0\},
\]
with constants independent of $h$. This is easily obtained by writing $\Phi_h$ as the subordination of the discrete heat kernel, $\Phi_h=\frac{h^{-d}}{\Gamma(s)}\int_0^\infty e^{t\Delta_h}\mathbf 1_{\{0\}}\,t^{s-1}\dd t$, which converges precisely when $d>2s$, and using the two-sided Gaussian bounds for the latter.}
\end{remark}

\subsection{Large-time behaviour for the time-fractional obstacle problem}

We now tackle the problem of characterizing the large-time behaviour of viscosity solutions to~\eqref{eq:main}--\eqref{eq:initial}. Mimicking what we did in the discrete case, the first step is to prove that the Caputo derivative of the solution goes to $0$ in the large-time limit. Some care is needed here: viscosity solutions need not have a Caputo derivative in the classical sense. However, this is guaranteed if the initial datum and the obstacle are smooth, as we see next.

\begin{proposition}\label{prop:CaputoDefined}
    Let $s\in(0,1]$ and $\alpha\in(0,1)$. Let $u_0,\varphi\in {\X^{2}_{\textup{b}}}(\R^d)$ be such that $u_0\ge\varphi$ in $\R^d$, and let $u$ be the viscosity solution of~\eqref{eq:main}--\eqref{eq:initial}. Then:
    \begin{enumerate}[\rm (a)]
        \item\label{prop:CaputoDefined-itema} $\partial_t^\alpha u(x,t)$ is well defined for all $(x,t)\in \R^d\times (0,\infty)$.
        \item\label{prop:CaputoDefined-itemb} For all $t>0$, we have that
        \begin{align*}
            \|\partial_t^\alpha u(\cdot,t)\|_{L^\infty(\R^d)}\leq K_{d,s,\alpha} \max\{\|u_0\|_{{\X^{2}}(\R^d)}\coma \|\varphi\|_{{\X^{2}}(\R^d)}\} \min\left\{1\coma \frac{1}{t^{\alpha(1-\alpha)}}\right\}.
        \end{align*}
        In particular, $\|\partial_t^\alpha u(\cdot,t)\|_{L^\infty(\R^d)}\to0$ as $t\to+\infty$.
    \end{enumerate}
\end{proposition}

\begin{proof}
By \Cref{main:ConvAndPropViscositySolutionParabolic}\eqref{main:ConvAndPropViscositySolutionParabolic-item4} with $a=b=2$ (and $\veps=2-2s$ if $s\in(0,1)$), we have
\begin{align}\label{eq:Lip-in-talpha}
    \|u(\cdot,t)-u(\cdot,\tilde{t})\|_{L^\infty(\R^d)} \leq \tilde{K}  |t^\alpha-\tilde{t}^\alpha|\quad \text{for all }t,\tilde{t}\in[0,\infty),
\end{align}
with $\tilde{K}=\widetilde{C}_{d,s,\alpha} \max\{\|u_0\|_{\X^{2}(\R^d)}\coma \|\varphi\|_{\X^{2}(\R^d)}\}$. Since $t^\alpha-s^\alpha\le\alpha s^{\alpha-1}(t-s)$ for $0<s<t$, the integral in~\eqref{eq:capdef} is absolutely convergent, which proves \eqref{prop:CaputoDefined-itema}. Using \eqref{eq:Lip-in-talpha} we get, on the one hand,
\begin{align*}
    |\partial_t^\alpha u(x,t)|&=\left|\frac{1}{\Gamma(1-\alpha)}\frac{u(x,t)-u(x,0)}{t^{\alpha}}+\frac{\alpha}{\Gamma(1-\alpha)}\int_0^t \frac{u(x,t)-u(x,s)}{(t-s)^{1+\alpha}} \dd s\right|\\
    &\leq \frac{1}{\Gamma(1-\alpha)}  \left(\frac{\|u(\cdot,t)-u(\cdot,0)\|_{L^\infty(\R^d)}}{t^{\alpha}}+\alpha\int_0^t \frac{\|u(\cdot,t)-u(\cdot,s)\|_{L^\infty(\R^d)}}{(t-s)^{1+\alpha}} \dd s\right)\\
    &\leq \frac{\tilde{K}}{\Gamma(1-\alpha)} \left(1+\alpha\int_0^t \frac{t^{\alpha}-s^\alpha}{(t-s)^{1+\alpha}} \dd s\right)\\
    &=  \frac{\tilde{K}}{\Gamma(1-\alpha)} \left(1+\alpha^2\int_0^t \frac{\int_s^t \eta^{\alpha-1}\dd \eta}{(t-s)^{1+\alpha}} \dd s\right)\\
    &= \frac{\tilde{K}}{\Gamma(1-\alpha)} \left(1+\alpha^2\int_0^t \eta^{\alpha-1} \int_0^\eta\frac{\dd s}{(t-s)^{1+\alpha}} \dd \eta\right)\\
    &=\frac{\tilde{K}}{\Gamma(1-\alpha)} \left(1+\alpha\int_0^t \eta^{\alpha-1} \left((t-\eta)^{-\alpha}-t^{-\alpha}\right) \dd \eta\right)\\
    &\leq \frac{\tilde{K}}{\Gamma(1-\alpha)} \left(1+\alpha\int_0^1 \sigma^{\alpha-1} (1-\sigma)^{-\alpha} \dd \sigma\right)\\
    &= \frac{\tilde{K}}{\Gamma(1-\alpha)}
    \left(1+\alpha\Gamma(\alpha)\Gamma(1-\alpha)\right),
\end{align*}
which proves that $\|\partial_t^\alpha u(\cdot,t)\|_{L^\infty(\R^d)}\leq K_{d,s,\alpha} \max\{\|u_0\|_{\X^{2}(\R^d)}\coma \|\varphi\|_{\X^{2}(\R^d)}\}$. On the other hand, for any $R\in (0,t)$, we also have
\begin{align*}
    |\partial_t^\alpha u(x,t)|&\leq \frac{1}{\Gamma(1-\alpha)}  \left(\frac{\|u(\cdot,t)-u(\cdot,0)\|_{L^\infty(\R^d)}}{t^{\alpha}}+\alpha\int_0^t \frac{\|u(\cdot,t)-u(\cdot,s)\|_{L^\infty(\R^d)}}{(t-s)^{1+\alpha}} \dd s\right)\\
    &\leq \frac{\tilde{K}}{\Gamma(1-\alpha)} \left(\frac{2}{t^{\alpha}}+2\alpha\int_0^{t-R}\frac{\dd s}{(t-s)^{1+\alpha}} +\alpha\int_{t-R}^t \frac{t^{\alpha}-s^\alpha}{(t-s)^{1+\alpha}} \dd s\right)\\
    &\leq \frac{2\tilde{K}}{\Gamma(1-\alpha)} \left(\frac{1}{t^{\alpha}}+\frac{1}{R^{\alpha}} +\alpha (t-R)^{\alpha-1}\int_{t-R}^t \frac{\dd s}{(t-s)^{\alpha}}\right)\\
    &\leq \frac{2\tilde{K}}{\Gamma(1-\alpha)} \left(\frac{1}{t^{\alpha}}+\frac{1}{R^{\alpha}} +\frac{\alpha}{1-\alpha} (t-R)^{\alpha-1}R^{1-\alpha}\right),
\end{align*}
where we have used that $\|u(\cdot,t)-u(\cdot,s)\|_{L^\infty(\R^d)}\le 2\max\{\|u_0\|_{L^\infty(\R^d)},\|\varphi\|_{L^\infty(\R^d)}\}\le 2\tilde K$ (enlarging $\widetilde C_{d,s,\alpha}$ if necessary) and $t^\alpha-s^\alpha\le\alpha(t-R)^{\alpha-1}(t-s)$ for $s\in(t-R,t)$. For $t\ge2$, the choice $R=t^{1-\alpha}=t^{-\alpha}t\le 2^{-\alpha}t$ gives $(t-R)^{\alpha-1}\le (1-2^{-\alpha})^{\alpha-1}t^{\alpha-1}$, and hence
\begin{align*}
    |\partial_t^\alpha u(x,t)|&\leq \frac{2\tilde{K}}{\Gamma(1-\alpha)} \left(2 +\frac{(1-2^{-\alpha})^{\alpha-1}\alpha}{1-\alpha} \right)\frac{1}{t^{\alpha(1-\alpha)}}.
\end{align*}
Together with the bound for $t\le2$ obtained above, this proves \eqref{prop:CaputoDefined-itemb}.
\end{proof}

When trying to adapt the proof of the large-time behaviour in the discrete case to the continuous one, we find a difficulty (cf.\ \cite[Proposition~2.4]{Erwin}, where a related property is established for sub- and supersolutions that are $\alpha'$-Hölder continuous in time for some $\alpha'>\alpha$, a slightly stronger requirement than the one available here): solutions are meant in the viscosity sense, and the decay to $0$ of the Caputo derivative does not by itself guarantee that the limit of $u(\cdot,t_k)$ along a subsequence is a solution of the elliptic obstacle problem, since the test functions in the definition of viscosity solution depend on both $x$ and $t$. Hence, we first need to show that viscosity solutions having a Caputo derivative are \lq\lq classical in time'': the Caputo derivative $\partial_t^\alpha u$ itself can be used in place of the time derivative of the test function, so that only test functions in space are needed. This property is inherited from the discrete approximation. The first step is to show that the discrete Caputo derivative of the numerical solution converges to the Caputo derivative of $u$.

\begin{lemma}\label{lem:convergence-discrete-Caputo}
    Let the assumptions of \Cref{prop:CaputoDefined} hold, and assume \eqref{as:CFL}. Let $U_{\tau,h}$ be the solution of the numerical scheme~\eqref{eq:numsch} with initial datum $u_0$ and obstacle $\varphi$, and write $D_\tau^\alpha U_{\tau,h}(x,t)\coloneqq D_\tau^\alpha (U_{\tau,h})^n(x)$ for $t=\tau n$, $n\in\N$. Then, for every $(x_0,t_0)\in\R^d\times(0,\infty)$ and every family $(x_{\tau,h},t_{\tau,h})\in\R^d\times\tau\N$ such that $(x_{\tau,h},t_{\tau,h})\to(x_0,t_0)$ as $\tau,h\to0^+$, we have
    \begin{align*}
        D_\tau^\alpha U_{\tau,h}(x_{\tau,h},t_{\tau,h})\to \partial_t^\alpha u(x_0,t_0) \quad \textup{as }  \tau,h \to 0^+.
    \end{align*}
    In particular, $\partial_t^\alpha u$ is continuous in $\R^d\times(0,\infty)$.
\end{lemma}

\begin{proof}
We write $(x_*,t_*)=(x_{\tau,h},t_{\tau,h})$ for short and use the decomposition
\begin{align*}
     |D_\tau^\alpha U_{\tau,h}(x_*,t_*)-\partial_t^\alpha u(x_0,t_0)| \leq &|D_\tau^\alpha U_{\tau,h}(x_*,t_*)-D_\tau^\alpha u(x_*,t_*)|+|D_\tau^\alpha u(x_*,t_*)-\partial_t^\alpha u(x_*,t_*)|\\
     &+|\partial_t^\alpha u(x_*,t_*)-\partial_t^\alpha u(x_0,t_0)|\eqqcolon I_1+I_2+I_3,
\end{align*}
where $D^\alpha_\tau u(x,t)$ denotes the L1 discretization~\eqref{eq:def.discretization.Caputo} applied to the sequence $\{u(x,\tau n)\}_{n\in\N}$. We check that each of the three terms tends to $0$ as $\tau,h\to0^+$.

Let us estimate $I_1$ first. Define $W_{\tau,h}(x,t)\coloneqq U_{\tau,h}(x,t)-u(x,t)$ for $t\in\tau\N$. By \Cref{lem:est-schm}, \Cref{lem:timeequi-nosmooth} (with $a=2$, so that $f_{2,s}(h)=1$) and \Cref{main:ConvAndPropViscositySolutionParabolic}\eqref{main:ConvAndPropViscositySolutionParabolic-itema}, \eqref{main:ConvAndPropViscositySolutionParabolic-itemb} and~\eqref{main:ConvAndPropViscositySolutionParabolic-item4}, we have, for all $t,\tilde{t}\in \tau \N$ with $t> \tilde{t}$,
\begin{align}\label{eq:threekeyestimates}
\begin{aligned}
    &\|W_{\tau,h}(\cdot,t)\|_{L^\infty(\R^d)}\leq 2\max\{\|u_0\|_{L^\infty(\R^d)},\|\varphi\|_{L^\infty(\R^d)}\},\\
    &\|W_{\tau,h}(\cdot,t)-W_{\tau,h}(\cdot,\tilde{t})\|_{L^\infty(\R^d)}\leq \tilde K (t^\alpha-\tilde t^\alpha)\le \alpha\tilde K\, \tilde t^{\,\alpha-1}(t-\tilde{t}),\\
    &W_{\tau,h} \to 0 \textup{ as } \tau,h\to0^+\textup{ locally uniformly in }  \R^d\times[0,\infty),
    \end{aligned}
\end{align}
with $\tilde K$ depending only on $d,s,\alpha$ and $\max\{\|u_0\|_{\X^{2}(\R^d)}\coma \|\varphi\|_{\X^{2}(\R^d)}\}$. Let us also recall the defining property~\eqref{eq:discretization.Caputo.2} of the L1 discretization,
\begin{align}\label{eq:propL1discretization}
D_\tau^\alpha W_{\tau,h}(x,t)=\partial_t^\alpha[\mathcal{I}_{\tau}[W_{\tau,h}]](x,t),\qquad t\in\tau\N,
\end{align}
where $\mathcal{I}_{\tau}[W_{\tau,h}]$ stands for the piecewise linear interpolation in time of $W_{\tau,h}$. Clearly, $\mathcal{I}_\tau[W_{\tau,h}]$ satisfies the estimates in \eqref{eq:threekeyestimates} for all $t>\tilde t\ge0$ (with $\tilde t^{\,\alpha-1}$ replaced by $\max\{\tilde t,\tau\}^{\alpha-1}$ in the second one, which is still bounded by $C\tilde t^{\,\alpha-1}$ for the values of $\tilde t$ that matter below) and converges to $0$ locally uniformly. Thus, by~\eqref{eq:capdef},
\begin{align*}
I_1&=|\partial_t^\alpha[\mathcal{I}_{\tau}[W_{\tau,h}]](x_*,t_*)|\\
&\leq \frac{1}{\Gamma(1-\alpha)}\frac{|\mathcal{I}_{\tau}[W_{\tau,h}](x_*,t_*)-\mathcal{I}_{\tau}[W_{\tau,h}](x_*,0)|}{t_*^{\alpha}} \\ & \quad +\frac{\alpha}{\Gamma(1-\alpha)}\int_0^{t_*} \frac{|\mathcal{I}_{\tau}[W_{\tau,h}](x_*,t_*)-\mathcal{I}_{\tau}[W_{\tau,h}](x_*,s)|}{(t_*-s)^{1+\alpha}} \dd s.
\end{align*}
The first term converges to $0$, since $(x_*,t_*)$ stays in a compact set and $\mathcal{I}_{\tau}[W_{\tau,h}]\to0$ locally uniformly. The second term goes to $0$ by the dominated convergence theorem (after the change of variables $s=t_*\sigma$, which makes the domain of integration independent of $\tau,h$): the integrand converges to $0$ for every $\sigma\in(0,1)$ by the locally uniform convergence of $\mathcal{I}_{\tau}[W_{\tau,h}]$ to $0$, while, by the second estimate in \eqref{eq:threekeyestimates}, it is dominated by
\[
C\,\frac{(t_*\sigma)^{\alpha-1}\,t_*(1-\sigma)}{(t_*(1-\sigma))^{1+\alpha}}\,t_*= C\,t_*^{\alpha-1}\sigma^{\alpha-1}(1-\sigma)^{-\alpha},
\]
which is integrable in $(0,1)$ and uniformly bounded in $\tau,h$, since $t_*\to t_0>0$.

For $I_2$ we note, using~\eqref{eq:propL1discretization} for the sequence $\{u(x,\tau n)\}_{n\in\N}$, that $I_2=|\partial_t^\alpha[\mathcal{I}_{\tau}[u]-u](x_*,t_*)|$. Since $u$ is continuous and satisfies~\eqref{eq:Lip-in-talpha}, the function $\mathcal{I}_{\tau}[u]-u$ converges to $0$ locally uniformly and satisfies the same kind of bounds as $W_{\tau,h}$ in \eqref{eq:threekeyestimates}. Hence $I_2\to0$ exactly as $I_1$.

Finally, $I_3\to0$ by the same dominated convergence argument applied to the continuous function $u$, which also shows that $\partial_t^\alpha u$ is continuous: the integrand in~\eqref{eq:capdef} at $(x_*,t_*)$, written in the variable $\sigma=s/t_*$, converges pointwise to the integrand at $(x_0,t_0)$ by continuity of $u$, and is dominated as above thanks to~\eqref{eq:Lip-in-talpha}.
\end{proof}

For regular enough data, the Caputo derivative $\partial_t^\alpha u$ of the viscosity solution is well defined by \Cref{prop:CaputoDefined}. We now show that the equation then holds classically in time and in the viscosity sense in space, which allows us to test only in space.

\begin{lemma}\label{lem:clasical-in-time-viscosity}
    Let the assumptions of \Cref{prop:CaputoDefined} hold. Let $(x_0,t_0)\in\R^d\times(0,\infty)$, $r>0$ and $\phi\in C^2(\overline{B_r(x_0)})\cap C_{\textup{b}}(\R^d\setminus \overline{B_r(x_0)})$ be such that
\begin{equation}\label{eq:touchabove}
            \phi(x)-u(x,t_0)\geq \phi(x_0)-u(x_0,t_0) \quad \textup{for all }x\in \R^d.
\end{equation}
    Then
    \begin{align}\label{eq:newviscsub}
        \min\{\partial_t^\alpha u(x_0,t_0)+\fl \phi(x_0) \coma u(x_0,t_0) - \varphi(x_0)\} &\leq 0.
    \end{align}
    Analogously, if \eqref{eq:touchabove} holds with the reverse sign, then
    \begin{align}\label{eq:newviscsuper}
        \min\{\partial_t^\alpha u(x_0,t_0)+\fl \phi(x_0) \coma u(x_0,t_0) - \varphi(x_0)\} &\geq 0.
    \end{align}
\end{lemma}

\begin{proof}
We only prove \eqref{eq:newviscsub}, since \eqref{eq:newviscsuper} follows in the same way. Let $\psi(x)\coloneqq|x-x_0|^4/(1+|x-x_0|^4)$ and $\phi_\varepsilon\coloneqq\phi+\varepsilon\psi$ for $\varepsilon>0$. Then $\phi_\varepsilon-u(\cdot,t_0)$ has a strict minimum at $x_0$ in $\R^d$, with a positive gap outside $\overline{B_r(x_0)}$. Since $\fl\phi_\varepsilon(x_0)=\fl\phi(x_0)+\varepsilon\fl\psi(x_0)\to\fl\phi(x_0)$ as $\varepsilon\to0^+$, it is enough to prove \eqref{eq:newviscsub} for $\phi_\varepsilon$, that is, to assume that the contact is strict, with a gap $\eta>0$ outside $\overline{B_r(x_0)}$. Since $u(\cdot,t_0)$ is continuous, by \Cref{rem:stability-elliptic} we may also assume that $\phi$ is a $C^2$ function in $\R^d$ with bounded derivatives up to order two. Let $U_{\tau,h}$ be the solution of the numerical scheme with initial datum $u_0$ and obstacle $\varphi$. By the locally uniform convergence $U_{\tau,h}\to u$, see \Cref{main:ConvAndPropViscositySolutionParabolic}\eqref{main:ConvAndPropViscositySolutionParabolic-itema}, we can find points $(x_{\tau,h},t_{\tau,h})\in\R^d\times\tau\N$ with $(x_{\tau,h},t_{\tau,h})\to(x_0,t_0)$ as $\tau,h\to0^+$ such that
\begin{align*}
\phi(x)-U_{\tau,h}(x,t_{\tau,h})&\geq \phi(x_{\tau,h})-U_{\tau,h}(x_{\tau,h},t_{\tau,h})\eqqcolon \xi_{\tau,h} \quad \textup{for all }  x\in \R^d.
\end{align*}
If $u(x_0,t_0)=\varphi(x_0)$, then \eqref{eq:newviscsub} holds trivially. If $u(x_0,t_0)>\varphi(x_0)$, then, by locally uniform convergence and the time equicontinuity of the scheme, $U_{\tau,h}(x_{\tau,h},t_{\tau,h}+\tau)>\varphi(x_{\tau,h})$ for $\tau,h$ small. Thus, the numerical scheme~\eqref{eq:numsch} at the point $x_{\tau,h}$ and time $t_{\tau,h}+\tau$ gives
\[
D_\tau^\alpha U_{\tau,h}(x_{\tau,h},t_{\tau,h}+\tau) + \fld U_{\tau,h}(x_{\tau,h},t_{\tau,h})=0.
\]
Now, since $U_{\tau,h}(\cdot,t_{\tau,h})+\xi_{\tau,h}\le\phi$ in $\R^d$ with equality at $x_{\tau,h}$,
\begin{align*}
\fld U_{\tau,h}(x_{\tau,h},t_{\tau,h})&= \frac{1}{h^{2s}}\sum_{\beta \not=0} \big((U_{\tau,h}(x_{\tau,h},t_{\tau,h})+\xi_{\tau,h}) - (U_{\tau,h}(x_{\tau,h}+x_\beta,t_{\tau,h})+\xi_{\tau,h})\big)\omega_\beta\\
&=\frac{1}{h^{2s}}\sum_{\beta \not=0} \big(\phi(x_{\tau,h}) - (U_{\tau,h}(x_{\tau,h}+x_\beta,t_{\tau,h})+\xi_{\tau,h})\big)\omega_\beta\\
&\geq \frac{1}{h^{2s}}\sum_{\beta \not=0} (\phi(x_{\tau,h}) - \phi(x_{\tau,h}+x_\beta))\omega_\beta\\
&= \fld \phi(x_{\tau,h}).
\end{align*}
Hence
\[
D_\tau^\alpha U_{\tau,h}(x_{\tau,h},t_{\tau,h}+\tau) + \fld \phi (x_{\tau,h})\leq 0.
\]
Passing to the limit as $\tau,h\to0^+$, using \Cref{lem:convergence-discrete-Caputo} for the first term (with the family $(x_{\tau,h},t_{\tau,h}+\tau)\to(x_0,t_0)$) and the consistency of the discretization of the fractional Laplacian for the second one, $\fld\phi(x_{\tau,h})\to\fl\phi(x_0)$ (see, e.g., \cite{Ciaurri2018,delTeso2018b}), we conclude that
\[
\partial_t^\alpha u(x_0,t_0) + \fl \phi (x_0)\leq 0,
\]
which is~\eqref{eq:newviscsub}.
\end{proof}

We now have all the ingredients to prove our main result.
\begin{theorem}\label{thm:AsBeParabObs}
    Let $s\in(0,1]$ and $\alpha\in(0,1)$. Let $u_0,\varphi\in \X_{\textup{b}}^a(\R^d)$ for some $a\in[0,2]$ be such that $u_0\geq \varphi$ in~$\R^d$. Let $u$ be the viscosity solution of \eqref{eq:main}--\eqref{eq:initial} with initial datum $u_0$ and obstacle $\varphi$. Then:
    \begin{enumerate}[\rm (a)]
        \item\label{thm:AsBeParabObs-item1} There exists a sequence $\{t_{k}\}_{k\in\N}\subset[0,\infty)$ with $t_k\to+\infty$ as $k\to+\infty$ such that
        \begin{align*}
            u(\cdot,t_k) \to v \quad \text{as }  k\to+\infty  \text{ locally uniformly in }  \R^d,
        \end{align*}
        where $v$ is a viscosity solution of the elliptic obstacle problem \eqref{eq:elliptobstacle}. Moreover, the limit of any locally uniformly convergent sequence $u(\cdot,t_k)$ with $t_k\to+\infty$ is a viscosity solution of \eqref{eq:elliptobstacle}.
        \item\label{thm:AsBeParabObs-item2} If, additionally, $u_0=\varphi$ in $\R^d$, then
        \begin{align*}
            u(\cdot,t) \to \underline{v} \quad \textup{as } t\to+\infty \textup{ locally uniformly in }  \R^d,
        \end{align*}
        where $\underline{v}$ is the minimal viscosity solution of the elliptic obstacle problem \eqref{eq:elliptobstacle}.
    \end{enumerate}
\end{theorem}

\begin{proof}
\emph{\eqref{thm:AsBeParabObs-item1} and \eqref{thm:AsBeParabObs-item2} for $u_0,\varphi\in {\X^{2}_\textup{b}}(\R^d)$.} \ By \Cref{main:ConvAndPropViscositySolutionParabolic}\eqref{main:ConvAndPropViscositySolutionParabolic-itemb} and~\eqref{main:ConvAndPropViscositySolutionParabolic-itemc}, the family $\{u(\cdot,t)\}_{t\ge0}$ is uniformly bounded and equicontinuous, so that, by the Arzelà--Ascoli theorem, every sequence $t_k\to+\infty$ has a subsequence along which $u(\cdot,t_k)$ converges locally uniformly in $\R^d$; this proves the existence part of~\eqref{thm:AsBeParabObs-item1}. Let now $v$ be the locally uniform limit of $u(\cdot,t_k)$ for some $t_k\to+\infty$. Clearly $v\in C_\textup{b}(\R^d)$ and $v\ge\varphi$. To see that $v$ is a viscosity subsolution of~\eqref{eq:elliptobstacle}, let $x_0\in\R^d$ with $v(x_0)>\varphi(x_0)$ and let $\phi$ be a test function touching $v$ from above at $x_0$, strictly in $\overline{B_r(x_0)}$ and with a positive gap outside, as in \Cref{rem:viscbetter}\eqref{item:moregeneraltestb}. By locally uniform convergence, there are points $x_k\to x_0$ such that $\phi-u(\cdot,t_k)$ attains its global minimum at $x_k$, and $u(x_k,t_k)>\varphi(x_k)$ for $k$ large. Since $\phi+\xi_k$, with $\xi_k\coloneqq u(x_k,t_k)-\phi(x_k)$, is a test function touching $u(\cdot,t_k)$ from above at $x_k$ in the sense of \Cref{lem:clasical-in-time-viscosity}, we get from~\eqref{eq:newviscsub}
\[
    \partial_t^\alpha u(x_k,t_k)+\fl\phi(x_k)\le0.
\]
By \Cref{prop:CaputoDefined}\eqref{prop:CaputoDefined-itemb}, $\partial_t^\alpha u(x_k,t_k)\to0$ as $k\to+\infty$, and $\fl\phi(x_k)\to\fl\phi(x_0)$ by continuity. Hence $\fl\phi(x_0)\le0$, as desired. The supersolution property is proved in the same way using~\eqref{eq:newviscsuper}; note that in this case there is no need to distinguish whether $v(x_0)=\varphi(x_0)$ or not, since $u\ge\varphi$ and only the inequality $\partial_t^\alpha u(x_k,t_k)+\fl\phi(x_k)\ge0$ is needed. This proves~\eqref{thm:AsBeParabObs-item1} for smooth data.

As for~\eqref{thm:AsBeParabObs-item2}, assume $u_0=\varphi$, and let $\underline v$ be the minimal viscosity solution of~\eqref{eq:elliptobstacle}. Regarded as a function of $(x,t)$, $\underline v$ is a (stationary) viscosity solution of the evolution obstacle problem~\eqref{eq:main}--\eqref{eq:initial} with initial datum $\underline v$, which belongs to $\X^1_\textup{b}(\R^d)$ by \Cref{lem:infsuper2}\eqref{lem:infsuper2-item3}. Since $u_0=\varphi\le\underline v$, the comparison principle of \Cref{main:ConvAndPropViscositySolutionParabolic}\eqref{main:ConvAndPropViscositySolutionParabolic-itemf} gives $u(\cdot,t)\le\underline v$ for all $t\ge0$. On the other hand, if $v$ is the limit of $u(\cdot,t_k)$ along any locally uniformly convergent sequence, then $v$ is a viscosity solution of~\eqref{eq:elliptobstacle} by~\eqref{thm:AsBeParabObs-item1}, so that $v\ge\underline v$ by minimality. Hence $v=\underline v$. Since every sequence $t_k\to+\infty$ has a subsequence along which $u(\cdot,t_k)$ converges locally uniformly, and the limit is always $\underline v$, we conclude that $u(\cdot,t)\to\underline v$ locally uniformly as $t\to+\infty$.

\noindent\emph{\eqref{thm:AsBeParabObs-item2} in the general case.} \
Let us now assume that $u_0=\varphi \in \X^a_\textup{b}(\R^d)$ for some $a\in(0,2]$, and denote by $\underline{u}$ the minimal viscosity solution of the elliptic obstacle problem with obstacle $\varphi$. For $\delta\in(0,1]$, let $\varphi_\delta\coloneqq\rho_\delta*\varphi-\eta_\delta$, where $\rho_\delta$ is a standard mollifier and $\eta_\delta\coloneqq\|\rho_\delta*\varphi-\varphi\|_{L^\infty(\R^d)}$, so that $\varphi_\delta\in C^\infty_{\textup{b}}(\R^d)$, $\varphi_\delta\leq
\varphi$ in $\R^d$ and $\varphi_\delta\to\varphi$ uniformly in $\R^d$ as
$\delta\to0^+$, since $\eta_\delta\le K\|\varphi\|_{\X^a(\R^d)}\delta^{a\wedge1}\to0$. Let $u_\delta$ be the viscosity solution of the parabolic obstacle
problem with initial datum and obstacle $\varphi_\delta$, and let $\underline{u}_\delta$
be the minimal viscosity solution of the elliptic obstacle problem with obstacle
$\varphi_\delta$. By Step~1,
\[
    u_\delta(\cdot,t) \to \underline{u}_\delta \quad \textup{as } t\to+\infty  \textup{ locally uniformly in } \R^d.
\]
Moreover, by the $L^\infty$-contraction property of \Cref{main:ConvAndPropViscositySolutionParabolic}\eqref{main:ConvAndPropViscositySolutionParabolic-iteme}, for all $t\geq0$ we have
\begin{equation}\label{eq:contraction-delta}
    \|u_\delta(\cdot,t)-u(\cdot,t)\|_{L^\infty(\R^d)}\leq \|\varphi_\delta-\varphi\|_{L^\infty(\R^d)},
\end{equation}
and, since $\varphi_\delta\leq \varphi$, the comparison principle gives $u_\delta \leq u$ in $\R^d\times[0,\infty)$. Finally, as in Step~1, $\underline{u}$ is a stationary solution of the evolution obstacle problem with initial datum $\underline u\in C^{a\wedge1}_\textup{b}(\R^d)$ (\Cref{lem:infsuper2}\eqref{lem:infsuper2-item3}) satisfying $\underline u\ge\varphi=u_0$, whence $u(\cdot,t)\leq \underline{u}$ for all $t\geq0$. Altogether,
\begin{align}\label{eq:keyasbeha}
    \underline{u}_\delta(x) = \lim_{t\to+\infty} {u}_\delta(x,t) \leq \liminf_{t\to+\infty}u(x,t)\leq \limsup_{t\to+\infty}u(x,t) \leq \underline{u}(x)\quad\textup{for all }x\in\R^d.
\end{align}

We now pass to the limit as $\delta\to0^+$. Since $\underline u_\delta$ is the locally uniform limit of $u_\delta(\cdot,t)$, it inherits the bounds of \Cref{main:ConvAndPropViscositySolutionParabolic}\eqref{main:ConvAndPropViscositySolutionParabolic-itemb} and~\eqref{main:ConvAndPropViscositySolutionParabolic-itemc}: for all $\delta\in(0,1]$ and $y\in\R^d$,
\begin{align*}
    \|\underline{u}_\delta\|_{L^\infty(\R^d)} &\leq  \|\varphi_\delta\|_{L^\infty(\R^d)} \leq \|\varphi\|_{L^\infty(\R^d)} +1\quad\textup{for $\delta$ small},\\
    \|\underline{u}_\delta(\cdot+y)-\underline{u}_\delta\|_{L^\infty(\R^d)} &\leq \|\varphi_\delta(\cdot+y)-\varphi_\delta\|_{L^\infty(\R^d)} \leq \|\varphi\|_{\X^a(\R^d)} |y|^{a\wedge1}.
\end{align*}
By Arzelà--Ascoli's theorem, along a subsequence $\delta_j\to0^+$ we have $\underline{u}_{\delta_j}\to v$ locally uniformly in $\R^d$ for some $v\in C_\textup{b}(\R^d)$, and, by \Cref{rem:stability-elliptic} (recall that $\varphi_\delta\to\varphi$ uniformly), $v$ is a viscosity solution of the elliptic obstacle problem with obstacle $\varphi$. Since $\underline{u}$ is the minimal solution, $\underline{u}\leq v$. Hence, from \eqref{eq:keyasbeha},
   \begin{align*}
    v(x) = \lim_{j\to+\infty} \underline{u}_{\delta_j}(x) \leq \liminf_{t\to+\infty}u(x,t)\leq \limsup_{t\to+\infty}u(x,t) \leq \underline{u}(x) \leq v(x)\quad\textup{for all }x\in\R^d,
\end{align*}
so that $v=\underline u$ and $\lim_{t\to+\infty}u(x,t)=\underline{u}(x)$ for all $x\in\R^d$. Since the limit $v$ does not depend on the subsequence, in fact $\underline u_\delta\to\underline u$ locally uniformly as $\delta\to0^+$.

It remains to show that the convergence of $u(\cdot,t)$ is locally uniform. Fix a compact set $K\subset\R^d$ and $\veps>0$. By~\eqref{eq:contraction-delta},
\begin{align*}
\|u(\cdot,t)-\underline{u}\|_{L^\infty(K)}&\leq \|u(\cdot,t)-u_\delta(\cdot,t)\|_{L^\infty(\R^d)}+\|\underline{u}_\delta-\underline{u}\|_{L^\infty(K)}+\|u_\delta(\cdot,t)-\underline{u}_\delta\|_{L^\infty(K)}\\
& \leq \|\varphi-\varphi_\delta \|_{L^\infty(\R^d)} + \|\underline{u}_\delta-\underline{u}\|_{L^\infty(K)}+\|u_\delta(\cdot,t)-\underline{u}_\delta\|_{L^\infty(K)}.
\end{align*}
We first choose $\delta>0$ small enough so that each of the first two terms is smaller than $\veps/3$, and then, for this fixed $\delta$, $T>0$ such that the third term is smaller than $\veps/3$ for all $t\ge T$, which is possible by Step~1. This concludes the proof of~\eqref{thm:AsBeParabObs-item2}.

\noindent\emph{\eqref{thm:AsBeParabObs-item1} in the general case.} \
Let now $u_0,\varphi \in \X^a_\textup{b}(\R^d)$ for some $a\in(0,2]$. The existence of locally uniformly convergent sequences $u(\cdot,t_k)$ with $t_k\to+\infty$ follows, as in Step~1, from the bounds obtained in \Cref{main:ConvAndPropViscositySolutionParabolic}\eqref{main:ConvAndPropViscositySolutionParabolic-itemb} and~\eqref{main:ConvAndPropViscositySolutionParabolic-itemc} and Arzelà--Ascoli's theorem. Let $v$ be the limit of one such sequence; we have to prove that $v$ is a viscosity solution of~\eqref{eq:elliptobstacle}. For $\delta\in(0,1]$, take $\varphi_\delta, \varphi^\delta, u_{0,\delta}, u_{0}^\delta\in C^\infty_{\textup{b}}(\R^d)$ such that
\[
    \varphi_\delta\leq \varphi\leq \varphi^\delta,\qquad u_{0,\delta}\leq u_0\leq u_0^\delta\qquad\textup{in }\R^d,
\]
with $u_{0,\delta}\ge\varphi_\delta$ and $u_0^\delta\ge\varphi^\delta$, and such that $\varphi_\delta,\varphi^\delta\to\varphi$ and $u_{0,\delta}, u_0^\delta\to u_0$ uniformly in $\R^d$ as $\delta\to0^+$ (for instance, $\varphi_\delta=\rho_\delta*\varphi-\eta_\delta$, $\varphi^\delta=\rho_\delta*\varphi+\eta_\delta$, $u_{0,\delta}=\rho_\delta*u_0-\eta_\delta$ and $u_0^\delta=\rho_\delta*u_0+\eta_\delta$, with $\eta_\delta\coloneqq\max\{\|\rho_\delta*\varphi-\varphi\|_{L^\infty(\R^d)},\|\rho_\delta*u_0-u_0\|_{L^\infty(\R^d)}\}\to0$; the use of the same shift $\eta_\delta$ for the data and the obstacles guarantees that $u_{0,\delta}-\varphi_\delta=\rho_\delta*(u_0-\varphi)\ge0$, and similarly $u_0^\delta\ge\varphi^\delta$). Let $u_\delta$ and $u^\delta$ be the viscosity solutions of the evolution obstacle problem with initial data $u_{0,\delta}$ and $u_{0}^\delta$ and obstacles $\varphi_\delta$ and $\varphi^\delta$, respectively. By equicontinuity, passing to a further subsequence of $\{t_k\}$ we may assume that
\begin{align*}
    u_{\delta}(\cdot,t_k) \to v_\delta, \quad u^{\delta}(\cdot,t_k) \to v^\delta  \quad \textup{as }  k\to+\infty \textup{ locally uniformly in } \R^d,
\end{align*}
and by the first step, $v_\delta$ and $v^\delta$ are viscosity solutions of the elliptic obstacle problem with obstacles $\varphi_\delta$ and $\varphi^\delta$, respectively. Moreover, by the comparison principle of \Cref{main:ConvAndPropViscositySolutionParabolic}\eqref{main:ConvAndPropViscositySolutionParabolic-itemf}, $u_\delta\le u\le u^\delta$, whence
\begin{align*}
    v_\delta(x) \leq v(x) \leq  v^\delta(x)\quad\textup{for all }x\in\R^d.
\end{align*}
Again by equicontinuity, now with respect to $\delta$, along a subsequence $\delta_j\to0^+$ we have
\begin{align*}
    v_{\delta_j} \to \underline{w}\textup{ and }  v^{\delta_j} \to \overline{w} \quad\textup{as }  j\to+\infty   \textup{ locally uniformly in } \R^d,
\end{align*}
where $\underline{w}$ and $\overline{w}$ are viscosity solutions of the elliptic obstacle problem with obstacle $\varphi$, by~\Cref{rem:stability-elliptic}. Thus,
\begin{align*}
    \underline{w}(x)\leq v(x) \leq   \overline{w}(x)\quad\textup{for all }x\in\R^d.
\end{align*}

We now show that $\underline{w}=\overline{w}$, which concludes the proof, since then $v=\underline w$ is a viscosity solution of the elliptic obstacle problem. Fix a compact set $K$ and $\veps>0$. Then,
\[
    \|\underline{w}-\overline{w}\|_{L^\infty(K)} \leq \|\underline{w}-v_{\delta_j}\|_{L^\infty(K)}  + \|\overline{w}-v^{\delta_j}\|_{L^\infty(K)} + \|v_{\delta_j}-v^{\delta_j}\|_{L^\infty(K)}.
\]
On the one hand, by locally uniform convergence, we can choose $j$ such that each of the first two terms on the right-hand side is smaller than $\veps/3$. On the other hand, by the $L^\infty$-contraction property, passed to the limit $k\to+\infty$,
\[
    \|v_{\delta_j}-v^{\delta_j}\|_{L^\infty(K)} \leq \max\{\|u_{0,\delta_j}-u_0^{\delta_j}\|_{L^\infty(\R^d)}\coma \|\varphi_{\delta_j}-\varphi^{\delta_j}\|_{L^\infty(\R^d)} \},
\]
which is smaller than $\veps/3$ for $j$ large. Since $K$ and $\veps$ were arbitrary, $\underline w=\overline w$.
\end{proof}

We finally consider initial data with a limit at infinity. Recall that, by \Cref{lem:minimal.with.limit} and \Cref{lem:small-dimension}, when $c_*\ge\limsup_{|x|\to\infty}\varphi(x)$ the minimal viscosity solution $\underline v^{c_*}$ of the elliptic obstacle problem with obstacle $\max\{\varphi,c_*\}$ is the unique solution of~\eqref{eq:elliptobstacle} with limit $c_*$ at infinity if $d>2s$, and the constant $\max\{\sup\varphi,c_*\}$ if $d\le2s$.

\begin{theorem}\label{thm:limit-at-infinity}
    Under the assumptions of \Cref{thm:AsBeParabObs}, assume in addition that the limit
    $$
        \lim_{|x|\to\infty}u_0(x)= c_*
    $$
    exists. Then $u(\cdot,t) \to \underline{v}^{c_*}$ as $t\to+\infty$ locally uniformly in $\R^d$, where $\underline{v}^{c_*}$ is the minimal viscosity solution of the elliptic obstacle problem \eqref{eq:elliptobstacle} with obstacle $\max\{\varphi,c_*\}$. Moreover, if $d>2s$, the convergence is uniform in $\R^d$, whereas if $d\le2s$ and $c_*<\sup\varphi$ it is not (see \Cref{prop:limit-infinity-fixed-time}).
\end{theorem}

\begin{proof}
Note first that $c_*\ge\limsup_{|x|\to\infty}\varphi(x)$, since $u_0\ge\varphi$. We treat separately the cases $d>2s$ and $d\le 2s$.

\noindent\underline{\emph{Case $d>2s$.}} \ Let $\hat{u}$ be the limit of
$u(\cdot,t)$ along an arbitrary sequence $t_k\to+\infty$ for which the limit exists (such sequences exist
by the equicontinuity and equiboundedness of $u$). By
\Cref{thm:AsBeParabObs}\eqref{thm:AsBeParabObs-item1}, $\hat{u}$ is a
viscosity solution of the elliptic obstacle problem~\eqref{eq:elliptobstacle}. We now identify its behaviour at infinity.

\noindent\textit{Lower bound.} Fix $\varepsilon>0$
and let $u_\varepsilon(\cdot,t)=Z(\cdot,t)*\min\{u_0,c_*-\varepsilon\}$ be the mild solution of the subdiffusive heat
equation $\partial_t^\alpha u_\varepsilon + \fl u_\varepsilon = 0$
with initial datum $\min\{u_0, c_*-\varepsilon\}$. Since $u$ is a
supersolution of the time-fractional heat equation with initial datum $
\min\{u_0,c_*-\varepsilon\}\le u_0$, \Cref{lem:mild-viscosity} gives $u_\varepsilon
\le u$ in $\R^d\times[0,\infty)$. On the other hand, since $\min\{u_0,c_*-\varepsilon\}
\to c_*-\varepsilon$ as $|x|\to\infty$, the representation formula through
the fundamental solution~\cite{Kochubei-1990,Eidelman-Kochubei-2004,Kemppainen-Siljander-Zacher-2017}
gives that $u_\varepsilon(\cdot,t)\to c_*-\varepsilon$ uniformly in $\R^d$
as $t\to+\infty$, as in the proof of \Cref{lem:minimal.with.limit}\eqref{lem:minimal.with.limit-item3}. Therefore, there is $T_\varepsilon>0$ such that
\begin{equation}\label{eq:lower-barrier}
    u(x,t)\ge c_*-2\varepsilon\quad\textup{for all }x\in\R^d\textup{ and }t\ge T_\varepsilon,
\end{equation}
and, passing to the limit along $\{t_k\}$, $\hat{u}\ge c_*-2\varepsilon$ in $\R^d$ for every $\varepsilon>0$, that is, $\hat u\ge c_*$ in $\R^d$.

\noindent\textit{Upper bound.} Fix $\varepsilon>0$ and, as in the proof of \Cref{lem:minimal.with.limit}\eqref{lem:minimal.with.limit-item2}, consider the bounded $s$-superharmonic function $U_\varepsilon(x)\coloneqq A_\varepsilon\min\{|x|^{2s-d},1\}+c_*+\varepsilon$, where now $A_\varepsilon$ is chosen so that $U_\varepsilon\ge \max\{u_0,\varphi\}$ in $\R^d$; this is possible because $u_0$ and $\varphi$ are bounded and $\max\{u_0,\varphi\}\le c_*+\varepsilon$ outside a large ball (here we use $d>2s$ to ensure that $|x|^{2s-d}\to0$ as $|x|\to\infty$). Regarded as a function of $(x,t)$, $U_\varepsilon$ is a (stationary) bounded continuous viscosity supersolution of the obstacle problem~\eqref{eq:main}--\eqref{eq:initial}, that is, $U_\varepsilon\in\mathcal V$, so that $u=\inf\limits_{w\in \mathcal{W}} w\le U_\varepsilon$ by \Cref{thm:uniqueness}. In particular, there is $R_\varepsilon>0$ such that 
\begin{equation}\label{eq:upper-barrier}
    u(x,t)\le c_*+2\varepsilon\quad\textup{for all }|x|\ge R_\varepsilon\textup{ and }t\ge0,
\end{equation}
and, passing to the limit along $\{t_k\}$, $\hat{u}\le U_\varepsilon$ in $\R^d$, whence $\limsup_{|x|\to\infty}\hat{u}(x)\le c_*+\varepsilon$ for every $\varepsilon>0$.

Combining both bounds, $\lim_{|x|\to\infty}\hat{u}(x)=c_*$. By \Cref{thm:uniqueness.elliptic}, $\hat{u}=\underline{v}^{c_*}$, independently of the sequence~$\{t_k\}$, which implies that $u(\cdot,t)\to \underline{v}^{c_*}$ locally uniformly in $\R^d$ as $t\to+\infty$.

To see that the convergence is in fact uniform, fix $\varepsilon>0$. By \Cref{lem:minimal.with.limit}\eqref{lem:minimal.with.limit-item2}, enlarging $R_\varepsilon$ if necessary, we have $|\underline v^{c_*}(x)-c_*|\le 2\varepsilon$ for $|x|\ge R_\varepsilon$. Together with~\eqref{eq:lower-barrier} and~\eqref{eq:upper-barrier}, this gives $|u(x,t)-\underline v^{c_*}(x)|\le4\varepsilon$ for all $|x|\ge R_\varepsilon$ and $t\ge T_\varepsilon$, while $\|u(\cdot,t)-\underline v^{c_*}\|_{L^\infty(B_{R_\varepsilon})}\le\varepsilon$ for $t$ large by the locally uniform convergence. Since $\varepsilon$ is arbitrary, the convergence is uniform in $\R^d$.

\noindent\underline{\emph{Case $d\le 2s$.}} \ Let $u^+(\cdot,t)=Z(\cdot,t)*\tilde u_0$ be the mild solution of the time-fractional heat equation with initial datum $\tilde u_0\coloneqq\max\{\sup \varphi,u_0\}$, which is a viscosity solution of the equation by \Cref{lem:mild-viscosity}. Since $u^+\ge\inf\tilde u_0\ge\sup\varphi\ge\varphi$, the function $u^+$ is also a viscosity solution of the obstacle problem with initial datum $\tilde u_0\ge u_0$, and hence $u\le u^+$ by the comparison principle for the obstacle problem, \Cref{main:ConvAndPropViscositySolutionParabolic}\eqref{main:ConvAndPropViscositySolutionParabolic-itemf}. Moreover, $\tilde u_0\to\max\{\sup\varphi,c_*\}$ as $|x|\to\infty$, so that $u^+(\cdot,t)\to\max\{\sup\varphi,c_*\}$ uniformly in $\R^d$ as $t\to+\infty$, again by the argument in the proof of \Cref{lem:minimal.with.limit}\eqref{lem:minimal.with.limit-item3}. On the other hand, arguing as in the case $d>2s$, $u\ge u_\varepsilon$, where $u_\varepsilon$ is the mild solution of the time-fractional heat equation with initial datum $\min\{u_0,c_*-\varepsilon\}$, which converges uniformly to $c_*-\varepsilon$. Passing to the limit along any sequence $t_k\to+\infty$ for which $u(\cdot,t_k)$ converges,
\[
    c_*-\varepsilon\le\lim_{k\to\infty} u(x,t_k)\le \max\{\sup \varphi, c_*\}\quad\textup{for all }x\in\R^d.
\]
We already know from \Cref{thm:AsBeParabObs}\eqref{thm:AsBeParabObs-item1} that the limit is a viscosity solution of the elliptic obstacle problem with obstacle $\varphi$, and, by the argument in the proof of \Cref{lem:small-dimension}, in this range of dimensions such solutions are constants greater than or equal to $\sup\varphi$. Therefore,
\[
    \max\{\sup \varphi,c_*-\varepsilon\}\le\lim_{k\to\infty} u(x,t_k)\le \max\{\sup \varphi, c_*\}.
\]
Letting $\varepsilon \to 0^+$, we get that $\lim_{k\to\infty} u(x,t_k)= \max\{\sup \varphi, c_*\}=\underline v^{c_*}$, where the last identity is \Cref{lem:small-dimension}. Since the limit does not depend on the sequence, and the family $\{u(\cdot,t)\}_{t\ge0}$ is equicontinuous and equibounded, the whole family converges locally uniformly to $\underline v^{c_*}$.
\end{proof}

\begin{remark}\label{rem:discrete-limit-at-infinity}
{\rm \Cref{thm:limit-at-infinity} has a discrete counterpart, with an analogous proof: under~\eqref{as:CFL}, if $u_0\in \X^a_\textup{b}(\R^d)$ has a limit $c_*$ at infinity, the solution of the scheme~\eqref{eq:numsch} converges, as $n\to\infty$, to the minimal solution of the discrete elliptic obstacle problem with obstacle $\max\{\varphi,c_*\}$ (see \Cref{rem:discrete-minimal-with-limit}). The mild solutions of the time-fractional heat equation are replaced by the solutions of the scheme without obstacle, which are given by convolution with its fundamental solution, a nonnegative kernel of unit mass that spreads out as $n\to\infty$, while the barrier for $d>2s$ is built from the fundamental solution $\Phi_h$ of the stationary discrete operator and its decay at infinity, see \Cref{rem:discrete-minimal-with-limit}.
}
\end{remark}

\begin{remark}
   {\rm If $c_*\ge\sup\varphi$, then $\max\{\varphi,c_*\}=c_*$ and $\underline v^{c_*}\equiv c_*$ in every dimension, so that \Cref{thm:limit-at-infinity} says that $u(\cdot,t)\to c_*$ as $t\to+\infty$: in this case the obstacle plays no role in the large-time behaviour, which is the same as for the time-fractional heat equation (for integrable data, a detailed description of the latter can be found in~\cite{CQW-JFA,CQW-VJM}). The obstacle only matters when $c_*<\sup\varphi$, in which case the limit lies strictly above $c_*$ somewhere.}
\end{remark}

When $d\le 2s$ and $c_*<\sup\varphi$, the limit $\underline v^{c_*}\equiv\sup\varphi$ is strictly larger than the limit $c_*$ of the initial datum at infinity. The next result shows that, nevertheless, the solution keeps the value $c_*$ at infinity for every fixed time, so that in this case the convergence in \Cref{thm:limit-at-infinity} cannot be uniform in $\R^d$. The argument does not depend on the dimension; for $d>2s$ it is consistent with the uniform convergence towards $\underline v^{c_*}$, since $\underline v^{c_*}(x)\to c_*$ as $|x|\to\infty$ in that case.

\begin{proposition}\label{prop:limit-infinity-fixed-time}
    Under the assumptions of \Cref{thm:limit-at-infinity}, for every $t>0$ we have
    \[
        \lim_{|x|\to\infty}u(x,t)=c_*.
    \]
    Consequently, if $d\le 2s$ and $c_*<\sup\varphi$, then $\|u(\cdot,t)-\underline v^{c_*}\|_{L^\infty(\R^d)}\ge \sup\varphi-c_*>0$ for every $t>0$, and the convergence $u(\cdot,t)\to\underline v^{c_*}$ as $t\to+\infty$ is not uniform in $\R^d$.
\end{proposition}

\begin{proof}
The lower bound follows from the proof of \Cref{thm:limit-at-infinity}: for every $\varepsilon>0$, 
\[
    u\ge u_\varepsilon=Z(\cdot,t)*\min\{u_0,c_*-\varepsilon\}, 
\]
and, since $Z(\cdot,t)$ is a probability density and $\min\{u_0,c_*-\varepsilon\}\to c_*-\varepsilon$ as $|x|\to\infty$, we have $u_\varepsilon(x,t)\to c_*-\varepsilon$ as $|x|\to\infty$ for every fixed $t>0$, by dominated convergence. Hence $$\liminf\limits_{|x|\to\infty}u(x,t)\ge c_*-\varepsilon$$ for every $\varepsilon>0$.

For the upper bound we use that $u$ is the smallest viscosity supersolution of the obstacle problem, so that it suffices to construct, for each fixed $T>0$, a viscosity supersolution $W$ in $\R^d\times[0,T]$ whose value at time $T$ approaches $c_*$ at infinity. Note that $W$ has to depend on time: when $d\le 2s$ there are no nonconstant bounded $s$-superharmonic functions (\Cref{prop:Liouville.superharmonic}), hence no stationary barrier decaying to $c_*$ at infinity, in contrast with the case $d>2s$, where $U_\varepsilon(x)\coloneqq A_\varepsilon\min\{|x|^{2s-d},1\}+c_*+\varepsilon$ was used for this purpose in the proof of \Cref{thm:limit-at-infinity}. The natural candidates are solutions $Z(\cdot,t)*W_0$ of the time-fractional heat equation with $W_0\ge u_0$, which are supersolutions of the obstacle problem as long as they stay above the obstacle. The difficulty is at infinity: if $W_0\to c_*$, then $Z(\cdot,t)*W_0\to c_*$ at infinity for every $t$, but $\varphi$ may exceed $c_*$ at arbitrarily far away points, so that one cannot guarantee $Z(\cdot,t)*W_0\ge\varphi$ outside a ball.

We overcome this by slightly enlarging the obstacle. Fix $T>0$ and $\delta>0$. Since $\limsup_{|x|\to\infty}\varphi\le c_*$ and $u_0\to c_*$ at infinity, there is $R>0$ such that $\max\{u_0,\varphi\}\le c_*+\delta$ in $\R^d\setminus B_R$. Consider the obstacle problem with the enlarged obstacle $\tilde\varphi\coloneqq\max\{\varphi,c_*+\delta\}$ and initial datum $\tilde u_0\coloneqq\max\{u_0,c_*+\delta\}\ge\tilde\varphi$, and let $\tilde u$ be its viscosity solution. By the comparison principle, \Cref{main:ConvAndPropViscositySolutionParabolic}\eqref{main:ConvAndPropViscositySolutionParabolic-itemf}, $u\le\tilde u$ in $\R^d\times[0,\infty)$. The gain is that $\tilde\varphi\equiv c_*+\delta$ in $\R^d\setminus B_R$: outside $B_R$, staying above the obstacle amounts to staying above a constant, which any solution of the time-fractional heat equation with initial datum $\ge c_*+\delta$ does, since the semigroup preserves the infimum of the datum.

Inside $B_R$ we have $\tilde\varphi\le\max\{\sup\varphi,c_*+\delta\}$, and we ensure that the barrier stays above this value up to time $T$ by adding a large bump to the initial datum. Let $\psi\in C_\textup{c}^\infty(\R^d)$ with $0\le\psi\le1$ and $\psi=1$ in $B_{2R}$, and, for $M>0$ to be chosen, let
\[
    W(x,t)\coloneqq Z(\cdot,t)*\big(\tilde u_0+M\psi\big)(x)=Z(\cdot,t)*\tilde u_0(x)+M\,Z(\cdot,t)*\psi(x).
\]
By \Cref{lem:mild-viscosity}, $W$ is a viscosity solution of the time-fractional heat equation with $W(\cdot,0)=\tilde u_0+M\psi\ge\tilde u_0$. We claim that $W\ge\tilde\varphi$ in $\R^d\times[0,T]$ for $M$ large. On the one hand, as just explained, 
\[
    W\ge Z(\cdot,t)*\tilde u_0\ge\inf\tilde u_0\ge c_*+\delta=\tilde\varphi\quad \textup{in }\R^d\setminus B_R. 
\]
On the other hand, if $|x|\le R$, then $B_R\subset B_{2R}(x)$. Thus, taking advantage of the self-similar form of the fundamental solution, 
\[
    Z(z,t)=t^{-\frac{\alpha d}{2s}}F\big(z\,t^{-\frac{\alpha}{2s}}\big),
\]
where $F>0$ is a radial profile with $\int_{\R^d}F=1$ (see~\cite{Eidelman-Kochubei-2004,Kemppainen-Siljander-Zacher-2017}), we get, 
for all $|x|\le R$ and $t\in[0,T]$,
\begin{align*}
    Z(\cdot,t)*\psi(x)&=\int_{\mathbb{R}^d}Z(y,t)\psi(x-y)\dd y\ge\int_{B_{2R}(x)}Z(y,t)\dd y\ge\int_{B_R}Z(y,t)\dd y\\
    &=\int_{B_{Rt^{-\alpha/(2s)}}}F(\xi)\dd\xi\ge \int_{B_{RT^{-\alpha/(2s)}}}F(\xi)\dd\xi\eqqcolon p_0>0.
\end{align*}
Choosing $M\ge(\sup\varphi-c_*)^+/p_0$ we get 
\[
    W\ge c_*+\delta+Mp_0\ge\max\{\sup\varphi,c_*\}+\delta\ge\max\{\sup\varphi,c_*+\delta\}\ge\tilde\varphi\quad\text{in }B_R\times(0,T], 
\]
which proves the claim. This is the point where the finite horizon is used in an essential way: $p_0$ depends on $T$ and tends to $0$ as $T\to\infty$, so that no choice of $M$ produces a barrier valid for all times, in accordance with the nonexistence of stationary barriers in this range of dimensions.

Therefore, $W$ is a viscosity supersolution of the obstacle problem with obstacle $\tilde\varphi$ in $\R^d\times(0,T]$, with $W(\cdot,0)\ge\tilde u_0=\tilde u(\cdot,0)$. Arguing exactly as in the proof of \Cref{thm:uniqueness}, with $\Omega=\{\tilde u>\tilde\varphi\}\cap(\R^d\times(0,T))$, where $\tilde u$ solves the time-fractional heat equation, and using that $\tilde u=\tilde\varphi\le W$ in the complement of $\Omega$, the finite-horizon version of the comparison principle of \Cref{lem:compErwin} yields $\tilde u\le W$ in $\R^d\times[0,T]$. In particular,
\[
    u(x,T)\le \tilde u(x,T)\le Z(\cdot,T)*\tilde u_0(x)+M\,Z(\cdot,T)*\psi(x)\quad\textup{for all }x\in\R^d.
\]
By dominated convergence, the first term on the right-hand side of the inequality tends to $c_*+\delta$ and the second one tends to $0$ as $|x|\to\infty$, since $Z(\cdot,T)\in L^1(\R^d)$, $\tilde u_0\to c_*+\delta$ at infinity and $\psi$ has compact support. The bump, which was needed to keep the barrier above the obstacle inside $B_R$, is not felt at infinity at any fixed time. Hence $$\limsup_{|x|\to\infty}u(x,T)\le c_*+\delta$$ for every $\delta>0$, and the first assertion follows since $T>0$ is arbitrary.

Finally, if $d\le2s$ and $c_*<\sup\varphi$, then $\underline v^{c_*}\equiv\sup\varphi$ by \Cref{lem:small-dimension}, and the first assertion gives $$\|u(\cdot,t)-\underline v^{c_*}\|_{L^\infty(\R^d)}\ge\lim_{|x|\to\infty}(\sup\varphi-u(x,t))=\sup\varphi-c_*>0$$ for every $t>0$.
\end{proof}

\section{Numerical simulations}
\label{num:section}

We present two sets of numerical experiments for problem~\eqref{eq:main}--\eqref{eq:initial} in dimension $d=1$. Both experiments use the explicit scheme~\eqref{eq:numsch}, with the L1 approximation~\eqref{eq:def.discretization.Caputo} of the Caputo derivative and the spatial discretization of the fractional or the local Laplacian based on powers of the discrete Laplacian~\eqref{eq:discfldef}. The discretization parameters $h$ and $\tau$ satisfy the stability condition~\eqref{as:CFL}. The reported rates are empirical slopes of the corresponding log--log plots.

The first experiment assesses the accuracy of the scheme under mesh refinement, using an explicit stationary solution of the obstacle problem~\eqref{eq:main} as a reference.

The second experiment investigates the large-time behaviour of the numerical solution with initial datum $u_0=\varphi$, by measuring its distance to a stationary profile at increasing times $t_n=n\tau$.

\subsection{An explicit stationary solution}
\label{num:benchmark}

We construct an obstacle for which both the minimal stationary
solution and its contact set can be computed explicitly. Fix
$s\in(0,1/2)$. The construction starts from a nonnegative function $\rho$
supported on $[-1,1]$. Its Riesz potential $u_*=(-\Delta)^{-s}\rho$ (here we need $s\in(0,1/2)$ to satisfy $d>2s$) is $s$-superharmonic
in $[-1,1]$ and $s$-harmonic in $\R\setminus[-1,1]$. We then
choose the obstacle to coincide with $u_*$ on $[-1,1]$
and to lie strictly below it elsewhere; the obstacle is moreover nonnegative and compactly supported, so that the results of Section~\ref{sec:asymptotics} with $c_*=0$ apply.

\begin{lemma}
    Let $s\in(0,1/2)$ and $d=1$. For $K_s=\frac{\Gamma(\tfrac12-s)}{2^{2s}\sqrt{\pi}\,\Gamma(s)2s(1+2s)}$ consider the functions
    \begin{equation}\label{def:explicitu}
        u_*(x)=K_s
     \left(
     |x+1|^{1+2s}-2|x|^{1+2s}+|x-1|^{1+2s}\right), \qquad \rho(x)=(1-|x|)_+.
    \end{equation}
    Then,  $(-\Delta)^s u_*=\rho$ in $\R$. Moreover, $u_*$ is positive, bounded and even, and $ u_*(x)\sim C_s|x|^{-1+2s}$ as $|x|\to\infty$, where $C_s=\Gamma(\frac12-s)/(2^{2s}\sqrt{\pi}\,\Gamma(s))$ is the normalizing constant of the Riesz potential of order~$2s$.
\end{lemma}

\begin{proof}
We will actually prove that $u_*=(-\Delta)^{-s}\rho$. First, note that
\begin{equation}\label{num:potential}
    (-\Delta)^{-s}\rho(x)=C_s\int_{-1}^{1}(1-|y|)|x-y|^{2s-1}\dd y.
\end{equation}
To evaluate this integral explicitly, for each $x$, consider the function $G_x(y)=\frac{|x-y|^{1+2s}}{2s(1+2s)}$. Integrating by parts on $(-1,0)$ and $(0,1)$ gives
\[
    \int_{-1}^{1}(1-|y|)|x-y|^{2s-1}\dd y =G_x(-1)-2G_x(0)+G_x(1).
\]
From here, the formula for $u_*$ follows. It is standard to check that $u_*$ is bounded and even. It is also positive by the positivity of the integrand in \eqref{num:potential}. Finally, since $\int_{-1}^{1}\rho(y)\dd y=1$, the asymptotic decay as $|x|\to\infty$ follows.
\end{proof}

Now, it is standard to construct a suitable obstacle $\varphi_*$ so that $u_*$ is a solution of the corresponding elliptic obstacle problem.
\begin{lemma}
    Let $s\in(0,1/2)$, $d=1$ and $u_*$ be given by \eqref{def:explicitu}. Consider the obstacle
    \begin{equation}\label{num:obstacle}
        \varphi_*(x)=\max\left\{u_*(x)-(x^2-1)_+,\,0\right\}.
    \end{equation}
    Then, $u_*$ is the minimal solution of the elliptic obstacle problem~\eqref{eq:elliptobstacle}, that is,
\[
    \min\{(-\Delta)^s u_*,\,u_*-\varphi_*\}=0\quad\text{in }\mathbb R.
\]
Moreover, the contact set $\{x\in\R:\ u_*(x)=\varphi_*(x)\}$ is exactly $[-1,1]$.
\end{lemma}

\begin{proof}
This obstacle is nonnegative, continuous and compactly supported since $u_*(x)\to0$ as $|x|\to+\infty$. By construction, $\varphi_*(x)=u_*(x)$ for $|x|\le1$, and $\varphi_*(x)<u_*(x)$ for $|x|>1$. On the contact set $[-1,1]$ we have $(-\Delta)^s u_*=\rho\ge0$, while $(-\Delta)^s u_*(x)=\rho(x)=0$ for $|x|>1$.

The function $u_*$ belongs to $\X^{1+2s}_{\textup{b}}(\R)$ (see \Cref{rk:obstacle.example} below). Hence $(-\Delta)^su_*$ is given pointwise by its integral formula and is continuous. Therefore, $u_*$ is a classical, and thus a viscosity, solution of the obstacle problem. Finally, $\varphi_*\ge0$ and $u_*$ tends to zero at infinity, so we can apply \Cref{lem:minimal.with.limit} and \Cref{thm:uniqueness.elliptic} with $c_*=0$: the minimal solution with obstacle $\max\{\varphi_*,0\}=\varphi_*$ has limit $0$ at infinity, and $u_*$ is the unique solution with this property. Therefore, $u_*$ is the minimal viscosity solution.
\end{proof}

\begin{remark}\label{rk:obstacle.example}
    {\rm (i) For every $\alpha\in(0,1)$, the function $u(x,t)=u_*(x)$ solves the obstacle problem~\eqref{eq:main}--\eqref{eq:initial} with obstacle $\varphi_*$ and initial datum $u_0=u_*$. Moreover, when $u_0=\varphi_*$ (which is Lipschitz continuous, so that \Cref{thm:intro-asymptotics} applies), the solution converges locally uniformly to $u_*$ as $t\to\infty$ by \Cref{thm:intro-asymptotics}\eqref{thm:intro-asymptotics-item2}. This provides an exact reference profile for both numerical experiments in this section, together with an exact stationary contact set.
    
    \noindent (ii) The explicit formula~\eqref{def:explicitu} shows that $u_*\in \X^{1+2s}_{\textup{b}}(\R)$ and is smooth away from $x=-1,0,1$. At each of these three points, its derivative has H\"older exponent exactly $2s$; in particular, $u_*$ is not $\X^2_{\textup{b}}$ there. This limited regularity should be taken into account when interpreting the errors observed under mesh refinement.    
    }
\end{remark}

\subsection{Discretization and implementation}
\label{num:implementation}

We use nodes $x_i=ih$ and $t_n=n\tau$. The time step is chosen separately for each $(h,s,\alpha)$:
\begin{equation}\label{num:timestep}
    \tau_{\rm CFL}=\left(\frac{(2-2^{1-\alpha})h^{2s}}{\Lambda\Gamma(2-\alpha)}\right)^{1/\alpha},\qquad N_t=\left\lceil\frac{T}{0.9\tau_{\rm CFL}}\right\rceil,\qquad \tau=\frac{T}{N_t}.
\end{equation}
This satisfies the \eqref{as:CFL} restriction and reaches $T$ exactly.

In one dimension, the weights of the powers-of-the-discrete-Laplacian discretization satisfy, for $k\geq1$ (see, e.g., \cite{Ciaurri2018}),
\begin{equation}\label{num:weights}
    \omega_{-k}=\omega_k =\frac{s\Gamma(2s+1)\Gamma(k-s)}{\Gamma(1-s)\Gamma(1+s)\Gamma(k+1+s)}.
\end{equation}
The weights are evaluated using the following formulas, which are straightforward to verify from the definition and avoid computing ratios of large gamma functions:
\[
    \omega_1= \frac{\Gamma(2s+1)}{\Gamma(s)\Gamma(2+s)}, \quad \omega_{k+1}= \frac{k-s}{k+1+s}\omega_k, \quad \textup{and} \quad \Lambda=2\sum_{k=1}^{\infty}\omega_k=\frac{\Gamma(2s+1)}{\Gamma(1+s)^2}.
\]

Given $V:h\Z \to \R$, the infinite-lattice operator is thus given by
\[
    (-\Delta_h)^s [V]_i=\frac{1}{h^{2s}}\sum_{k\ne0}\omega_{|k|}(V_i-V_{i+k}).
\]
In practice, we need to work with a truncated-sum version $\mathcal{L}_{h,R}$ of $(-\Delta_h)^s$ for a feasible finite implementation. We describe below how this is done.

Writing $b_j=(j+1)^{1-\alpha}-j^{1-\alpha}$, the update at the computed nodes takes the form
\begin{equation}\label{num:update}
    U_i^{n+1}=\max\left\{b_nU_i^0+\sum_{k=1}^{n}(b_{n-k}-b_{n+1-k})U_i^k-\Gamma(2-\alpha)\tau^\alpha\mathcal L_{h,R}[U^n_\cdot]_i\, , \,  \varphi_*(x_i)\right\}.
\end{equation}

Thus, the entire L1 history is retained, without a short-memory approximation.

\paragraph{Exterior values and spatial summation.}
In our experiments we only update the numerical solution according to the scheme~\eqref{num:update} on the finite domain $[-L,L]$. At every lattice node $x_i$ with $|x_i|>L$ we prescribe $U_i^n=u_*(x_i)$ at all time steps $n$. This is important for numerical accuracy, because $u_*$ has an algebraic tail and does not vanish outside a bounded interval. For $K=\lceil R/h\rceil$, with $R>2L+2$, the implemented operator is
\begin{equation}\label{num:finite-operator}
    \mathcal L_{h,R}[V]_i=\frac{1}{h^{2s}}\left(\Lambda V_i-\sum_{0<|k|\leq K}\omega_{|k|}\widetilde V_{i+k}\right),\quad |x_i|\leq L,\qquad\widetilde V_j=
    \begin{cases}
        V_j,& |x_j|\leq L,\\
        u_*(x_j),& |x_j|>L.
    \end{cases}
\end{equation}
Thus, the full diagonal coefficient $\Lambda$ is retained.  The nearby exterior convolution is evaluated by FFTs. The exterior
contribution is independent of time and is computed only once. 

\begin{remark} 
    {\rm The exterior cutoff remains a separate source of error. Since $u_*(x_k)$ decays like $|x_k|^{2s-1}$ and the weights $\omega_k$ like $|x_k|^{-1-2s}$, the omitted exterior contribution decays proportionally to $R^{-1}$. Increasing $R$ is therefore necessary if this error becomes visible under mesh refinement.}
\end{remark}

\subsection{Experiment 1: mesh refinement around a stationary solution} 
\label{num:mesh}

We perform the experiments for every pair $(\alpha,s)$ with $s\in\{1/8,1/4,3/8\}$ and $\alpha\in\{0.5,0.8\}$. We initialize all computed nodes with $U_i^0=u_*(x_i)$ and use
\[
    T=10,\qquad L=10,\qquad R=10^6,\qquad h\in\{1/4,1/8,1/16,1/32\}.
\]
The exact continuous solution is stationary. We measure the error on the observation domain $[-5,5]$ as
\begin{equation}\label{num:mesh-error}
 E(h;s,\alpha)=
 \max_{0\leq n\leq N_t}\ \max_{|x_i|\leq5}
       |U_i^n-u_*(x_i)|.
\end{equation}
This is a maximum over the full discrete space--time observation set,
not merely the error at time $T$. The sampled continuous profile need not be
an equilibrium of the discrete problem, so a small numerical evolution is
expected even with exact initialization and exterior values.

\begin{figure}[h!]
 \centering
 \includegraphics[width=\linewidth]{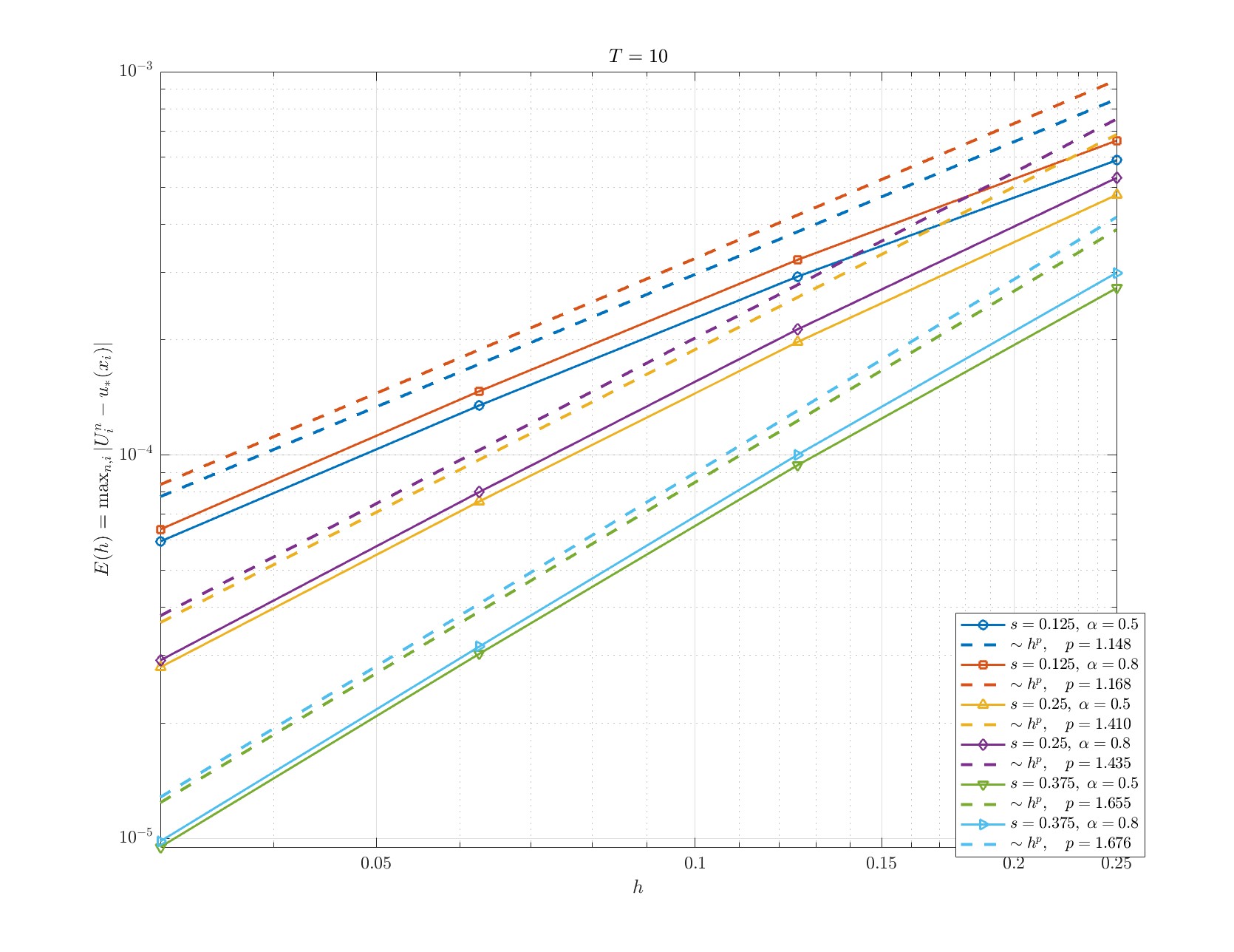}
 \caption{Mesh-refinement experiment: the error~\eqref{num:mesh-error}
 against $h$, with $T=10$. Solid curves are the computed errors. Dashed
 lines are the fitted powers $\sim h^{p}$.}
 \label{num:fig-mesh}
\end{figure}

Figure~\ref{num:fig-mesh} shows a decrease in the error for all six
pairs of parameters $(\alpha,s)$. We fit $\log E=\log C+p\log h$ by least squares using
the three finest meshes, $h=1/8,1/16,1/32$.

The observed order increases with $s$ and varies only moderately with
$\alpha$ in these tests. All slopes are below two. This is compatible
with the limited regularity $\X^{1+2s}_{\textup{b}}$ of the benchmark and does not contradict the
second-order consistency of the spatial operator on sufficiently smooth
functions: for $\X^{1+2s}_{\textup{b}}$ functions the pointwise consistency error of $\fld$ at the three singular points is only of order $h$, so that orders between one and two are to be expected.

\subsection{Experiment 2: large-time convergence towards the stationary profile} 
\label{num:relaxation}

We now fix
\[
    s=1/16,\qquad h=1/128,\qquad T=200,\qquad L=15,\qquad R=10^5,
\]
and take $\alpha\in\{0.5,0.7,0.9\}$. At the computed nodes, the initial condition is $U_i^0=\varphi_*(x_i)$, while the exterior values remain $u_*(x_i)$ at all times. The initial state therefore differs from the stationary profile and leads to nontrivial large-time convergence. We compute 
\begin{equation}\label{num:time-distance}
    D_{h,s,\alpha}(t_n)=\max_{|x_i|\leq5}|U_i^n-u_*(x_i)|.
\end{equation}

\begin{figure}[h!]
    \centering
    \includegraphics[width=\linewidth]{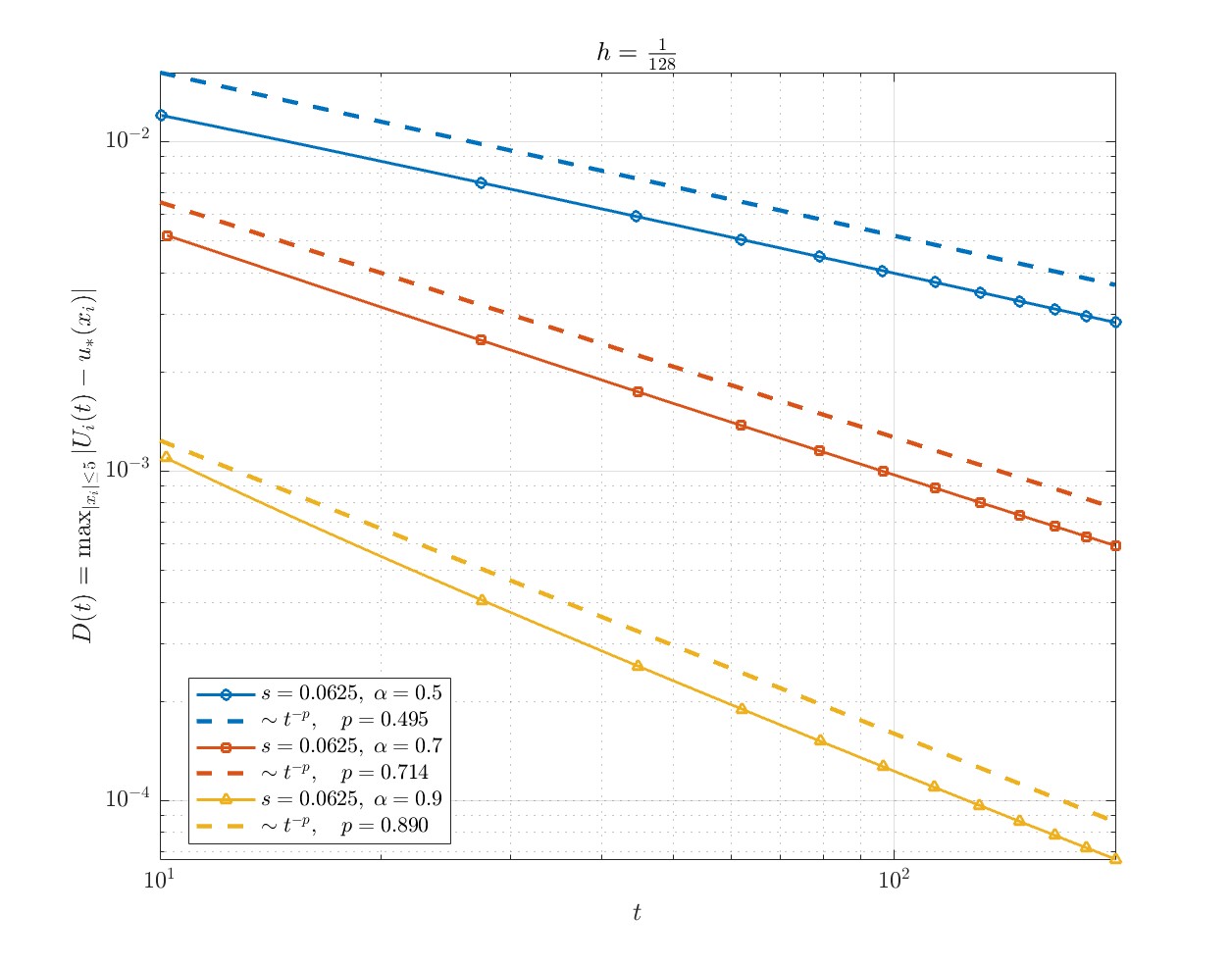}
    \caption{Large-time convergence towards the explicit stationary profile. Solid curves show the
    distance~\eqref{num:time-distance}. Dashed lines show the fitted powers $t^{-p}$.}
    \label{num:fig-time}
\end{figure}

Figure~\ref{num:fig-time} displays this distance for $10\leq t_n\leq200$. We fit $\log D=\log C-p \log t$ n the late-time interval $[80,200]$ using the sampled distances. The distances decrease throughout the displayed interval, with faster decay for larger $\alpha$ in this experiment. The fitted exponents are close to $\alpha$, suggesting a $t^{-\alpha}$ decay over the observed late-time range. 

\section{Complements and open problems}
\label{sec:comments}

This final section contains a complement to the analysis of the numerical scheme, which relates its large-time limits to the continuous elliptic obstacle problem (\Cref{subsec:discrete-approx}), and some open problems (\Cref{subsec:open}).

\subsection{Discrete approximation of the elliptic obstacle problem}\label{subsec:discrete-approx}
As anticipated in \Cref{rem:uniform-in-h}, the large-time limits of the numerical scheme provide approximations of viscosity solutions of the continuous elliptic obstacle problem~\eqref{eq:elliptobstacle}, since the bounds of \Cref{thm:limitscheme2}\eqref{thm:limitscheme2-item1} do not depend on the discretization parameters. This applies in particular to the minimal solutions $\underline U_h$, by \Cref{lem:minimalsoldisc}.

\begin{proposition}\label{prop:discrete-to-continuous}
    Let $s\in(0,1]$ and $u_0,\varphi\in C^0_{\textup{b}}(\R^d)$ be such that $u_0\ge\varphi$ in $\R^d$. For each $h>0$, let $V_h$ be a limit of the scheme~\eqref{eq:numsch} along a time subsequence, as in \Cref{thm:limitscheme2}\eqref{thm:limitscheme2-item1}, with $\tau$ satisfying~\eqref{as:CFL}; for instance, $V_h=\underline U_h$ if $u_0=\varphi$. Then there exist a viscosity solution $v$ of~\eqref{eq:elliptobstacle} and a sequence $h_k\to 0^+$ such that
\[
V_{h_k} \to v \quad \textup{as }  k\to+\infty\textup{ locally uniformly in }   \R^d.
\]
Moreover, $v$ satisfies the same bounds as $V_h$ in \Cref{thm:limitscheme2}\eqref{thm:limitscheme2-item1}.
\end{proposition}
\begin{proof}
By \Cref{thm:limitscheme2}\eqref{thm:limitscheme2-item1}, the family $\{V_h\}_{h>0}$ is equibounded and equicontinuous in $\R^d$, with bounds independent of $h$, so that compactness follows from the Arzelà--Ascoli theorem, and the limit inherits the bounds. If, in addition, $u_0\in \X^a_\textup{b}(\R^d)$ for some $a\in(0,2]$, each $V_h$ is a solution of the discrete elliptic obstacle problem~\eqref{eq:ellipticObs} by \Cref{thm:limitscheme2}\eqref{thm:limitscheme2-item2}, and convergence towards a viscosity solution of~\eqref{eq:elliptobstacle} follows from the locally uniform convergence, the monotonicity of the discrete elliptic operator and the consistency of the discretization of $\fl$, by the classical argument of~\cite{BarlesSouganidis1991}. For general $u_0\in C^0_\textup{b}(\R^d)$, we approximate $u_0$ by $u_{0,\delta}\in C^\infty_\textup{b}(\R^d)$ as in the proof of \Cref{thm:AsBeParabObs}, and use the $L^\infty$-contraction property of the scheme, \Cref{lem:contrac-equicont}\eqref{lem:contrac-equicont-itemcontrac}, passed to the limit in time, together with the stability of viscosity solutions, \Cref{rem:stability-elliptic}.
\end{proof}

\subsection{Open problems}\label{subsec:open}
Our analysis leaves the following questions open.

\begin{enumerate}[\rm (a)]
\item \emph{Rates of convergence.} \Cref{thm:limitscheme2}\eqref{lem:convergence-parabolic-minimal-discrete} and \Cref{thm:AsBeParabObs}\eqref{thm:AsBeParabObs-item2} give convergence towards the minimal solution without a rate. Since the decay of the (discrete) Caputo derivative is quantified in \Cref{prop:limitscheme} and \Cref{prop:CaputoDefined}, a natural approach would be to regard it as a right-hand side in the elliptic obstacle problem and to use some regularity or stability theory for the latter. We further note that, in view of \Cref{thm:intro-asymptotics}, no general result on convergence rates in the uniform norm on $\R^d$ can hold when $d\leq 2s$, since convergence need not even be uniform. 

\item \emph{Convergence of the discrete minimal solutions.} By \Cref{prop:discrete-to-continuous}, every sequence $\underline U_{h_k}$ has a subsequence converging to a viscosity solution $v$ of~\eqref{eq:elliptobstacle}, and $v\ge\underline u$ by minimality of the continuous minimal solution $\underline u$. Therefore, we have that
$$
    \liminf_{h\to 0^+} \underline{U}_h \geq \underline{u}.
$$

In addition, when $\underline{u}$ is constant (this happens for $d\leq 2s$), we get that the constant $\underline{u}$ is a solution to the discrete obstacle problem and hence we have $\underline{u}\geq \underline{U}_h$ in this case (notice that we have not proved that $\underline{U}_h$ is constant when $d\leq 2s$). This implies
$$
    \limsup_{h\to 0^+} \underline{U}_h \leq \underline{u},
$$
and we conclude that
$$
    \lim_{h\to 0^+} \underline{U}_h = \underline{u}.
$$
Whether $\underline U_h\to\underline u$ as $h\to0^+$ when $d> 2s$ is not settled by our arguments.

\item \emph{More general operators and discretizations.} As observed in Section~\ref{subsec:scheme}, the analysis of the scheme only uses that the weights $\omega_\beta$ are nonnegative, symmetric, summable and satisfy~\eqref{eq:weightsbound}. It would be interesting to extend the theory to more general nonlocal operators in space, such as generators of symmetric Lévy processes, and to the corresponding discretizations.
\end{enumerate}

\appendix
\addcontentsline{toc}{section}{Appendix: Comparison and Liouville results for $s$-superharmonic functions}
\section*{Appendix. Comparison and Liouville results for $s$-superharmonic functions}
\setcounter{section}{1}
\setcounter{theorem}{0}
\renewcommand{\thesection}{A}
\label{sec:appendix}

In this appendix we collect the results on $s$-superharmonic functions used in Section~\ref{sec:asymptotics}, namely: (i)~the fact that continuous bounded viscosity $s$-superharmonic functions in an open set are distributionally \hbox{$s$-superharmonic} there (\Cref{prop:equivalence.superharmonic}); (ii) a comparison principle between viscosity sub- and supersolutions, used in the uniqueness result of \Cref{thm:uniqueness.elliptic} when $d>2s$ (\Cref{lem:strong-max-difference}); and (iii) a Liouville theorem for bounded distributional $s$-superharmonic functions in low dimensions, $d\le2s$ (\Cref{prop:Liouville.superharmonic}), which combined with the first result yields the fact, used in \Cref{lem:small-dimension}, that bounded viscosity $s$-superharmonic functions in $\R^d$ are constant when $d\le2s$. None of them is surprising: the first one is a nonlocal counterpart of the classical relation between viscosity and distributional supersolutions of linear elliptic equations; the comparison principle is the strong minimum principle for $s$-superharmonic functions, once the difference of a supersolution and a subsolution is known to be a supersolution; and the Liouville theorem is the analytic expression of the recurrence of the $2s$-stable process (or of Brownian motion, if $s=1$) in dimensions $d\le 2s$, see \Cref{rem:recurrence} below.

For $s=1$ the results are classical. \Cref{prop:equivalence.superharmonic} is \cite[Theorems~1 and~2]{Ishii-1995}, and \Cref{lem:strong-max-difference} reduces to the classical strong minimum principle. As for \Cref{prop:Liouville.superharmonic}, a function $v\in L^\infty(\R^d)$ with $-\Delta v\ge0$ in $\mathcal D'(\R^d)$ coincides almost everywhere with a bounded superharmonic function in the classical sense of potential theory (lower semicontinuous and with the super-mean-value property on balls), see \cite[Chapter~II, Section~4, Propositions~20 and~21]{Dautray-Lions-1990}, and such functions are constant when $d\le2$ \cite[Chapter~5, Theorem~1.20]{Port-Stone-1978}.

For $s\in(0,1)$, the first result follows from the equivalence between weak and viscosity supersolutions of fractional $p$-Laplace type equations proved in \cite[Theorem~1.2]{Korvenpaa-Kuusi-Lindgren-2019}. The comparison principle is then obtained from it by a regularization argument. The Liouville theorem, which we have not found in the literature in the form needed here, relies on the fact that $\int_1^\infty p_t(0)\dd t$ diverges when $d\le 2s$, where $p_t$ is the heat kernel associated with the equation $\partial_tw+(-\Delta)^sw=0$.

We now state and prove these results in turn.

\begin{lemma}\label{prop:equivalence.superharmonic}
    Let $s\in(0,1]$, let $\Omega\subset\R^d$ be open, and let $u\in C_\textup{b}(\R^d)$ be a viscosity supersolution of $\fl u\ge0$ in $\Omega$, in the sense of \Cref{def:visc-elliptic} without obstacle term, the inequality being required only at points $x_0\in\Omega$ and for balls $\overline{B_r(x_0)}\subset\Omega$. Then $\fl u\ge0$ in $\mathcal D'(\Omega)$. The analogous statement holds for subsolutions.
\end{lemma}

\begin{proof}
    For $s=1$ this is \cite[Theorem~1]{Ishii-1995}, which applies to continuous subsolutions of linear elliptic equations in an arbitrary open set, and hence to supersolutions by considering $-u$. 
    
    For $s\in(0,1)$ we use \cite[Theorem~1.2]{Korvenpaa-Kuusi-Lindgren-2019} with $p=2$. Since the distributional inequality is a local property, it is enough to prove it in an arbitrary bounded open set $D\subset\Omega$, which is the setting of that result. Moreover, the notion of viscosity supersolution used in \cite{Korvenpaa-Kuusi-Lindgren-2019}, in which the test function touches $u$ from below in a ball and $u$ itself is used outside, is implied by ours for continuous functions, as in \Cref{rem:stability-elliptic}. Therefore, $u$ is a weak supersolution in $D$, that is, $u\in W^{s,2}_\textup{loc}(D)$ and
    \begin{equation}\label{eq:double.integral}
        \int_{\R^d}\int_{\R^d}(u(x)-u(y))(\varphi(x)-\varphi(y))|x-y|^{-d-2s}\dd x\dd y\ge0
    \end{equation}
    for every nonnegative $\varphi\in C^\infty_\textup{c}(D)$. Finally, for bounded $u$ the double integral in~\eqref{eq:double.integral} is a positive multiple of $\int_{\R^d}u\,\fl\varphi$, as is seen by truncating the domain of integration to $|x-y|>\varepsilon$, using Fubini's theorem and letting $\varepsilon\to0^+$. Hence $\fl u\ge0$ in $\mathcal D'(D)$. Subsolutions are treated by considering~$-u$.
\end{proof}

\begin{proposition}\label{lem:strong-max-difference}
    Let $s\in(0,1]$, let $\Omega\subset\R^d$ be open, and let $u,v\in C_\textup{b}(\R^d)$ be such that $\fl u\ge0$ in $\R^d$ and $\fl v\le0$ in $\Omega$ in the viscosity sense (that is, as in \Cref{def:visc-elliptic} without the obstacle term, the inequality for $v$ being required only at points $x_0\in\Omega$ and for balls $\overline{B_r(x_0)}\subset\Omega$). If $v\le u$ in $\R^d\setminus\Omega$ and $\limsup_{|x|\to\infty}(v-u)(x)\le0$, then $v\le u$ in $\R^d$.
\end{proposition}

\begin{proof}
    By \Cref{prop:equivalence.superharmonic}, applied to $u$ in $\R^d$ and to $v$ in $\Omega$, we have $\fl u\ge0$ in $\mathcal D'(\R^d)$ and $\fl v\le0$ in $\mathcal D'(\Omega)$, so that $w:= u-v$ satisfies $\fl w\ge0$ in $\mathcal D'(\Omega)$. Moreover, $w\ge0$ in $\R^d\setminus\Omega$ and $\liminf_{|x|\to\infty}w\ge0$.

    Assume by contradiction that $m:=\inf_{\R^d}w<0$. By the two properties just mentioned, the infimum is attained at some point of $\Omega$. For $\varepsilon>0$ let $w_\varepsilon:=\rho_\varepsilon*w$, where $\rho_\varepsilon$ is a standard mollifier. Then $w_\varepsilon$ is smooth and bounded, $\fl w_\varepsilon=\rho_\varepsilon*\fl w\ge0$ pointwise in $\Omega_\varepsilon:=\{x:\ B_\varepsilon(x)\subset\Omega\}$, $w_\varepsilon\to w$ locally uniformly, and $w_\varepsilon(x)\ge\inf_{B_\varepsilon(x)}w$. The last inequality implies $\liminf_{|x|\to\infty}w_\varepsilon\ge0$. It also shows that $w_\varepsilon\ge-\omega(2\varepsilon)$ at any point $x$ of a given compact set with $\operatorname{dist}(x,\R^d\setminus\Omega)\le\varepsilon$, where $\omega$ is a modulus of continuity of $w$ on a neighbourhood of that set: indeed, every $y\in B_\varepsilon(x)$ is at distance at most $2\varepsilon$ from a point where $w\ge0$. Since $\inf w_\varepsilon\le m+o(1)<0$ as $\varepsilon\to0^+$, the infimum of $w_\varepsilon$ is attained at some $x_\varepsilon$, which lies in $\Omega_\varepsilon$ for $\varepsilon$ small.

    Fix such an $\varepsilon$. If $s\in(0,1)$, we have $0\le\fl w_\varepsilon(x_\varepsilon)=C_{d,s}\int_{\R^d}\big(w_\varepsilon(x_\varepsilon)-w_\varepsilon(y)\big)|x_\varepsilon-y|^{-d-2s}\dd y\le0$, with a nonpositive integrand, so that $w_\varepsilon$ is constant, equal to $w_\varepsilon(x_\varepsilon)<0$, in the whole $\R^d$, which contradicts that $\liminf_{|x|\to\infty}w_\varepsilon\ge0$. If $s=1$, the smooth superharmonic function $w_\varepsilon$ attains its minimum at the interior point $x_\varepsilon$ of $\Omega_\varepsilon$, so that it is constant in the connected component of $\Omega_\varepsilon$ containing $x_\varepsilon$ by the classical strong minimum principle. This component is the whole $\R^d$, since on its boundary, if nonempty, $w_\varepsilon\ge-\omega(2\varepsilon)>\inf w_\varepsilon$ for $\varepsilon$ small, and we reach the same contradiction.
\end{proof}

\begin{proposition}\label{prop:Liouville.superharmonic}
    Let $s\in(0,1]$ and let $d\leq 2s$. Suppose that $u\in L^\infty(\mathbb R^d)$ satisfies $(-\Delta)^s u\geq 0$ in $\mathcal D'(\mathbb R^d)$. Then $u$ is constant a.e.~in $\mathbb R^d$ (and, in particular, constant everywhere if $u$ is continuous).
\end{proposition}

\begin{proof}
    Set $\mu:=(-\Delta)^s u$. Then $\mu$ is a nonnegative Radon measure. We first prove that $\mu\equiv 0$. Let $\rho_\varepsilon$ be a standard mollifier and define $u_\varepsilon=\rho_\varepsilon*u$. Then
    \[
        u_\varepsilon\in C^\infty(\mathbb R^d)\cap L^\infty(\mathbb R^d), \qquad \|u_\varepsilon\|_{L^\infty}\leq \|u\|_{L^\infty}, \qquad (-\Delta)^s u_\varepsilon = \rho_\varepsilon*((-\Delta)^s u) = \rho_\varepsilon*\mu =:f_\varepsilon\geq 0.
    \]
    It is enough to show that $f_\varepsilon\equiv 0$ for every $\varepsilon>0$. Fix $\varepsilon>0$, and write $v=u_\varepsilon$, $f=f_\varepsilon$. Thus
    \[
        (-\Delta)^s v=f\geq 0, \qquad v\in C^\infty(\mathbb R^d)\cap L^\infty(\mathbb R^d).
    \]
    Let $p_t$ be the heat kernel of the equation $\partial_tw+(-\Delta)^sw=0$, which is the Gaussian kernel if $s=1$ and the fractional heat kernel if $s\in(0,1)$. In both cases $p_t$ is positive, has unit mass, and is self-similar, $p_t(z)=t^{-d/(2s)}p_1(t^{-1/(2s)}z)$. Since $v$ is smooth and bounded with bounded derivatives, $w(x,t):=(p_t*v)(x)$ is a smooth bounded solution of the equation with $w(\cdot,0)=v$, and
    \[
        \partial_t (p_t*v) = -(-\Delta)^s(p_t*v)=-p_t*((-\Delta)^sv)=-p_t*f,
    \]
    so that, for every $T>0$,
    \[
        v(x)-(p_T*v)(x) = \int_0^T (p_t*f)(x)\dd t.
    \]
    If $\|v\|_{L^\infty}\leq M$, then
    \[
        0\leq \int_0^T (p_t*f)(x)\dd t \leq 2M \quad \textup{for every }x\in\mathbb R^d\textup{ and every }T>0.
    \]
    We claim that this is impossible unless $f\equiv0$. Indeed, assume by contradiction that $f\not\equiv0$. Then there exists a compact set $K\subset\mathbb R^d$ such that
    \[
        \int_K f(y)\dd y>0.
    \]
    Fix $x\in\mathbb R^d$. By the self-similarity, positivity and continuity of $p_1$, the heat kernel satisfies
    \[
        p_t(x-y)\geq c\,t^{-d/(2s)} \quad\textup{for all }y\in K\textup{ and all }t\geq t_0,
    \]
    with constants $c>0$ and $t_0>0$ depending on $x$ and $K$. Therefore, for $t\geq t_0$,
    \[
        (p_t*f)(x) = \int_{\mathbb R^d}p_t(x-y)f(y)\dd y \geq c\,t^{-d/(2s)}\int_K f(y)\dd y.
    \]
    Hence
    \[
        \int_0^T (p_t*f)(x)\dd t \geq c\left(\int_K f(y)\dd y\right) \int_{t_0}^T t^{-d/(2s)}\dd t.
    \]
    Since $d\leq 2s$, one has
    \[
        \int_{t_0}^\infty t^{-d/(2s)}\dd t=+\infty.
    \]
    This contradicts the uniform bound
    \[
        \int_0^T (p_t*f)(x)\dd t\leq 2M.
    \]
    Thus $f\equiv0$. Consequently, $\rho_\varepsilon*\mu=0$ for every $\varepsilon>0$, and letting $\varepsilon\to0$ in $\mathcal D'(\mathbb R^d)$ we obtain $\mu=0$.

    Therefore, $u$ is a bounded function with $(-\Delta)^su=0$ in $\mathcal D'(\R^d)$. For $s=1$, by Weyl's lemma $u$ coincides a.e.\ with a harmonic function in the classical sense, which is constant by Liouville's theorem. For $s\in(0,1)$, every function in $L^1_s(\R^d)\coloneqq\{w:\int_{\R^d}|w(x)|(1+|x|)^{-d-2s}\dd x<\infty\}\supset L^\infty(\R^d)$ with $(-\Delta)^sw=0$ in $\mathcal D'(\R^d)$ is affine by \cite[Theorem~1.1]{Fall-PAMS-2016}, and hence constant if it is bounded. In both cases $u$ is constant a.e.~in $\R^d$.
\end{proof}

\begin{remark}[Recurrence, transience, and the critical dimension $d=2s$]\label{rem:recurrence}
    {\rm The proof of \Cref{prop:Liouville.superharmonic} makes transparent why the dimension $d=2s$ is critical throughout the paper. The quantity that decides the argument is
    \[
        \int_{t_0}^\infty p_t(0)\dd t = c\int_{t_0}^\infty t^{-\frac{d}{2s}}\dd t,
    \]
    which is infinite if and only if $d\le 2s$. In probabilistic terms, $p_t$ is the transition density of the symmetric $2s$-stable Lévy process (Brownian motion if $s=1$), and the divergence of $\int_0^\infty p_t(x)\dd t$ is precisely the criterion for the process to be \emph{recurrent}, see \cite[Chapter~7]{Sato-1999}; thus $d\le2s$ is the recurrent regime and $d>2s$ the \emph{transient} one. Recall that a nonnegative function $w$ is \emph{excessive} for the process if $p_t*w\le w$ for all $t>0$; bounded distributional $s$-superharmonic functions correspond, up to an additive constant, to the excessive functions of the process (for $s=1$ this is classical, see \cite[Chapter~5, Theorem~1.19]{Port-Stone-1978}). It is a general fact of potential theory that a recurrent process admits no nonconstant bounded excessive functions, see \cite{Getoor-1980}, and \Cref{prop:Liouville.superharmonic} is the analytic expression of this fact in our setting. For Brownian motion on Riemannian manifolds, the equivalence between recurrence, the divergence of $\int^\infty p_t\dd t$ and this Liouville property is part of the classical characterization of parabolic manifolds, see \cite[Theorem~5.1]{Grigoryan-1999}.

    In the transient regime $d>2s$, on the contrary, $G(x):=\int_0^\infty p_t(x)\dd t=c_{d,s}|x|^{2s-d}$ is finite for $x\neq0$ and is, up to a constant, the fundamental solution of $(-\Delta)^s$. In particular, $\min\{|x|^{2s-d},1\}$ is a bounded, nonconstant $s$-superharmonic function, which shows that \Cref{prop:Liouville.superharmonic} fails for $d>2s$. It is precisely this function the one that we use as a barrier in Section~\ref{sec:asymptotics} to construct solutions of the elliptic obstacle problem with a prescribed limit at infinity and to prove their uniqueness (\Cref{lem:minimal.with.limit} and \Cref{thm:uniqueness.elliptic}), whereas in the recurrent regime $d\le 2s$ such solutions are necessarily constant (\Cref{lem:small-dimension}). The dichotomy in the large-time behaviour of the obstacle problem described in \Cref{thm:limit-at-infinity} is therefore a direct consequence of the recurrence/transience dichotomy for the underlying process.}
\end{remark}

\phantomsection
\addcontentsline{toc}{section}{Acknowledgements}
\section*{Acknowledgements} 

\textsc{F.\,del Teso} is financially supported by MICIU/AEI/10.13039/501100011033 (Spain), through the Ram\'on y Cajal contract RYC2020-029589-I and grants PID2021-127105NB-I00, CNS2024-154515 and RED2024-153842-T. \textsc{F.\,Quir\'os} is financially supported by MICIU/AEI/10.13039/501100011033 (Spain), through grants PID2023-146931NB-I00, RED2024-153842-T, and CEX2023-001347-S. \textsc{J.\,D.\,Rossi} is financially supported by CONICET grants PIP GI No 11220150100036CO, UBACyT grant 20020160100155BA, Argentina, and by MICIU/AEI/10.13039/501100011033, through grants PID2023-146931NB-I00 and CEX2023-001347-S.

\phantomsection
\addcontentsline{toc}{section}{AI disclaimer}
\section*{AI disclaimer}

The authors acknowledge the use of  Claude Opus 5.5 and Claude Sonnet 5.5 (Anthropic), and GPT-6 Astra (OpenAI)  as AI tools for proofreading and for bibliographic verification. The authors independently checked the references against the original sources. All mathematical arguments and proofs in the final manuscript were checked and written by the authors, who retain full responsibility for the whole content.

\phantomsection
\small
\addcontentsline{toc}{section}{References}
\bibliographystyle{abbrv}
\bibliography{bibliography}

\end{document}